\documentclass[10pt]{article}
\usepackage{amsmath}
{
    \allowdisplaybreaks 
}
\usepackage{amsthm}
{
    \theoremstyle{definition}
    {
        \newtheorem{definition}{Definition}[section]

        \newtheorem{condition}[definition]{Condition}
    }
    \theoremstyle{plain}
    {
        \newtheorem{theorem}[definition]{Theorem}
        \newtheorem{proposition}[definition]{Proposition}
        \newtheorem{corollary}[definition]{Corollary}
        \newtheorem{lemma}[definition]{Lemma}
        \newtheorem*{lemma*}{Lemma}
    }
    \theoremstyle{remark}
    {
        \newtheorem{remark}[definition]{Remark}
        
    }
}
\usepackage{amsfonts}
\usepackage{amssymb}
\usepackage{mathtools}
\usepackage{setspace}
\usepackage{tikz-cd}
\usepackage[hidelinks]{hyperref}
\usepackage[style=alphabetic]{biblatex}
{
    \addbibresource{./main.bib}
}
\usepackage{geometry}
{
    \geometry{left=1cm,right=1cm,top=1cm,bottom=1.5cm}

}
\title{GIT for actions of graded unipotent groups on algebraic $\Bbbk$-schemes}

\author{Yikun Qiao\thanks{This work was supported by the Engineering and Physical Sciences Research Council [EP/L015811/1].}}
\date{}
\begin{document}

\setcounter{secnumdepth}{3}
\setcounter{tocdepth}{2}

\maketitle                  
\begin{abstract}
    Let $\Bbbk$ be an algebraically closed field of characteristic zero. We extend the so-called $\hat U$-theorem for projective varieties over $\Bbbk$ proved in \cite{berczi2018constructing} and \cite{berczi2016projective} to projective schemes over $\Bbbk$. The group $\hat U:=U\rtimes \mathbb G_m$ is a semi-direct product of a unipotent group $U$ and a multiplicative group $\mathbb G_m$ acting on $U$ by conjugation, such that the weights of $\mathbb G_m\curvearrowright\mathrm{Lie}(U)$ are all positive. The condition that ``semistability coincides with stability'' (``ss=s'') required in \cite{berczi2018constructing} and \cite{berczi2016projective} is reformulated algebraically for our setting. We also describe a blow-up  procedure to produce a situation where ``ss=s'' holds from one in which it does not. This procedure is different from the sequence of blow-ups in \cite{berczi2016projective}, and produces ``ss=s'' in one blow-up. 
\end{abstract}
\tableofcontents            
\section{Introduction}
{
    Let $\Bbbk$ be an algebraically closed field of characteristic zero. Let $U$ be a unipotent group with Lie algebra $\mathfrak u$. A \emph{grading 1 parameter subgroup} is $\lambda:\mathbb G_m\to \mathrm{Aut}(U)$ such that the weight of $\lambda$ on $\mathfrak u$ are all positive. We call $U$ \emph{graded} if there exists a grading 1PS. Given a grading 1PS $\lambda$, let $\hat U$ denote the semi-direct product $U\rtimes_{\lambda}\mathbb G_m$, which is sometimes denoted by $U\rtimes\lambda$ alternatively. Geometric invariant theory for linearisations of $\hat U$ on projective varieties was studied in \cite{berczi2018geometric}, \cite{berczi2017graded}, \cite{berczi2017gradedlinearisation}, \cite{berczi2016projective}, \cite{hoskins2021quotients}, with conditions on dimensions of $U$-stabilisers or conditions on an extension of the $U$-action to a suitable reductive group. 
    
    Geometric invariant theory for $\hat U$-linearisations has been applied in the construction of moduli of unstable objects. Recently, moduli of pure sheaves of a fixed Harder-Narasimhan type of length 2 on projective varieties are constructed in \cite{jackson2021moduli}. More applications can be found in \cite{berczi2018geometric}, \cite{hamilton2020non}, \cite{hoskins2021quotients}. 

    Let $X$ be a projective variety and let $L$ be an ample line bundle on $X$. Let $\hat U$ act on $X$ linearly with respect to $L$, denoted by $\hat U\curvearrowright(X,L)$. Let $X^0_{\lambda-\min}\subseteq X$ be the non-vanishing locus of $\bigoplus_{d>0}H^0(X,L^d)_{\lambda=\max}$ and let $Z_{\lambda-\min}\hookrightarrow X$ be the vanishing locus of $\bigoplus_{d>0}H^0(X,L^d)_{\lambda<\max}$, where $H^0(X,L^d)_{\lambda=\max}$ is the weight space of $\lambda\curvearrowright H^0(X,L^d)$ whose weight is the maximal, and $H^0(X,L^d)_{\lambda<\max}$ is the sum of weight spaces whose weights are not maximal. In \cite{berczi2018geometric}, under the condition ``ss=s'' (Condition \ref{condition of ss=s}), it is proved that $X^0_{\lambda-\min}\setminus UZ_{\lambda-\min}$ is open and it has a projective geometric quotient by $\hat U$. In \cite{berczi2016projective}, the condition ``ss=s'' is replaced by a weaker version that $\mathrm{Stab}_U(z)=\{e\}$ for all closed points $z\in Z_{\lambda-\min}$. These two versions of ``ss=s'' are actually equivalent for $\hat U$-linearisations on projective varieties (See Corollary \ref{corollary of equivalences of ss=s and stab=e}). 
    
    Let $w_1>\cdots>w_n$ be the positive weights of $\lambda\curvearrowright\mathfrak u$. Let $U_i\subseteq U$ be the subgroup such that $\mathrm{Lie}(U_i)=\mathfrak u_{\lambda\geq w_i}$, which is the sum of weight spaces whose $\lambda$-weights are greater than or equal to $w_i$. As suggested in \cite{berczi2016projective} Remark 2.8 and Remark 7.12, the condition ``ss=s'' can be further weakened to the condition that the minimal dimension of $U_i$-stabilisers on $Z_{\lambda-\min}$ and the minimal dimension of $U_i$-stabilisers on $X^0_{\lambda-\min}$ coincide for all $1\leq i\leq n$, when $X$ is a variety. We split this relaxation into two steps. The first relaxation of ``ss=s'' is the condition that the dimension of $U_i$-stabilisers are constant on $X^0_{\lambda-\min}$ for $1\leq i\leq n$, which is assumed in Section \ref{section of CDRS}. The second relaxation is that the minimal dimensions of $U_i$-stabilisers on $Z_{\lambda-\min}$ and $X^0_{\lambda-\min}$ coincide, which requires a blow-up to satisfy the first relaxation of ``ss=s'' above, which is assumed in Section \ref{section of blowup}. 

    In Section \ref{section of CDRS}, we follow ideas in \cite{berczi2016projective}, Section 7, and generalise the theorem to projective schemes over $\Bbbk$, and complete the proof of the results stated in \cite{berczi2016projective} Remark 2.8 and Remark 7.12, (See Theorem \ref{theorem of universal geometric quotient by Un with proj CDRS}, Theorem \ref{theorem of U-hat quotient with CDRS}). The condition analogous to ``ss=s'' is Condition \ref{condition of CDRS for projective} (the `constant dimension of relative stabilisers' or CDRS condition), which says cokernels of certain \emph{infinitesimal actions} are locally free. More precisely, for an action $G\curvearrowright Y$ over $\Bbbk$, the \emph{infinitesimal action} is a morphism of sheaves $\phi:\Omega_{Y/\Bbbk}\to \mathfrak g^*\otimes_{\Bbbk}\mathcal O_Y$, and for a closed point $y\in Y$, we have $\mathrm{coker}(\phi)\otimes_{\mathcal O_Y}\kappa(y)\cong\mathrm{Stab}_{\mathfrak g}(y)^*$. Therefore the condition that $\dim\mathrm{Stab}_{\mathfrak g}$ is constant can be replaced by the condition that $\mathrm{coker}(\phi)$ is locally free of constant rank. If $X^0_{\lambda-\min}$ is reduced, then Condition \ref{condition of CDRS for projective} is equivalent to the condition that $U_i$-stabilisers have constant dimensions on $X^0_{\lambda-\min}$ (See Lemma \ref{lemma of equivalence of coker being locally free and dim stab is constant when reduced}). 
    
    In Section \ref{section of blowup}, we give a blow-up construction when Condition \ref{condition of CDRS for projective} fails but Condition \ref{condition before blow-up for CDRS} is satisfied. Condition \ref{condition of CDRS for projective} requires certain sheaves to be locally free. Fitting ideals can measure the failure of a sheaf to be locally free in the sense that if $\mathrm{Fit}_{r-1}(\mathcal F)=0$, then $\mathcal F$ is locally free of rank $r$ on the non-vanishing locus of $\mathrm{Fit}_r(\mathcal F)$. They play an important role in determining the centre of the blow-up. In \cite{berczi2016projective}, a sequence of equivariant blow-ups are described to achieve the desired ``ss=s''. Here we provide a different construction, which involves only one blow-up and our argument is more algebraic. Our main theorem (Theorem \ref{theorem of blow-up for CDRS}) says that the result of this blow-up satisfies Condition \ref{condition of CDRS for projective} and thus Theorem \ref{theorem of universal geometric quotient by Un with proj CDRS} and Theorem \ref{theorem of U-hat quotient with CDRS} apply. 
}
\section{Geometric Invariant Theory for \texorpdfstring{$\hat U$}{U-hat}-actions}\label{section of CDRS}
{
    Let $\hat U:=U\rtimes_\lambda\mathbb G_m$ be a semi-direct product of a unipotent group $U$ with a grading 1PS $\lambda:\mathbb G_m\to \mathrm{Aut}(U)$. Let $X$ be a projective scheme over $\Bbbk$ and let $L$ be an ample line bundle on $X$. Let $\hat U$ act linearly on $X$ with respect to $L$. Let $X^0_{\lambda-\min}\subseteq X$ be the non-vanishing locus of $\bigoplus_{d>0}H^0(X,L^d)_{\lambda=\max}$, where $H^0(X,L^d)_{\lambda=\max}\subseteq H^0(X,L^d)$ is the weight space of $\lambda$ with its maximal weight. Without assuming that $X$ is reduced, it is possible that $X^0_{\lambda-\min}=\emptyset$ when sections in $\bigoplus_{d>0}H^0(X,L^d)_{\lambda=\max}$ are all nilpotent.
    
    Let $w_1>\cdots>w_n$ be the weights of $\lambda\curvearrowright\mathfrak u:=\mathrm{Lie}(U)$. The grading 1PS $\lambda$ on $U$ induces a filtration of subgroups $\{U_i\}_{i=0}^n$, where $U_i\subseteq U$ is the subgroup with Lie algebra $\mathfrak u_{\lambda\geq w_i}\subseteq \mathfrak u$, where $\mathfrak u_{\lambda\geq w_i}$ is the sum of weight spaces whose $\lambda$-weights are greater than or equal to $w_i$. The strategy is to construct quotients by $U_i/U_{i-1}$ successively. For this process, we assume Condition \ref{condition of CDRS for projective}, which generalises the relaxed condition ``ss=s'' in \cite{berczi2016projective}, Remark 2.8 to non-reduced situation. When $X$ is reduced, Condition \ref{condition of CDRS for projective} is equivalent to a condition on dimensions of $U_i$-stabilisers (Lemma \ref{lemma of equivalence of coker being locally free and dim stab is constant when reduced}). With this condition, we prove that the universal geometric quotient by $U$ of $X^0_{\lambda-\min}$ exists and the quotient is quasi-projective (Theorem \ref{theorem of universal geometric quotient by Un with proj CDRS}). 

    \subsection{Fitting ideals and stabilisers}
    {
        In this subsection, we introduce some useful exact sequences for group actions (Lemma \ref{Lemma of snake sequence for infinitesimal actions}). We can recover stabilisers of closed point from these sequences (Lemma \ref{lemma of stab as ann and relative stab as ann when flat}). Therefore conditions on unipotent stabilisers can be formulated in terms of these exact sequences. As a result, two versions of the condition ``semistability coincides with stability'' for $\hat U$-linearisations, one in \cite{berczi2018geometric}, another in \cite{berczi2016projective}, are proved equivalent (Corollary \ref{corollary of equivalences of ss=s and stab=e}). 
        
        \subsubsection{Fitting ideals}
        {
            We introduce definitions and properties of Fitting ideals (\cite[\href{https://stacks.math.columbia.edu/tag/07Z6}{Section 07Z6}]{stacks-project}, \cite[\href{https://stacks.math.columbia.edu/tag/0C3C}{Section 0C3C}]{stacks-project}). 
            
            Let $R$ be a commutative ring and let $M$ be a finite $R$-module. There exists a presentation of $M$ 
            \begin{equation}
                \begin{tikzcd}
                    R^{\oplus J}\ar[r]&R^{\oplus n}\ar[r]&M\ar[r]&0
                \end{tikzcd}
            \end{equation}
            where $n\in\mathbb N$ is the number of generators, and $J$ is the cardinality of generators of the kernel of $R^{\oplus n}\to M$. The $R$-linear map $R^{\oplus J}\to R^{\oplus n}$ is represented by a matrix $(a_{ij})_{\substack{1\leq i\leq n\\j\in J}}$ 
            \begin{equation}
                R^{\oplus J}\to R^{\oplus n},\quad e_j\mapsto \sum_{i=1}^na_{ij}e_i. 
            \end{equation}
            
            \begin{definition}
                For $k\in\mathbb Z$, the $k$th Fitting ideal $\mathrm{Fit}_k(M)\subseteq R$ of $M$ is the ideal generated by $(n-k)\times(n-k)$-minors of the the matrix $(a_{ij})_{\substack{1\leq i\leq n\\j\in J}}$. The ideals $\mathrm{Fit}_k(M)$ are independent of the presentation. 
            \end{definition}
    
            Let $S$ be a scheme, and let $\mathcal F$ be a quasi-coherent sheaf of finite type over $S$. 
            
            \begin{definition}
                For $k\in\mathbb Z$, the $k$th Fitting ideal $\mathrm{Fit}_k(\mathcal F)\subseteq \mathcal O_S$ is the quasi-coherent sheaf of ideals such that for any affine open $U\subseteq S$ 
                \begin{equation}
                    \mathrm{Fit}_k(\mathcal F)(U)\cong \mathrm{Fit}_k(\mathcal F(U)). 
                \end{equation}
            \end{definition}
            
            The Fitting ideals form increasing chains 
            \begin{equation}
                \begin{split}
                    0=\mathrm{Fit}_{-1}(M)\subseteq \cdots\subseteq\mathrm{Fit}_{k-1}(M)\subseteq \mathrm{Fit}_k(M)\subseteq\cdots\subseteq R\\
                    0=\mathrm{Fit}_{-1}(\mathcal F)\subseteq \cdots\subseteq\mathrm{Fit}_{k-1}(\mathcal F)\subseteq \mathrm{Fit}_k(\mathcal F)\subseteq\cdots\subseteq \mathcal O_S. 
                \end{split}
            \end{equation}
    
            {
                The following proposition lists some properties of Fitting ideals. 
                \begin{proposition}\label{proposition of properties of Fit}
                    Let $S,T$ be schemes. Let $\mathcal F,\mathcal F'$ be quasi-coherent sheaves of finite type over $S$. Then: 
                    \begin{itemize}
                        \item $\mathcal F$ can be generated by $k$ sections in a neighbourhood of $s\in S$ if and only if $\mathrm{Fit}_k(\mathcal F)_s=\mathcal O_{S,s}$ (\cite[\href{https://stacks.math.columbia.edu/tag/0C3F}{Lemma 0C3F}]{stacks-project}); 
                        \item The closed subset $\mathrm{Supp}(\mathcal F)$ is cut out by $\mathrm{Fit}_0(\mathcal F)$ (\cite[\href{https://stacks.math.columbia.edu/tag/0CYX}{Lemma 0CYX}]{stacks-project}); 
                        \item $\mathcal F$ is locally free of rank $k$ if and only if $\mathrm{Fit}_{k-1}(\mathcal F)=0$ and $\mathrm{Fit}_k(\mathcal F)=\mathcal O_S$ (\cite[\href{https://stacks.math.columbia.edu/tag/0C3G}{Lemma 0C3G}]{stacks-project}); 
                        \item If $f:T\to S$ is a morphism, then $\mathrm{Fit}_k(f^*\mathcal F)=f^{-1}\mathrm{Fit}_k(\mathcal F)\mathcal O_T$ (\cite[\href{https://stacks.math.columbia.edu/tag/0C3D}{Lemma 0C3D}]{stacks-project}); 
                        \item If $\mathcal F\to\mathcal F'$ is surjective, then $\mathrm{Fit}_k(\mathcal F)\subseteq\mathrm{Fit}_k(\mathcal F')$ (\cite[\href{https://stacks.math.columbia.edu/tag/07ZA}{Lemma 07ZA}]{stacks-project}). 
                    \end{itemize}
                \end{proposition}
            }
        }
        
        \subsubsection{Infinitesimal actions}
        {
            Let $G$ be an affine algebraic group with Lie algebra $\mathfrak g$. Let $X$ be a scheme of finite type over $\Bbbk$. Let $G$ act on $X$ with the action morphism $\sigma:G\times_\Bbbk X\to X$. 
            
            Consider the following commutative diagram 
            \begin{equation}
                \begin{tikzcd}
                    X\ar[r,"\gamma"]\ar[d,equal]&G\times_\Bbbk X\ar[d,"\rho"]\\
                    X\ar[r,"\delta"]&X\times_\Bbbk X
                \end{tikzcd},\quad 
                \begin{cases}
                    \gamma:&x\mapsto (e,x)\\
                    \rho:&(g,x)\mapsto \big(\sigma(g,x),x\big)\\
                    \delta:&x\mapsto (x,x). 
                \end{cases}
            \end{equation}
            
            Note that in the diagram above, the morphisms $\gamma,\delta$ are immersions, so by \cite[\href{https://stacks.math.columbia.edu/tag/01R4}{Lemma 01R4}]{stacks-project}, there is a canonical morphism 
            \begin{equation}
                \mathcal C_\delta\to \mathcal C_\gamma
            \end{equation}
            where $\mathcal C_\delta$ and $\mathcal C_\gamma$ are the conormal sheaves of $\delta$ and $\gamma$ respectively. We have that $\mathcal C_\delta\cong\Omega_{X/\Bbbk}$ by \cite[\href{https://stacks.math.columbia.edu/tag/08S2}{Lemma 08S2}]{stacks-project} and $\mathcal C_\gamma\cong \mathfrak g^*\otimes_\Bbbk\mathcal O_X$. Therefore we have a morphism on $X$ 
            \begin{equation}
                \phi:\Omega_{X/\Bbbk}\to \mathfrak g^*\otimes_\Bbbk\mathcal O_X. 
            \end{equation}
            
            {
                \begin{definition}\label{definition of infinitesimal actions}
                    For the action $\sigma:G\times_\Bbbk X\to X$, we call the morphism $\phi:\Omega_{X/\Bbbk}\to \mathfrak g^*\otimes_\Bbbk\mathcal O_X$ the \emph{infinitesimal action of $\mathfrak g$ on $X$}. For a vector subspace $\mathfrak v\subseteq \mathfrak g$, we call the composition $\Omega_{X/\Bbbk}\xrightarrow{\phi}\mathfrak g^*\otimes_\Bbbk\mathcal O_X\to \mathfrak v^*\otimes_\Bbbk \mathcal O_X$ the \emph{infinitesimal action of $\mathfrak v$ on $X$}. 
                \end{definition}
            }

            For a subspace $\mathfrak v\subseteq \mathfrak g$, let $K(\mathfrak v)$ and $Q(\mathfrak v)$ be the kernel and cokernel of the infinitesimal action of $\mathfrak v$ on $X$ 
            \begin{equation}
                \begin{tikzcd}
                    0\ar[r]&K(\mathfrak v)\ar[r]&\Omega_{X/\Bbbk}\ar[r]&\mathfrak v^*\otimes_\Bbbk\mathcal O_X\ar[r]&Q(\mathfrak v)\ar[r]&0. 
                \end{tikzcd}
            \end{equation}
            {
                \begin{lemma}\label{Lemma of snake sequence for infinitesimal actions}
                    Consider the action $\sigma:G\times_\Bbbk X\to X$. Let $\mathfrak w\subseteq\mathfrak v$ be subspaces of $\mathfrak g$. Then there is an exact sequence 
                    \begin{equation}
                        \begin{tikzcd}
                            0\ar[r]&K(\mathfrak v)\ar[r]&K(\mathfrak w)\ar[r]&(\mathfrak v/\mathfrak w)^*\otimes_\Bbbk\mathcal O_X\ar[r]&Q(\mathfrak v)\ar[r]&Q(\mathfrak w)\ar[r]&0. 
                        \end{tikzcd}
                    \end{equation}
                    
                    Moreover, if $Q(\mathfrak w)$ is flat over $X$, then for any scheme $Y$ of finite type over $\Bbbk$ and any morphism $f:Y\to X$, the following is exact 
                    \begin{equation}
                        \begin{tikzcd}
                            f^*K(\mathfrak w)\ar[r]&(\mathfrak v/\mathfrak w)^*\otimes_\Bbbk\mathcal O_Y\ar[r]&f^*Q(\mathfrak v)\ar[r]&f^*Q(\mathfrak w)\ar[r]&0. 
                        \end{tikzcd}
                    \end{equation}
                \end{lemma}
                \begin{proof}
                    Consider the following diagram with splitting exact rows 
                    \begin{equation}
                        \begin{tikzcd}
                            0\ar[r]&0\ar[r]\ar[d]&\Omega_{X/\Bbbk}\ar[r,equal]\ar[d]&\Omega_{X/\Bbbk}\ar[d,"\varphi"]\ar[r]&0\\
                            0\ar[r]&(\mathfrak v/\mathfrak w)^*\otimes_\Bbbk \mathcal O_X\ar[r]&\mathfrak v^*\otimes_\Bbbk \mathcal O_X\ar[r]&\mathfrak w^*\otimes_\Bbbk \mathcal O_X\ar[r]&0. 
                        \end{tikzcd}
                    \end{equation}
                    Apply the Snake lemma to get the long exact sequence in the lemma. Let $\phi_{(\mathfrak v:\mathfrak w)}:K(\mathfrak w)\to(\mathfrak v/\mathfrak w)^*\otimes_\Bbbk\mathcal O_X$ denote the boundary morphism. 
                    
                    Now assume $Q(\mathfrak w)$ is flat. The above commutative diagram pulls back along $f:Y\to X$ to a commutative diagram on $Y$ with exact rows 
                    \begin{equation}
                        \begin{tikzcd}
                            0\ar[r]&0\ar[r]\ar[d]&f^*\Omega_{X/\Bbbk}\ar[r,equal]\ar[d]&f^*\Omega_{X/\Bbbk}\ar[r]\ar[d,"f^*\varphi"]&0\\
                            0\ar[r]&(\mathfrak v/\mathfrak w)^*\otimes_\Bbbk\mathcal O_Y\ar[r]&\mathfrak v^*\otimes_\Bbbk\mathcal O_Y\ar[r]&\mathfrak w^*\otimes_\Bbbk\mathcal O_Y\ar[r]&0. 
                        \end{tikzcd}
                    \end{equation}
                    By the Snake lemma, we have that the following sequence is exact at $f^*Q(\mathfrak v)$ and $f^*Q(\mathfrak w)$
                    \begin{equation}
                        \begin{tikzcd}
                            f^*K(\mathfrak w)\ar[r]&(\mathfrak v/\mathfrak w)^*\otimes_\Bbbk\mathcal O_Y\ar[r]&f^*Q(\mathfrak v)\ar[r]&f^*Q(\mathfrak w)\ar[r]&0. 
                        \end{tikzcd}
                    \end{equation}
                    
                    It remains to show the sequence above is exact at $(\mathfrak v/\mathfrak w)^*\otimes_\Bbbk\mathcal O_Y$. From the first part of the lemma, the following is exact  
                    \begin{equation}
                        \begin{tikzcd}
                            0\ar[r]&\mathrm{coker}(\phi_{(\mathfrak v:\mathfrak w)})\ar[r]&Q(\mathfrak v)\ar[r]&Q(\mathfrak w)\ar[r]&0
                        \end{tikzcd}
                    \end{equation}
                    which pulls back along $f:Y\to X$ to the following exact sequence since $Q(\mathfrak w)$ is flat (\cite[\href{https://stacks.math.columbia.edu/tag/00HL}{Lemma 00HL}]{stacks-project})
                    \begin{equation}
                        \begin{tikzcd}
                            0\ar[r]&f^*\mathrm{coker}(\phi_{(\mathfrak v:\mathfrak w)})\ar[r]&f^*Q(\mathfrak v)\ar[r]&f^*Q(\mathfrak w)\ar[r]&0. 
                        \end{tikzcd}
                    \end{equation}
                    Also, the following is exact since $f^*$ is right exact
                    \begin{equation}
                        \begin{tikzcd}
                            f^*K(\mathfrak w)\ar[r,"f^*\phi_{(\mathfrak v:\mathfrak w)}"]&(\mathfrak v/\mathfrak w)^*\otimes_\Bbbk\mathcal O_Y\ar[r]&f^*\mathrm{coker}(\phi_{(\mathfrak v:\mathfrak w)})\ar[r]&0. 
                        \end{tikzcd}
                    \end{equation}
                    
                    The two short exact sequences above show that the following is exact 
                    \begin{equation}
                        \begin{tikzcd}
                            f^*K(\mathfrak w)\ar[r]&(\mathfrak v/\mathfrak w)^*\otimes_\Bbbk\mathcal O_Y\ar[r]&f^*Q(\mathfrak v)
                        \end{tikzcd}
                    \end{equation}
                    which finishes the proof. 
                \end{proof}
            }
            
            Let $x\in X$ be a closed point. Let $\xi\in\mathfrak g$ be a vector. We say the vector $\xi\in\mathfrak g$ \emph{stabilises} the closed point $x\in X$ if the following composition is zero 
            \begin{equation}
                \begin{tikzcd}
                    \Omega_{X/\Bbbk}\ar[r,"\phi"]&\mathfrak g^*\otimes_\Bbbk\mathcal O_X\ar[r,"\xi\otimes x^\#"]&\Bbbk\otimes_\Bbbk x_*(\Bbbk)\cong x_*(\Bbbk)
                \end{tikzcd}
            \end{equation}
            where we think of $\xi\in\mathfrak g$ as a linear map $\xi:\mathfrak g^*\to \Bbbk$ and we view the closed point $x\in X$ as a morphism $x:\mathrm{Spec}(\Bbbk)\to X$ and $x^\#:\mathcal O_X\to x_*(\Bbbk)$ is the associated morphism of sheaves of rings. 
            
            Let $\mathfrak w\subseteq\mathfrak v$ be subspaces of $\mathfrak g$. For the closed point $x\in X$, its \emph{stabiliser in $\mathfrak v$}, denoted by $\mathrm{Stab}_{\mathfrak v}(x)$, and its \emph{stabiliser in $\mathfrak v$ relative to $\mathfrak w$}, denoted by $\mathrm{Stab}_{(\mathfrak v:\mathfrak w)}(x)$, are defined as follows 
            \begin{equation}
                \begin{split}
                    &\mathrm{Stab}_{\mathfrak v}(x):=\{\xi\in\mathfrak v:\xi\textrm{ stabilises }x\}\subseteq \mathfrak v\\
                    &\mathrm{Stab}_{(\mathfrak v:\mathfrak w)}(x):=\frac{\mathrm{Stab}_{\mathfrak v}(x)}{\mathrm{Stab}_{\mathfrak w}(x)}=\frac{\mathrm{Stab}_{\mathfrak v}(x)}{\mathfrak w\cap\mathrm{Stab}_{\mathfrak v}(x)}\subseteq \mathfrak v/\mathfrak w. 
                \end{split}
            \end{equation}
            
            {
                \begin{lemma}\label{lemma of stab as ann and relative stab as ann when flat}
                    Consider the action $\sigma:G\times_\Bbbk X\to X$. Let $\mathfrak w\subseteq\mathfrak v$ be subspaces of $\mathfrak g$. Let $x:\mathrm{Spec}(\Bbbk)\to X$ be a closed point. Then: 
                    \begin{itemize}
                        \item[(1)] The stabiliser of $x$ in $\mathfrak v$ is 
                        \begin{equation}
                            \mathrm{Stab}_{\mathfrak v}(x)=\mathrm{im}(x^*\phi_{\mathfrak v})^\perp\subseteq \mathfrak v
                        \end{equation}
                        where $\phi_{\mathfrak v}:\Omega_{X/\Bbbk}\to\mathfrak v^*\otimes_\Bbbk\mathcal O_X$ is the infinitesimal action of $\mathfrak v$ on $X$; 
                        \item[(2)] If $Q(\mathfrak w)$ is flat over $X$, then the stabiliser of $x$ in $\mathfrak v$ relative to $\mathfrak w$ is
                        \begin{equation}
                            \mathrm{Stab}_{(\mathfrak v:\mathfrak w)}(x)=\mathrm{im}(x^*\phi_{(\mathfrak v:\mathfrak w)})^\perp\subseteq \mathfrak v/\mathfrak w
                        \end{equation}
                        where $\phi_{(\mathfrak v:\mathfrak w)}:K(\mathfrak w)\to (\mathfrak v/\mathfrak w)^*\otimes_\Bbbk\mathcal O_X$. 
                    \end{itemize}
                \end{lemma}
                \begin{proof}
                    Let $x^\#:\mathcal O_X\to x_*(\Bbbk)$ be the morphism of sheaves of rings associated to $x:\mathrm{Spec}(\Bbbk)\to X$. 
                    
                    For $(1)$, let $\xi\in\mathfrak v$ be a vector, which corresponds to a linear map $\xi:\mathfrak v^*\to\Bbbk$. Recall that $\xi$ stabilises $x$ if and only if the following composition is zero 
                    \begin{equation}
                        \begin{tikzcd}
                            \Omega_{X/\Bbbk}\ar[r,"\phi_{\mathfrak v}"]&\mathfrak v^*\otimes_\Bbbk\mathcal O_X\ar[r,"\xi\otimes x^\#"]&\Bbbk\otimes_\Bbbk x_*(\Bbbk)\cong x_*(\Bbbk)
                        \end{tikzcd}
                    \end{equation}
                    that is $(\xi\otimes x^\#)\circ\phi_{\mathfrak v}=0$ in $\mathrm{Hom}_{\mathcal O_X}(\Omega_{X/\Bbbk},x_*(\Bbbk))$. Since $x_*$ has a left adjoint $x^*$, we have 
                    \begin{equation}
                        \mathrm{Hom}_{\mathcal O_X}(\Omega_{X/\Bbbk},x_*(\Bbbk))\cong \mathrm{Hom}_\Bbbk(x^*\Omega_{X/\Bbbk},\Bbbk).
                    \end{equation}
                    The element $(\xi\otimes x^\#)\circ \phi_{\mathfrak v}\in\mathrm{Hom}_{\mathcal O_X}(\Omega_{X/\Bbbk},x_*(\Bbbk))$ corresponds to the element $\xi\circ (x^*\phi_{\mathfrak v})\in\mathrm{Hom}_\Bbbk(x^*\Omega_{X/\Bbbk},\Bbbk)$, which is the following composition 
                    \begin{equation}
                        \begin{tikzcd}
                            x^*\Omega_{X/\Bbbk}\ar[r,"x^*\phi_{\mathfrak v}"]&\mathfrak v^*\ar[r,"\xi"]&\Bbbk. 
                        \end{tikzcd}
                    \end{equation}
                    Then we have 
                    \begin{equation}
                        \begin{split}
                            \mathrm{Stab}_{\mathfrak v}(x)=&\big\{\xi\in\mathfrak v:(\xi\otimes x^\#)\circ\phi_{\mathfrak v}=0\big\}\\
                            =&\big\{\xi\in\mathfrak v:\xi\circ (x^*\phi_{\mathfrak v})=0\big\}\\
                            =&\mathrm{im}(x^*\phi_{\mathfrak v})^\perp
                        \end{split}
                    \end{equation}
                    which proves $(1)$. 
                    
                    For $(2)$, assume $Q(\mathfrak w)$ is flat and by Lemma \ref{Lemma of snake sequence for infinitesimal actions} the following sequence is exact 
                    \begin{equation}
                        \begin{tikzcd}
                            x^*K(\mathfrak w)\ar[r,"x^*\phi_{(\mathfrak v:\mathfrak w)}"]&(\mathfrak v/\mathfrak w)^*\ar[r]&x^*Q(\mathfrak v)\ar[r]&x^*Q(\mathfrak w)\ar[r]&0. 
                        \end{tikzcd}
                    \end{equation}
                    
                    Note that $(1)$ is equivalent to that the two surjective linear maps $\mathfrak v^*\to x^*Q(\mathfrak v)$ and $\mathfrak v^*\to \big(\mathrm{Stab}_{\mathfrak v}(x)\big)^*$ have the same kernel, and the kernel is $\mathrm{im}(x^*\phi_{\mathfrak v})\subseteq \mathfrak v^*$. In particular $x^*Q(\mathfrak v)\cong \big(\mathrm{Stab}_{\mathfrak v}(x)\big)^*$, and similarly $x^*Q(\mathfrak w)\cong \big(\mathrm{Stab}_{\mathfrak w}(x)\big)^*$. The above exact sequence becomes 
                    \begin{equation}
                        \begin{tikzcd}
                            x^*K(\mathfrak w)\ar[r,"x^*\phi_{(\mathfrak v:\mathfrak w)}"]&(\mathfrak v/\mathfrak w)^*\ar[r]&\big(\mathrm{Stab}_{\mathfrak v}(x)\big)^*\ar[r]&\big(\mathrm{Stab}_{\mathfrak w}(x)\big)^*\ar[r]&0. 
                        \end{tikzcd}
                    \end{equation}
                    Since the kernel of $\big(\mathrm{Stab}_{\mathfrak v}(x)\big)^*\to \big(\mathrm{Stab}_{\mathfrak w}(x)\big)^*$ is $\big(\mathrm{Stab}_{(\mathfrak v:\mathfrak w)}(x)\big)^*$, the following is exact 
                    \begin{equation}
                        \begin{tikzcd}
                            x^*K(\mathfrak w)\ar[r,"x^*\phi_{(\mathfrak v:\mathfrak w)}"]&(\mathfrak v/\mathfrak w)^*\ar[r]&\big(\mathrm{Stab}_{(\mathfrak v:\mathfrak w)}(x)\big)^*\ar[r]&0. 
                        \end{tikzcd}
                    \end{equation}
                    
                    It is easy to see the linear map $(\mathfrak v/\mathfrak w)^*\to\big(\mathrm{Stab}_{(\mathfrak v:\mathfrak w)}(x)\big)^*$ is dual to the inclusion $\mathrm{Stab}_{(\mathfrak v:\mathfrak w)}(x)\subseteq \mathfrak v/\mathfrak w$. Then the above exact sequence is equivalent to 
                    \begin{equation}
                        \mathrm{Stab}_{(\mathfrak v:\mathfrak w)}(x)=\mathrm{im}(x^*\phi_{(\mathfrak v:\mathfrak w)})^\perp
                    \end{equation}
                    which proves $(2)$. 
                \end{proof}
            }
            
            {
                \begin{lemma}\label{lemma of dim of relative stab and Fit when flat}
                    Consider the action $\sigma:G\times_\Bbbk X\to X$. Let $\mathfrak w\subseteq\mathfrak v$ be subspaces of $\mathfrak g$. Let $\phi_{(\mathfrak v:\mathfrak w)}:K(\mathfrak w)\to(\mathfrak v/\mathfrak w)^*\otimes_\Bbbk\mathcal O_X$. Let $x:\mathrm{Spec}(\Bbbk)\to X$ be a closed point and let $k\in\mathbb Z$. Then the following are equivalent: 
                    \begin{itemize}
                        \item[(1)] $\dim \mathrm{coker}(x^*\phi_{(\mathfrak v:\mathfrak w)})>k$; 
                        \item[(2)] $\mathrm{Fit}_k(\phi_{(\mathfrak v:\mathfrak w)})\subseteq \ker x^\#$, where $\mathrm{Fit}_k(\phi_{(\mathfrak v:\mathfrak w)})\subseteq\mathcal O_X$ denotes the $k$th Fitting ideal of $\mathrm{coker}(\phi_{(\mathfrak v:\mathfrak w)})$;
                    \end{itemize}
                    
                    Moreover, when $Q(\mathfrak w)$ is flat, the above are equivalent to the following: 
                    \begin{itemize}
                        \item[(3)] $\dim\mathrm{Stab}_{(\mathfrak v:\mathfrak w)}(x)>k$. 
                    \end{itemize}
                \end{lemma}
                \begin{proof}
                    We have that $(1),(2),(3)$ are all true when $k\leq -1$ and they are all false when $k\geq \dim(\mathfrak v/\mathfrak w)$. Moreover, when $Q(\mathfrak w)$ is flat, by Lemma \ref{lemma of stab as ann and relative stab as ann when flat}, we have $\dim\mathrm{Stab}_{(\mathfrak v:\mathfrak w)}(x)=\dim\mathrm{coker}(x^*\phi_{(\mathfrak v:\mathfrak w)})$, so $(1)$ and $(3)$ are equivalent. Therefore we only need to prove the equivalence of $(1)$ and $(2)$ when $0\leq k<\dim (\mathfrak v/\mathfrak w)$. 
                    
                    Let $r:=\dim(\mathfrak v/\mathfrak w)$. Let $\mathrm{Spec}(A)\subseteq X$ be an affine open subset such that $x:\mathrm{Spec}(\Bbbk)\to X$ factors through $\mathrm{Spec}(A)\subseteq X$. Let $\varphi_x:A\to\Bbbk$ be the ring map corresponding to the closed point $x$ in $\mathrm{Spec}(A)$. Since $K(\mathfrak w)$ is coherent, there is a surjective map $A^n\to K(\mathfrak w)(\mathrm{Spec}(A))$ for some $n\in\mathbb N$. Then the following sequence of $A$-modules is exact 
                    \begin{equation}
                        \begin{tikzcd}
                            A^n\ar[r,"\psi"]&(\mathfrak v/\mathfrak w)^*\otimes_\Bbbk A\ar[r]&\mathrm{coker}(\phi_{(\mathfrak v:\mathfrak w)})(\mathrm{Spec}(A))\ar[r]&0
                        \end{tikzcd}
                    \end{equation}
                    where $\psi$ is the composition $A^n\to K(\mathfrak w)(\mathrm{Spec}(A))\to(\mathfrak v/\mathfrak w)^*\otimes_\Bbbk A$. The base change of the exact sequence above along $\varphi_x:A\to \Bbbk$ is exact 
                    \begin{equation}
                        \begin{tikzcd}
                            \Bbbk^n\ar[r,"\psi\otimes 1_\Bbbk"]&(\mathfrak v/\mathfrak w)^*\ar[r]&\mathrm{coker}(\phi_{(\mathfrak v:\mathfrak w)})(\mathrm{Spec}(A))\otimes_{A,\varphi_x}\Bbbk \ar[r]&0. 
                        \end{tikzcd}
                    \end{equation}
                    The map $\psi$ is between free $A$-modules. Choose a basis of $(\mathfrak v/\mathfrak w)^*$ and let $M=(a_{ij})$ denote the matrix representing $\psi$. Then the matrix representing $\psi\otimes 1_\Bbbk$ is $\overline{M}:=(\overline{a_{ij}})$ for $\overline{a_{ij}}:=\varphi_x(a_{ij})\in\Bbbk$.  

                    Observe that: 
                    \begin{itemize}
                        \item $\dim\mathrm{coker}(x^*\phi_{(\mathfrak v:\mathfrak w)})=\dim\mathrm{coker}(\psi\otimes 1_\Bbbk)$ since 
                        \begin{equation}
                            \begin{split}
                                \mathrm{coker}(x^*\phi_{(\mathfrak v:\mathfrak w)})\cong\;&x^*\mathrm{coker}(\phi_{(\mathfrak v:\mathfrak w)})\\
                                \cong\;&\mathrm{coker}(\phi_{(\mathfrak v:\mathfrak w)})(\mathrm{Spec}(A))\otimes_{A,\varphi_x}\Bbbk\\
                                \cong&\mathrm{coker}(\psi\otimes 1_\Bbbk); 
                            \end{split}
                        \end{equation}
                        
                        \item $\mathrm{Fit}_k(\phi_{(\mathfrak v:\mathfrak w)})(\mathrm{Spec}(A))\subseteq A$ is the ideal generated by $(r-k)\times(r-k)$-minors of $M$; 
                        
                        \item $\mathrm{Fit}_k(\phi_{(\mathfrak v:\mathfrak w)})\subseteq \ker x^\#$ if and only if $\mathrm{Fit}_k(\phi_{(\mathfrak v:\mathfrak w)})(\mathrm{Spec}(A))\subseteq \ker\varphi_x$, since $(\ker x^\#)(\mathrm{Spec}(A))=\ker\varphi_x$ and the support of $x_*(\Bbbk)$ is contained in $\mathrm{Spec}(A)$. 
                    \end{itemize}
                    
                    Then we have:  
                    \begin{equation}
                        \begin{split}
                            \dim\mathrm{coker}(x^*\phi_{(\mathfrak v:\mathfrak w)})>k\quad \textrm{if and only if}\quad& \dim\mathrm{coker}(\psi\otimes 1_\Bbbk)>k\\
                            \textrm{if and only if}\quad&\dim \mathrm{im}(\psi\otimes 1_\Bbbk)<r-k\\
                            \textrm{if and only if}\quad&\mathrm{rank}(\overline{M})<r-k\\
                            \textrm{if and only if}\quad& (r-k)\times(r-k)\textrm{-minors of }(\overline{a_{ij}})\textrm{ vanish}\\
                            \textrm{if and only if}\quad&(r-k)\times(r-k)\textrm{-minors of }(a_{ij}) \textrm{ are in }\ker\varphi_x\\
                            \textrm{if and only if}\quad&\mathrm{Fit}_k(\phi_{(\mathfrak v:\mathfrak w)})(\mathrm{Spec}(A))\subseteq \ker(\varphi_x)\\
                            \textrm{if and only if}\quad&\mathrm{Fit}_k(\phi_{(\mathfrak v:\mathfrak w)})\subseteq \ker x^\#
                        \end{split}
                    \end{equation}
                    which proves the equivalence of $(1)$ and $(2)$. 
                \end{proof}
                
                When $X$ is reduced, we have the following point-wise criterion for flatness of $\mathrm{coker}(\phi_{(\mathfrak v:\mathfrak w)})$. 
                \begin{corollary}\label{corollary of relative stab and coker when reduced and flat}
                    Let $\sigma:G\times_\Bbbk X\to X$, $\mathfrak w\subseteq\mathfrak v$, and $\phi_{(\mathfrak v:\mathfrak w)}$ be as in Lemma \ref{lemma of dim of relative stab and Fit when flat}. Assume $Q(\mathfrak w)$ is flat. Let $k\in\mathbb Z$. Consider the following conditions: 
                    \begin{itemize}
                        \item[(1)] $\mathrm{coker}(\phi_{(\mathfrak v:\mathfrak w)})$ is locally free of rank $k$; 
                        \item[(2)] $\dim\mathrm{Stab}_{(\mathfrak v:\mathfrak w)}(x)=k$ for all closed points $x\in X$. 
                    \end{itemize}
                    Then $(1)\implies (2)$. When $X$ is reduced, they are equivalent. 
                \end{corollary}
                \begin{proof}
                    Recall that for a commutative ring, its \emph{nilpotent radical} is defined as the ideal of all nilpotent elements, equivalently the radical ideal of the zero ideal. Its \emph{Jacobson radical} is defined as the intersection of all maximal ideals. We can easily extend these notions to schemes in an obvious way. Let $\mathrm{nil}(X)\subseteq\mathcal O_X$ denote the nilpotent radical of $X$ and let $\mathrm{rad}(X)\subseteq\mathcal O_X$ denote its Jacobson radical. Recall that $X$ is a scheme of finite type over $\Bbbk$. By \cite[\href{https://stacks.math.columbia.edu/tag/00FV}{Theorem 00FV}]{stacks-project}, we have $\mathrm{rad}(X)=\mathrm{nil}(X)$. The condition that $X$ is reduced is defined as $\mathrm{nil}(X)=0$, which is equivalent to $\mathrm{rad}(X)=0$. 

                    Since $Q(\mathfrak w)$ is flat, we have that $(2)$ is equivalent to the following by Lemma \ref{lemma of dim of relative stab and Fit when flat}: 
                    \begin{itemize}
                        \item[(2')] $\mathrm{Fit}_{k-1}(\phi_{(\mathfrak v:\mathfrak w)})\subseteq \ker x^\#$ and $\mathrm{Fit}_k(\phi_{(\mathfrak v:\mathfrak w)})\not\subseteq\ker x^\#$ for all closed points $x:\mathrm{Spec}(\Bbbk)\to X$
                    \end{itemize}
                    which is equivalent to 
                    \begin{itemize}
                        \item[(2'')] $\mathrm{Fit}_{k-1}(\phi_{(\mathfrak v:\mathfrak w)})\subseteq \mathrm{rad}(X)$ and $\mathrm{Fit}_k(\phi_{(\mathfrak v:\mathfrak w)})=\mathcal O_X$. 
                    \end{itemize}
                    
                    By Proposition \ref{proposition of properties of Fit}, we have that $(1)$ is equivalent to 
                    \begin{itemize}
                        \item[(1')] $\mathrm{Fit}_{k-1}(\phi_{(\mathfrak v:\mathfrak w)})=0$ and $\mathrm{Fit}_k(\phi_{(\mathfrak v:\mathfrak w)})=\mathcal O_X$. 
                    \end{itemize}
                    Obviously $(1')\implies(2'')$, i.e. $(1)\implies (2)$. When $X$ is reduced, we have $\mathrm{rad}(X)=0$ and then $(2'')\iff (1')$, i.e. $(1)\iff (2)$. 
                \end{proof}
            }
            
            {
                In our non-reductive GIT, we often consider a filtration $\{\mathfrak u_i\}_{i=0}^n$ of the Lie algebra $\mathfrak u$ of the unipotent radical, and require the condition (See Condition \ref{condition of CDRS for projective} below) that $\mathrm{coker}(\phi_{(\mathfrak u_i:\mathfrak u_{i-1})})$ is locally free of constant rank for $1\leq i\leq n$. Here is a corollary about this condition. In particular, the condition is equivalent to a condition on unipotent stabilisers of closed points when $X$ is reduced. 
                \begin{corollary}\label{corollary of flatness of coker(phi_i), Q(g_i) and dim of stab when reduced}
                    Consider the action $\sigma:G\times_\Bbbk X\to X$. Let $\mathfrak g$ be the Lie algebra of $G$. Let the following be a filtration of subspaces of $\mathfrak g$ 
                    \begin{equation}
                        0=\mathfrak g_0\subseteq \mathfrak g_1\subseteq \cdots\subseteq \mathfrak g_n=\mathfrak g. 
                    \end{equation}
                    Recall that $K(\mathfrak g_i),Q(\mathfrak g_i)$ are the kernel and cokernel of $\Omega_{X/\Bbbk}\to (\mathfrak g_i)^*\otimes_\Bbbk\mathcal O_X$. For $1\leq i\leq n$, let $\phi_i:K(\mathfrak g_{i-1})\to (\mathfrak g_i/\mathfrak g_{i-1})^*\otimes_\Bbbk\mathcal O_X$ denote the morphism. Let $(k_1,\cdots,k_n)\in\mathbb N^n$ Then the following are equivalent: 
                    \begin{itemize}
                        \item[(1)] $\mathrm{coker}(\phi_i)$ is locally free of rank $k_i$ for each $1\leq i\leq n$; 
                        \item[(2)] $Q(\mathfrak g_i)$ is locally free of rank $k_1+\cdots+k_i$ for each $1\leq i\leq n$. 
                    \end{itemize}
                    
                    The above conditions imply the following: 
                    \begin{itemize}
                        \item[(3)] $\dim\mathrm{Stab}_{\mathfrak g_i}(x)=k_1+\cdots+k_i$ for $1\leq i\leq n$ and all closed points $x\in X$. 
                    \end{itemize}
                    When $X$ is reduced, all three conditions above are equivalent. 
                \end{corollary}
                \begin{proof}
                    Let $\mathfrak v=\mathfrak g_i$ and $\mathfrak w=0$ and apply Corollary \ref{corollary of relative stab and coker when reduced and flat} to claim that $(2)\implies(3)$ in general and $(2)\iff(3)$ when $X$ is reduced. 

                    We then prove the equivalence of $(1)$ and $(2)$. Let $1\leq i\leq n$. By Lemma \ref{Lemma of snake sequence for infinitesimal actions}, the following is exact 
                    \begin{equation}
                        \begin{tikzcd}
                            0\ar[r]&\mathrm{coker}(\phi_i)\ar[r]&Q(\mathfrak g_i)\ar[r]&Q(\mathfrak g_{i-1})\ar[r]&0. 
                        \end{tikzcd}
                    \end{equation}
                    By \cite[\href{https://stacks.math.columbia.edu/tag/05NK}{Lemma 05NK}]{stacks-project}, when $Q(\mathfrak g_{i-1})$ is flat, the flatness of $\mathrm{coker}(\phi_i)$ is equivalent to the flatness of $Q(\mathfrak g_i)$. Initially we have that $Q(\mathfrak g_0)=0$ is flat. Then inductively we can prove that $\mathrm{coker}(\phi_i)$ are flat for all $1\leq i\leq n$ if and only if $Q(\mathfrak g_i)$ are flat for all $1\leq i\leq n$. 
                    
                    By \cite[\href{https://stacks.math.columbia.edu/tag/05P2}{Lemma 05P2}]{stacks-project}, finite locally free coherent sheaves are equivalent to flat coherent sheaves. Ranks of $\mathrm{coker}(\phi_i)$ and $Q(\mathfrak g_i)$ can be read from the short exact sequences for $1\leq i\leq n$. 
                \end{proof}
            }

        }

        \subsubsection{``Semistability coincides with stability''}
        {
            Let $\lambda:\mathbb G_m\to\mathrm{Aut}(U)$ be a grading 1PS with weights $w_1>\cdots>w_n$. Let $\mathfrak u$ be the Lie algebra of $U$ and let $\mathfrak u_{\lambda=w_i}\subseteq\mathfrak u$ be the weight space whose weight is $w_i$. Let $X$ be a projective scheme over $\Bbbk$ and let $L$ be an ample line bundle on $X$. Let $\hat U$ act on $X$ linearly with respect to $L$. Let $X^0_{\lambda-\min}\subseteq X$ be the non-vanishing locus of $\bigoplus_{d>0}H^0(X,L^d)_{\lambda=\max}$ and let $Z_{\lambda-\min}\hookrightarrow X$ be the vanishing locus of $\bigoplus_{d>0}H^0(X,L^d)_{\lambda<\max}$. Then $Z_{\lambda-\min}\hookrightarrow X^0_{\lambda-\min}$ is the \emph{fixed point subscheme} with respect to $\lambda\curvearrowright X^0_{\lambda-\min}$ in the sense of \cite{fogarty1973fixed}. 
            
            For the $\hat U$-linearisation on $(X,L)$, we will call the following condition ``semistability coincides with stability'' (``ss=s''). 
            
            {
                \begin{condition}[cf. Definition 2.1, \cite{berczi2018geometric}]\label{condition of ss=s}
                    There exists a positive integer $n\in\mathbb N_+$ such that $H^0(X,L^n)$ generates $\bigoplus_{d\geq0}H^0(X,L^{nd})$ and the representation of $\hat U$ on $H^0(X,L^n)$ satisfies that if $\xi,\mathfrak u'$ satisfy the following conditions: 
                    \begin{itemize}
                        \item[(1)] $\mathfrak u'\subseteq\mathfrak u$ is a $\lambda$-invariant Lie ideal; 
                        \item[(2)] $\xi\in\mathfrak u\setminus\mathfrak u'$ is a weight vector of $\lambda$; 
                        \item[(3)]  the weights of $\lambda$ on $\mathfrak u'$ are all greater than or equal to the weight of $\xi$. 
                    \end{itemize}
                    then 
                    \begin{equation}
                        H^0(X,L^n)_{\lambda=\max}\subseteq \xi.H^0(X,L^n)^{\mathfrak u'}. 
                    \end{equation}
                \end{condition}
                \begin{remark}
                    Definition 2.1 in \cite{berczi2018geometric} says that if $\xi,\mathfrak u'$ satisfy $(1)$ and $(2)$ above, then\newline $H^0(X,L^n)_{\lambda=\max}\subseteq\xi.H^0(X,L^n)^{\mathfrak u'}$. Therefore Definition 2.1 in \cite{berczi2018geometric} implies Condition \ref{condition of ss=s}. However, they both suffice to prove Theorem 2, the main theorem in \cite{berczi2018geometric}. 
                \end{remark}
            }
            
            For simplicity we often ignore $\lambda$ and write $X^0_{\min}$, $Z_{\min}$, etc. Consider the diagram 
            \begin{equation}
                \begin{tikzcd}
                    X^0_{\min}\ar[r,"\gamma"]\ar[rd,"\delta"]&U\times X^0_{\min}\ar[d,"\sigma"]\\
                    &X^0_{\min}\times X^0_{\min}. 
                \end{tikzcd},\quad \begin{cases}\gamma:x\mapsto (e,x)\\\sigma:(u,x)\mapsto (u.x,x)\\\delta:x\mapsto (x,x). \end{cases}
            \end{equation}
            The following is exact by \cite[\href{https://stacks.math.columbia.edu/tag/01UX}{Lemma 01UX}]{stacks-project}
            \begin{equation}
                \begin{tikzcd}
                    \sigma^*\Omega_{X^0_{\min}\times X^0_{\min}/\Bbbk}\ar[r]&\Omega_{U\times X^0_{\min}/\Bbbk}\ar[r]&\Omega_\sigma\ar[r]&0
                \end{tikzcd}
            \end{equation}
            which pulls back along $\gamma$ to the following exact sequence 
            \begin{equation}
                \begin{tikzcd}
                    \delta^*\Omega_{X^0_{\min}\times X^0_{\min}/\Bbbk}\ar[r]&\gamma^*\Omega_{U\times X^0_{\min}/\Bbbk}\ar[r]&\gamma^*\Omega_\sigma\ar[r]&0
                \end{tikzcd}
            \end{equation}
    
            The morphisms $\gamma,\delta$ are closed immersions with left inverses locally. By \cite[\href{https://stacks.math.columbia.edu/tag/0474}{Lemma 0474}]{stacks-project}, the following rows are exact 
            \begin{equation}
                \begin{tikzcd}
                    0\ar[r]&\mathcal C_\delta\ar[r]\ar[d]&\delta^*\Omega_{X^0_{\min}\times X^0_{\min}/\Bbbk}\ar[r]\ar[d]&\Omega_{X^0_{\min}/\Bbbk}\ar[r]\ar[d,equal]&0\\
                    0\ar[r]&\mathcal C_\gamma\ar[r]&\gamma^*\Omega_{U\times X^0_{\min}/\Bbbk}\ar[r]&\Omega_{X^0_{\min}/\Bbbk}\ar[r]&0. 
                \end{tikzcd}
            \end{equation}
            By the Snake lemma, we have an exact sequence
            \begin{equation}
                \begin{tikzcd}
                    \mathcal C_\delta\ar[r]&\mathcal C_\gamma\ar[r]&\gamma^*\Omega_\sigma\ar[r]&0. 
                \end{tikzcd}
            \end{equation}
            Since $\mathcal C_\delta\cong \Omega_{X^0_{\min}/\Bbbk}$ and $\mathcal C_\gamma\cong \mathfrak u^*\otimes\mathcal O_{X^0_{\min}}$, the above exact sequence is 
            \begin{equation}
                \begin{tikzcd}
                    \Omega_{X^0_{\min}/\Bbbk}\ar[r]&\mathfrak u^*\otimes\mathcal O_{X^0_{\min}/\Bbbk}\ar[r]&\gamma^*\Omega_\sigma\ar[r]&0. 
                \end{tikzcd}
            \end{equation}
    
            {
                \begin{proposition}\label{proposition of equivalence of ss=s and pull-back of differential vanishing}
                    Let $X$ be a projective scheme over $\Bbbk$ and let $L$ be an ample line bundle on $X$. Let $\hat U\curvearrowright (X,L)$. The following are equivalent: 
                    \begin{itemize}
                        \item[(1)] Condition \ref{condition of ss=s} (``ss=s''); 
                        \item[(2)] $\gamma^*\Omega_\sigma=0$. 
                    \end{itemize}
                \end{proposition}
                \begin{proof}
                    Let $m\in\mathbb N_+$ be an integer such that $H^0(X,L^m)$ generates $\bigoplus_{d\geq0}H^0(X,L^{md})$. Then $X^0_{\min}$ is the union of the open affine subschemes $X_f$ for $f\in H^0(X,L^m)_{\max}$. The following is exact 
                    \begin{equation}
                        \begin{tikzcd}
                            \Omega_{X^0_{\min}/\Bbbk}\ar[r]&\mathfrak u^*\otimes\mathcal O_{X^0_{\min}/\Bbbk}\ar[r]&\gamma^*\Omega_\sigma\ar[r]&0
                        \end{tikzcd}
                    \end{equation}
                    whose local version is 
                    \begin{equation}
                        \begin{tikzcd}
                            \Omega_{\mathcal O(X_f)/\Bbbk}\ar[r]&\mathfrak u^*\otimes\mathcal O(X_f)\ar[r]&\gamma^*\Omega_\sigma(X_f)\ar[r]&0. 
                        \end{tikzcd}
                    \end{equation}
                    So the condition $\gamma^*\Omega_\sigma=0$ is equivalent to the surjectivity of $\Omega_{\mathcal O(X_f)/\Bbbk}\to\mathfrak u^*\otimes\mathcal O(X_f)$ for all $f\in H^0(X,L^m)_{\max}$. 
        
                    For $(1)\implies (2)$, assume Condition \ref{condition of ss=s} (``ss=s''). Let $m\in\mathbb N_+$ be such that $H^0(X,L^m)$ generates $\bigoplus_{d\geq0}H^0(X,L^{md})$ and $H^0(X,L^m)_{\max}\subseteq \xi.H^0(X,L^m)^{\mathfrak u'}$ for all $\xi,\mathfrak u'$ satisfying $(1)-(3)$ in Condition \ref{condition of ss=s}. We want to show the following map is surjective 
                    \begin{equation}
                        \begin{split}
                            \phi:\Omega_{\mathcal O(X_f)/\Bbbk}&\to \mathfrak u^*\otimes\mathcal O(X_f)\\
                            a\mathrm{d}b&\mapsto \sum_{i,j}u^{(i)}_j\otimes a\big(\xi^{(i)}_j.b\big)
                        \end{split}
                    \end{equation}
                    where $\big(\xi^{(i)}_1,\cdots,\xi^{(i)}_{r_i}\big)$ is a basis of $\mathfrak u_{w_i}$, the weight space of $\mathfrak u$ whose weight is $w_i$, and $\big(u^{(i)}_1,\cdots,u^{(i)}_{r_i}\big)$ is the set of basis of $(\mathfrak u_{w_i})^*$ dual to $\big(\xi^{(i)}_1,\cdots,\xi^{(i)}_{r_i}\big)$. 
                    
                    Take $\xi:=\xi^{(i)}_j$ and 
                    \begin{equation}
                        \mathfrak u':=\mathrm{span}\big\{\xi^{(i)}_{j'}:1\leq j'\leq r_i,\;j'\ne j\big\}+\bigoplus_{i'<i}\mathfrak u_{w_{i'}}. 
                    \end{equation}
                    Then $\xi,\mathfrak u'$ satisfy $(1)-(3)$ in Condition \ref{condition of ss=s}, so $H^0(X,L^m)_{\max}\subseteq\xi.H^0(X,L^m)^{\mathfrak u'}$. In particular for $f\in H^0(X,L^m)_{\max}$, there exists $g\in H^0(X,L^m)$, which is fixed by $\mathfrak u'$, and $\xi^{(i)}_j.g=f$. Without loss of generality assume $g\in H^0(X,L^m)_{\max-w_i}$, since we can decompose $g$ according to weights and replace $g$ by its eigen-component in $H^0(X,L^m)_{\max-w_i}$. Then $\frac{g}{f}\in\mathcal O(X_f)_{-w_i}$ and 
                    \begin{equation}
                        \begin{split}
                            \mathrm{d}\big(\frac{g}{f}\big)\mapsto&\sum_{i',j'}u^{(i')}_{j'}\otimes \frac{\xi^{(i')}_{j'}.g}{f}\\
                            =&u^{(i)}_j\otimes\frac{\xi^{(i)}_j.g}{f}+\sum_{i'>i}\sum_{j'=1}^{r_{i'}}u^{(i')}_{j'}\otimes \frac{\xi^{(i')}_{j'}.g}{f}\\
                            =&u^{(i)}_j\otimes 1+\sum_{i'>i}\sum_{j'=1}^{r_{i'}}u^{(i')}_{j'}\otimes \frac{\xi^{(i')}_{j'}.g}{f}. 
                        \end{split}
                    \end{equation}
                    The first equality is because that $\xi^{(i')}_{j'}.g=0$ when $i'<i$ or $\begin{cases}i'=i\\j'\ne j\end{cases}$, and the second is because $\xi^{(i)}_j.g=f$. The above for $1\leq j\leq r_i$ together show that 
                    \begin{equation}
                        (\mathfrak u_{w_i})^*\otimes\mathcal O(X_f)\subseteq\mathrm{im}(\phi)+\sum_{i'>i}(\mathfrak u_{w_{i'}})^*\otimes\mathcal O(X_f). 
                    \end{equation}
                    In particular $(\mathfrak u_{w_n})^*\otimes\mathcal O(X_f)\subseteq\mathrm{im}(\phi)$, and if $(\mathfrak u_{w_{i'}})^*\otimes\mathcal O(X_f)\subseteq\mathrm{im}(\phi)$ for all $i'>i$, then $(\mathfrak u_{w_i})^*\otimes\mathcal O(X_f)\subseteq\mathrm{im}(\phi)$. By induction, this proves $(\mathfrak u_{w_i})^*\otimes\mathcal O(X_f)\subseteq\mathrm{im}(\phi)$ for all $1\leq i\leq n$, i.e. $\phi$ is surjective. This proves $(1)\implies (2)$. 
                    
                    For $(2)\implies (1)$, assume $\gamma^*\Omega_\sigma=0$. Let $m\in\mathbb N_+$ be such that $H^0(X,L^m)$ generates $\bigoplus_{d\geq0}H^0(X,L^{md})$. For any non-nilpotent $f\in H^0(X,L^m)_{\max}$, the following map is surjective 
                    \begin{equation}
                        \begin{split}
                            \phi:\Omega_{\mathcal O(X_f)/\Bbbk}&\to \mathfrak u^*\otimes\mathcal O(X_f)\\
                            a\mathrm{d}b&\mapsto \sum_{i,j}u^{(i)}_j\otimes a\big(\xi^{(i)}_j.b\big)
                        \end{split}
                    \end{equation}
                    where $\big(\xi^{(i)}_1,\cdots,\xi^{(i)}_{r_i}\big)$ is a basis of $\mathfrak u_{w_i}$, and $\big(u^{(i)}_1,\cdots,u^{(i)}_{r_i}\big)$ is the dual basis of $(\mathfrak u_{w_i})^*$. 
        
                    For $u^{(i)}_j$, there exists $\omega^{(i)}_j\in\Omega_{\mathcal O(X_f)/\Bbbk}$ such that $\phi(\omega^{(i)}_j)=u^{(i)}_j\otimes 1$. Write 
                    \begin{equation}
                        \omega^{(i)}_j=\sum_\mu a^{(i)}_{j,\mu}\mathrm{d}b^{(i)}_{j,\mu},\quad a^{(i)}_{j,\mu},b^{(i)}_{j,\mu}\in\mathcal O(X_f). 
                    \end{equation}
                    Recall $\phi(\omega^{(i)}_j)=u^{(i)}_j\otimes 1$, that is 
                    \begin{equation}
                        \delta^{ii'}\delta_{jj'}=\sum_\mu a^{(i)}_{j,\mu}\big(\xi^{(i')}_{j'}.b^{(i)}_{j,\mu}\big),\quad 1\leq i,i'\leq n,\;1\leq j\leq r_i,\;1\leq j'\leq r_{i'}. 
                    \end{equation}
                    
                    Let $\lambda:\mathbb G_m\to\hat U$ be the grading 1PS. Then $\lambda$ acts on $\mathcal O(X_f)$. Decompose $a^{(i)}_{j,\mu},b^{(i)}_{j,\mu}$ according to $\lambda$-weights 
                    \begin{equation}
                        \begin{split}
                            a^{(i)}_{j,\mu}=&a^{(i)}_{j,\mu,\lambda=0}+a^{(i)}_{j,\mu,\lambda<0}\\
                            b^{(i)}_{j,\mu}=&b^{(i)}_{j,\mu,\lambda=-w_i}+b^{(i)}_{j,\mu,\lambda<-w_i}+b^{(i)}_{j,\mu,\lambda>-w_i}. 
                        \end{split}
                    \end{equation}
                    Then 
                    \begin{equation}
                        \begin{split}
                            \delta_{jj'}=&\sum_\mu a^{(i)}_{j,\mu}\big(\xi^{(i)}_{j',\mu}.b^{(i)}_{j,\mu}\big)\\
                            =&\sum_\mu \big(a^{(i)}_{j,\mu,\lambda=0}+a^{(i)}_{j,\mu,\lambda<0}\big)\;\xi^{(i)}_{j'}.\big(b^{(i)}_{j,\mu,\lambda=-w_i}+b^{(i)}_{j,\mu,\lambda<-w_i}+b^{(i)}_{j,\mu,\lambda>-w_i}\big)\\
                            =&\sum_\mu a^{(i)}_{j,\mu,\lambda=0}\big(\xi^{(i)}_{j'}.b^{(i)}_{j,\mu,\lambda=-w_i}\big)+\mathrm{err}
                        \end{split}
                    \end{equation}
                    where $\mathrm{err}\in\mathcal O(X_f)_{<0}$. Then $\delta_{jj'}=\xi^{(i)}_{j'}.h^{(i)}_j$ for 
                    \begin{equation}
                        h^{(i)}_j:=\sum_\mu a^{(i)}_{j,\mu,\lambda=0}b^{(i)}_{j,\mu,\lambda=-w_i}\in\mathcal O(X_f)_{-w_i}. 
                    \end{equation}
                    Since $h^{(i)}_j\in\mathcal O(X_f)_{-w_i}$, it is fixed by $\mathfrak u^{(i)}_j$ for 
                    \begin{equation}
                        \mathfrak u^{(i)}_j:=\mathrm{span}\big\{\xi^{(i)}_{j'}:1\leq j'\leq r_i,\;j'\ne j\big\}+\sum_{i'<i}\mathfrak u_{w_{i'}}. 
                    \end{equation}
                    
                    For $h^{(i)}_j$, we can write $h^{(i)}_j=\frac{g^{(i)}_j}{f^N}$ for $N\in\mathbb N_+$ and $g^{(i)}_j\in H^0\big(X,L^{mN}\big)_{\max-w_i}$. We have: 
                    \begin{itemize}
                        \item $\xi^{(i)}_j.g^{(i)}_j-f^N$ is annihilated by a power of $f$; 
                        \item for any $\xi\in\mathfrak u^{(i)}_j$, the section $\xi.g^{(i)}_j\in H^0\big(X,L^{mN}\big)_{\max}$ is annihilated by a power of $f$. 
                    \end{itemize}
                    We can enlarge $N$ and assume $\xi^{(i)}_j.g^{(i)}_j=f^N$ and $g^{(i)}_j\in H^0\big(X,L^{mN}\big)^{\mathfrak u^{(i)}_j}$. We can enlarge $N$ further such that it is uniform for all $1\leq i\leq n$ and $1\leq j\leq r_i$. In other words we have 
                    \begin{equation}
                        f^N\in \xi^{(i)}_j.H^0\big(X,L^{mN}\big)^{\mathfrak u^{(i)}_j},\quad 1\leq i\leq n,\;1\leq j\leq r_i. 
                    \end{equation}
                    
                    Choose a basis $f_1,\cdots,f_p$ of $H^0(X,L^m)_{\max}$, then there exists $N_t$ for each $1\leq t\leq p$ such that $f_t^{N_t}\in \xi^{(i)}_j.H^0(X,L^{mN_t})^{\mathfrak u^{(i)}_j}$. Take $N=\max\{N_1,\cdots,N_p\}$ and then $f_t^N\in \xi^{(i)}_j.H^0(X,L^{mN})^{\mathfrak u^{(i)}_j}$ for any $1\leq t\leq p$ and any $i,j$. This implies that $H^0(X,L^{mpN})_{\max}\subseteq \xi^{(i)}_j.H^0(X,L^{mpN})^{\mathfrak u^{(i)}_j}$ for all $i,j$, which is Condition \ref{condition of ss=s} (``ss=s''). 
                \end{proof}
                
                \begin{corollary}\label{corollary of equivalences of ss=s and stab=e}
                    Let $X$ be a projective scheme over $\Bbbk$ and let $L$ be an ample line bundle on $X$. Let $\hat U\curvearrowright (X,L)$. Then the following conditions are equivalent: 
                    \begin{itemize}
                        \item[(1)] Condition \ref{condition of ss=s} (``ss=s''); 
                        \item[(2)] $\gamma^*\Omega_\sigma=0$; 
                        \item[(3)] $\mathrm{Fit}_0(\gamma^*\Omega_\sigma)=\mathcal O_{X^0_{\min}}$; 
                        \item[(4)] $\mathrm{Stab}_U(z)=\{e\}$ for all $\Bbbk$-points $z\in Z_{\min}$. 
                    \end{itemize}
                \end{corollary}
                \begin{proof}
                    It is by definition that $(2)\iff(3)$. Proposition \ref{proposition of equivalence of ss=s and pull-back of differential vanishing} proves $(1)\iff(2)$. 
                    
                    The application of Lemma \ref{lemma of dim of relative stab and Fit when flat} to $k=0$ and $(\mathfrak v:\mathfrak w)=(\mathfrak u:0)$ shows $(3)$ is equivalent to that $\mathrm{Stab}_U(x)=\{e\}$ for all $\Bbbk$-points $x\in X^0_{\min}$. In particular $(3)\implies (4)$. 
                    
                    By the upper semi-continuity of $\dim\mathrm{Stab}_U(x)$, we have that $(4)$ implies that $\mathrm{Stab}_U(x)=\{e\}$ for all $\Bbbk$-points $x\in X^0_{\min}$, which is equivalent to $(3)$. This proves $(4)\implies(3)$. 
                \end{proof}
                
                \begin{remark}
                    When ``ss=s'' fails, the locus in $Z_{\min}$ where it fails is cut out by the sheaf of ideals $\mathrm{Fit}_0(\gamma^*\Omega_\sigma)$. This gives a subscheme structure on the locus. 
                \end{remark}
                
                \begin{remark}
                    The corollary has proved that for $d=0$ the following are equivalent:
                    \begin{itemize}
                        \item[(1)] $\gamma^*\Omega_\sigma$ is locally free of rank $d$ on $X^0_{\min}$; 
                        \item[(2)] $\dim\mathrm{Stab}_U(x)=d$ for all closed points $x\in X^0_{\min}$.
                    \end{itemize}
                    We can also prove that they are equivalent for general $d\in\mathbb N$ when $X$ is reduced. However, we only have $(1)\implies (2)$ for general $d\in\mathbb N_+$ if $X$ is not reduced. 
                \end{remark}
            }
        }
    }
    
    \subsection{U-hat theory for affine schemes}
    {
        Let $U$ be a unipotent group with Lie algebra $\mathfrak u$ and let $\lambda:\mathbb G_m\to \mathrm{Aut}(U)$ be a grading with weights $w_1>\cdots>w_n$. Consider a $\lambda$-invariant filtration of Lie ideals 
        \begin{equation}
            0=\mathfrak u_0\subseteq\mathfrak u_1\subseteq\cdots\subseteq\mathfrak u_n=\mathfrak u,\quad \mathfrak u_i:=\mathfrak u_{\lambda\geq w_i}. 
        \end{equation}
        Then each $\mathfrak u_i/\mathfrak u_{i-1}$ is Abelian, since it is graded with one weight $w_i$. Let $U_i\subseteq \hat U $ be the normal Lie subgroup corresponding to $\mathfrak u_i\subseteq\mathfrak u$. 
        
        Let $Y:=\mathrm{Spec}(A)$ be an affine scheme of finite type over $\Bbbk$. Let $\hat U $ act on $Y$. We prove under certain conditions, a universal geometric quotient of $Y$ by $U$ exists and it is affine of finite type over $\Bbbk$. 
        
        If $\mathbb G_m$ acts on $Y$, with the dual action $\rho^*:A\to\Bbbk[t,t^{-1}]\otimes A$, then the weight space of weight $w\in\mathbb Z$ is defined as $A_w:=\{a\in A:\rho^*(a)=t^{-w}\otimes a\}$. If the action of $\mathbb G_m$ on $Y$ is from $\mathbb G_m\to U\rtimes_\lambda\mathbb G_m$, then denote by $A_{\lambda=w}$ the weight space. 
        
        Define $K_0,K_1,\cdots,K_n$ as the kernels of infinitesimal actions of $\mathfrak u_i$ on $Y$ (Infinitesimal actions are defined in Definition \ref{definition of infinitesimal actions})
        \begin{equation}
            \begin{tikzcd}
                0\ar[r]&K_i\ar[r]&\Omega_{A/\Bbbk}\ar[r]&(\mathfrak u_i)^*\otimes A. 
            \end{tikzcd}
        \end{equation}
        Consider the following diagram for $1\leq i\leq n$ 
        \begin{equation}
            \begin{tikzcd}
                &&K_{i-1}\ar[d]\ar[rd,"0"]\ar[ld,dashed,"\phi_i"]\\
                0\ar[r]&(\mathfrak u_i/\mathfrak u_{i-1})^*\otimes A\ar[r]&(\mathfrak u_i)^*\otimes A\ar[r]&(\mathfrak u_{i-1})^*\otimes A\ar[r]&0
            \end{tikzcd}
        \end{equation}
        where the dashed arrow exists and denote it by $\phi_i:K_{i-1}\to(\mathfrak u_i/\mathfrak u_{i-1})^*\otimes A$. 
        
        {
            \begin{condition}\label{condition of affine CDRS by semi-gradings}
                Let $\lambda$ grade $U$ with weights $w_1>\cdots>w_n$, and let $\{\mathfrak u_i\}_{i=0}^n$ be the filtration $\mathfrak u_i:=\mathfrak u_{\lambda\geq w_i}$. Let $\hat U \curvearrowright Y=\mathrm{Spec}(A)$. This condition for the action $\hat U \curvearrowright Y$ assumes that 
                \begin{itemize}
                    \item $A=\bigoplus_{w\in\mathbb Z_{\leq 0}}A_{\lambda=w}$; 
                    \item For each $1\leq i\leq n$, there exists a surjective linear map $\mathfrak u_i/\mathfrak u_{i-1}\to\mathfrak s_i$ such that $(\mathfrak s_i)^*\otimes A$ is isomorphic to the cokernel of $\phi_i:K_{i-1}\to (\mathfrak u_i/\mathfrak u_{i-1})^*\otimes A$ in the following sense 
                    \begin{equation}
                        \begin{tikzcd}
                            &(\mathfrak s_i)^*\otimes A\ar[d]\ar[rd,"\cong"]\\
                            K_{i-1}\ar[r,"\phi_i"]&(\mathfrak u_i/\mathfrak u_{i-1})^*\otimes A\ar[r]&\mathrm{coker}(\phi_i)\ar[r]&0. 
                        \end{tikzcd}
                    \end{equation}
                \end{itemize}
            \end{condition}
            \begin{remark}
                By Corollary \ref{corollary of flatness of coker(phi_i), Q(g_i) and dim of stab when reduced}, when $Y$ is reduced, the condition above is equivalent to that for $1\leq i\leq n$ and all $\Bbbk$-points $x\in Y$, the composition $\mathrm{Stab}_{(\mathfrak u_i:\mathfrak u_{i-1})}(x)\to \mathfrak u_i/\mathfrak u_{i-1}\to \mathfrak s_i$ is an isomorphism, which implies that for $1\leq i\leq n$, the dimensions of $\mathrm{Stab}_{(\mathfrak u_i:\mathfrak u_{i-1})}(x)$ are constant for $\Bbbk$-points $x\in Y$.  
            \end{remark}
        }
        
        \subsubsection{Length one filtrations}
        {
            To construct the quotient of $Y$ by $U$, the first step is to construct the quotient of $Y$ by $U_1\subseteq \hat U $. It is equivalent to consider the case when the filtration has length one. In this subsection, Condition \ref{condition of affine CDRS by semi-gradings} in the length one case for $\hat U \curvearrowright Y$ is always assumed: 
            \begin{itemize} 
                \item $\lambda$ acts on $\mathfrak u$ with one positive weight; 
                \item $A=\bigoplus_{w\in\mathbb Z_{\leq 0}}A_{\lambda=w}$; 
                \item There exists a surjective linear map $\mathfrak u\to\mathfrak s$ such that 
                \begin{equation}
                    \begin{tikzcd}
                        &\mathfrak s^*\otimes A\ar[d]\ar[rd,"\cong"]\\
                        \Omega_{A/\Bbbk}\ar[r]&\mathfrak u^*\otimes A\ar[r]&\mathrm{coker}\ar[r]&0. 
                    \end{tikzcd}
                \end{equation}
            \end{itemize}
            
            Let $w_1\in\mathbb N_+$ be the unique weight of $\lambda$ on $\mathfrak u$. Let $\mathfrak u'\subseteq \mathfrak u$ be the kernel 
            \begin{equation}
                \begin{tikzcd}
                    0\ar[r]&\mathfrak u'\ar[r]&\mathfrak u\ar[r]&\mathfrak s\ar[r]&0
                \end{tikzcd}
            \end{equation}
            and choose a splitting $\mathfrak u=\mathfrak u'\oplus\mathfrak s$. Let $\big(\xi_1,\cdots,\xi_r\big)$ be a basis of $\mathfrak u$ such that $\big(\xi_1,\cdots,\xi_{r-k}\big)$ is a basis of $\mathfrak u'$. Let $\big(u_1,\cdots ,u_r\big)$ be the basis of $\mathfrak u^*$ dual to $\big(\xi_1,\cdots,\xi_r\big)$. 
            
            {
                \begin{lemma}\label{lemma of coincidence of invariants for affine CDRS length one filtration}
                    Assume Condition \ref{condition of affine CDRS by semi-gradings}. Let $\mathfrak u'$ be as above. Then $A^{\mathfrak u}=A^{\mathfrak u'}$. 
                \end{lemma}
                \begin{proof}
                    Since $\mathfrak u'\subseteq\mathfrak u$, we have $A^{\mathfrak u}\subseteq A^{\mathfrak u'}$. 
                    
                    Let $\mathrm{d}:A\to\Omega_{A/\Bbbk}$ be the derivation. Let $i:\mathfrak u'\to \mathfrak u$ be the inclusion and let $q:\mathfrak u\to\mathfrak s$ be the quotient. Let $i^*:\mathfrak u^*\to(\mathfrak u')^*$ and $q^*:\mathfrak s^*\to\mathfrak u^*$ be the dual maps. 
                    
                    Consider the following diagram with arrows labelled as $\phi,\varphi$ 
                    \begin{equation}\label{equation of diagram for affine CDRS length one filtration}
                        \begin{tikzcd}
                            &\mathfrak s^*\otimes A\ar[d,"q^*\otimes 1"]\ar[rd,"\cong"]\\
                            \Omega_{A/\Bbbk}\ar[r,"\phi"]&\mathfrak u^*\otimes A\ar[r,"\varphi"]&\mathrm{coker}\ar[r]&0. 
                        \end{tikzcd}
                    \end{equation}
                    
                    We prove $\ker\big((i^*\otimes1)\circ\phi\big)\subseteq\ker\phi$. Let $\omega\in\ker\big((i^*\otimes 1)\circ\phi\big)\subseteq\Omega_{A/\Bbbk}$. Then $\phi(\omega)\in\ker(i^*\otimes 1)=\mathrm{im}(q^*\otimes 1)$. There exists $g\in \mathfrak s^*\otimes A$ such that $\phi(\omega)=(q^*\otimes 1)(g)$. Since $\varphi\circ (q^*\otimes 1)$ is an isomorphism and $\varphi\circ(q^*\otimes 1)(g)=\varphi\circ \phi(\omega)=0$, we have $g=0$. Then $\phi(\omega)=(q^*\otimes 1)(g)=0$. This proves $\ker\big((i^*\otimes1)\circ\phi\big)\subseteq\ker\phi$. 
                    
                    Let $f\in A^{\mathfrak u'}$, i.e. $(i^*\otimes 1)\circ\phi(\mathrm{d}f)=0$. Since $\ker\big((i^*\otimes1)\circ\phi\big)\subseteq\ker\phi$, we have $\phi(\mathrm{d}f)=0$, i.e. $f\in A^{\mathfrak u}$. This proves $A^{\mathfrak u'}\subseteq A^{\mathfrak u}$. 
                \end{proof}
            }
            
            {
                \begin{lemma}\label{lemma of functions derivative equal delta_ij for affine CDRS length one filtration}
                    Assume Condition \ref{condition of affine CDRS by semi-gradings}. Let $w_1$, $\mathfrak u'$, $\{\xi_i\}$ and $\{u_i\}$ be as above. Then there exist $f_1,\cdots,f_{r-k}\in A_{\lambda=-w_1}$ such that 
                    \begin{equation}
                        \xi_i.f_j=\delta_{ij},\quad 1\leq i,j\leq r-k. 
                    \end{equation}
                \end{lemma}
                \begin{proof}
                    For any $w\in\mathbb Z$, and any $\xi\in\mathfrak u$, we have $\xi.A_{\lambda=w}\subseteq A_{\lambda=w+w_1}$, where $w_1$ is the $\lambda$-weight on $\mathfrak u$. In particular elements in $A_{\lambda=0}$ are fixed by $\mathfrak u$. 
    
                    Consider diagram \eqref{equation of diagram for affine CDRS length one filtration}. We have $\mathfrak u^*\otimes A=\mathrm{im}\phi+\mathrm{im}(q^*\otimes 1)$. In particular, for $1\leq i\leq r-k$, there exists 
                    \begin{equation}
                        \sum_\mu a_{i,\mu}\mathrm{d}b_{i,\mu }\in\Omega_{A/\Bbbk}
                    \end{equation}
                    such that 
                    \begin{equation}
                        u_i\otimes 1-\sum_{i'=1}^ru_{i'}\otimes \Big(\sum_\mu a_{i,\mu}\big(\xi_{i'}.b_{i,\mu}\big)\Big)\in \mathfrak s^*\otimes A. 
                    \end{equation}
                    
                    Compare the coefficients of $u_{i'}\otimes 1$ for $1\leq i'\leq r-k$
                    \begin{equation}
                        \delta_{ii'}=\sum_\mu a_{i,\mu}\big(\xi_{i'}.b_{i,\mu }\big),\quad 1\leq i,i'\leq r-k. 
                    \end{equation}
                    
                    Decompose $a_{i,\mu},b_{i,\mu}$ according to $\lambda$-weights 
                    \begin{equation}
                        \begin{split}
                            a_{i,\mu}=&a_{i,\mu,\lambda=0}+a_{i,\mu,\lambda<0}\\
                            b_{i,\mu}=&b_{i,\mu,\lambda>-w_1}+b_{i,\mu,\lambda=-w_1}+b_{i,\mu,\lambda<-w_1}. 
                        \end{split}
                    \end{equation}
                    Then 
                    \begin{equation}
                        \begin{split}
                            \delta_{ii'}=&\sum_\mu a_{i,\mu}\big(\xi_{i'}.b_{i,\mu}\big)\\
                            =&\sum_\mu a_{i,\mu,\lambda=0}\big(\xi_{i'}.b_{i,\mu,\lambda=-w_1}\big)+\sum_\mu a_{i,\mu,\lambda=0}\big(\xi_{i'}.b_{i,\mu,\lambda<-w_1}\big)+\sum_\mu a_{i,\mu,\lambda<0}\big(\xi_{i'}.b_{i,\mu}\big). 
                        \end{split}
                    \end{equation}
                    Compare the components of $\lambda$-weight zero 
                    \begin{equation}
                        \delta_{ii'}=\sum_\mu a_{i,\mu,\lambda=0}\big(\xi_{i'}.b_{i,\mu,\lambda=-w_1}\big)=\xi_{i'}.\Big(\sum_\mu a_{i,\mu,\lambda=0}b_{i,\mu,\lambda=-w_1}\Big). 
                    \end{equation}
                    Let $f_i:=\sum_\mu a_{i,\mu,\lambda=0}b_{i,\mu,\lambda=-w_1}\in A_{\lambda=-w_1}$ and the lemma is proved. 
                \end{proof}
            }
    
            The following lemma is essentially from \cite{dixmier1996enveloping} Lemma 4.7.5, adapted to a set of commuting derivations. 
            {
                \begin{lemma}[cf. \cite{berczi2016projective} Lemma 7.2, \cite{dixmier1996enveloping} Lemma 4.7.5]\label{lemma of additive action with delta-preimages of 1 under derivations imply trivial quotient}
                    Let $R$ be an algebra over $\Bbbk$ with a dual action of $(\mathbb G_a)^r$. Let $\xi_1,\cdots,\xi_r$ be commuting derivations on $R$ generating the dual action of $(\mathbb G_a)^r$. Assume there exists $f_1,\cdots,f_r\in R$ such that $\xi_i.f_j=\delta_{ij}$, then: 
                    \begin{itemize}
                        \item $f_1,\cdots,f_r$ are algebraically independent over $R^{(\mathbb G_a)^r}$; 
                        \item $R=R^{(\mathbb G_a)^r}[f_1,\cdots,f_t]$; 
                        \item $\mathrm{Spec}(R)\to\mathrm{Spec}\big(R^{(\mathbb G_a)^r}\big)$ is a trivial $(\mathbb G_a)^r$-quotient 
                        \begin{equation}
                            \begin{tikzcd}
                                \mathrm{Spec}(R)\ar[r,"\cong"]\ar[rd]&(\mathbb G_a)^r\times_{\mathrm{Spec}(\Bbbk)}\mathrm{Spec}\big(R^{(\mathbb G_a)^r}\big)\ar[d,"\mathrm{pr}"]\\
                                &\mathrm{Spec}\big(R^{(\mathbb G_a)^r}\big). 
                            \end{tikzcd}
                        \end{equation}
                    \end{itemize}
                \end{lemma}
                \begin{proof}
                    Let $I\subseteq R$ be the ideal generated by $f_1,\cdots,f_r$. Let $\alpha:R\to R/I$ be the quotient map. Consider the following linear map 
                    \begin{equation}
                        \begin{split}
                            \phi:R&\to (R/I)[t_1,\cdots,t_r]\\
                            g&\mapsto \sum_{n\in\mathbb N^r}\frac{1}{n!}\alpha(\xi^n.g)t^n
                        \end{split}
                    \end{equation}
                    where 
                    \begin{equation}
                        \begin{split}
                            n:=&(n_1,\cdots,n_r)\\
                            n!:=&n_1!\cdots n_r!\\
                            \xi^n.g:=&\xi_1^{n_1}\cdots\xi_r^{n_r}.g\\
                            t^n:=&t_1^{n_1}\cdots t_r^{n_r}. 
                        \end{split}
                    \end{equation}
                    
                    We need to show that the summation defining $\phi$ is essentially finite. For any $g\in R$, there exists a finite dimensional invariant subspace $V\subseteq R$ containing $g$. Since $(\mathbb G_m)^r$ is unipotent, and $V$ is a finite dimensional representation, $\xi^n|_V=0$ if $|n|:=n_1+\cdots+n_r>\dim V$. So $\xi^n.g=0$ for all but finite $n\in\mathbb N^r$. 
                    
                    The linear map $\phi$ preserves multiplication. For $f,g\in R$, we have 
                    \begin{equation}
                        \begin{split}
                            \phi(fg)=&\sum_{n\in\mathbb N^r}\frac{1}{n!}\alpha\big(\xi^n.(fg)\big)t^n\\
                            =&\sum_{n\in\mathbb N^r}\frac{1}{n!}t^n\alpha\bigg(\sum_{0\leq m\leq n}\binom{n}{m}(\xi^{n-m}.f)(\xi^{m}.g)\bigg)\\
                            =&\sum_{u,v\in\mathbb N^r}\frac{1}{u!}\frac{1}{v!}t^{u+v}\alpha(\xi^u.f)\alpha(\xi^v.g)\\
                            =&\phi(f)\phi(g). 
                        \end{split}
                    \end{equation}
                    So $\phi$ is a map of $\Bbbk$-algebras. It is easy to check $\phi(f_i)=t_i$. 
                    
                    We show $\phi$ is an isomorphism by constructing its inverse. Consider 
                    \begin{equation}
                        \begin{split}
                            \varphi':R[t_1,\cdots,t_r]&\to R\\
                            gt^m&\mapsto \sum_{n\in\mathbb N^r}\frac{(-1)^{|n|}}{n!}(\xi^n.g)f_1^{m_1+n_1}\cdots f_r^{m_r+n_r}. 
                        \end{split}
                    \end{equation}
                    By the same argument for $\phi$, we have that the summation defining $\varphi'$ is finite, and $\varphi'$ is a ring map. It is easy to see $\varphi'(f_i)=0$ for $i=1,\cdots,r$. So $\varphi'$ factors through $R[t_1,\cdots,t_r]\to (R/I)[t_1,\cdots,t_r]$, inducing 
                    \begin{equation}
                        \varphi:(R/I)[t_1,\cdots,t_r]\to R
                    \end{equation}
                    which is the inverse of $\phi$. 
                    
                    Via the isomorphism $R\cong (R/I)[t_1,\cdots,t_r]$, the derivation $\xi_i$ on $R$ corresponds to the derivation $\frac{\partial}{\partial t_i}$ on $(R/I)[t_1,\cdots,t_r]$. Then 
                    \begin{equation}
                        R^{(\mathbb G_a)^r}\cong R/I,\quad g\mapsto \alpha(g). 
                    \end{equation}
                    The isomorphism $\varphi:R/I[t_1,\cdots,t_r]\cong R$ then becomes 
                    \begin{equation}
                        R^{(\mathbb G_a)^r}[t_1,\cdots,t_r]\cong R,\quad gt^n\mapsto gf_1^{n_1}\cdots f_r^{n_r}; 
                    \end{equation}
                    that is 
                    \begin{equation}
                        R=R^{(\mathbb G_a)^r}[f_1,\cdots,f_r]
                    \end{equation}
                    and $f_1,\cdots,f_r$ are algebraically independent over $R^{(\mathbb G_a)^r}$. 
        
                    There is an isomorphism $\Bbbk[t_1,\cdots,t_r]\cong\mathcal O\big((\mathbb G_a)^r\big)$. The left action of $(\mathbb G_a)^r$ on itself induces a dual action of $(\mathbb G_a)^r$ on $\mathcal O\big((\mathbb G_a)^r\big)$. Through the identification $\mathcal O\big((\mathbb G_a)^r\big)\cong \Bbbk[t_1,\cdots,t_r]$, the induced dual action of $(\mathbb G_a)^r$ on $\Bbbk[t_1,\cdots,t_r]$ is generated by derivations $\frac{\partial}{\partial t_1},\cdots,\frac{\partial}{\partial t_r}$. 
                    
                    The following isomorphism is $(\mathbb G_a)^r$-equivariant 
                    \begin{equation}
                        \begin{split}
                            R\cong &R^{(\mathbb G_a)^r}[t_1,\cdots,t_r]\\
                            \cong &\Bbbk[t_1,\cdots,t_r]\otimes_{\Bbbk}R^{(\mathbb G_a)^r}\\
                            \cong &\mathcal O\big((\mathbb G_a)^r\big)\otimes_{\Bbbk}R^{(\mathbb G_a)^r}. 
                        \end{split}
                    \end{equation}
                    Geometrically, there is a $(\mathbb G_a)^r$-equivariant isomorphism 
                    \begin{equation}
                        \begin{tikzcd}
                            \mathrm{Spec}(R)\ar[r,"\cong"]\ar[rd]&(\mathbb G_a)^r\times_{\mathrm{Spec}(\Bbbk)}\mathrm{Spec}\big(R^{(\mathbb G_a)^r}\big)\ar[d,"\mathrm{pr}"]\\
                            &\mathrm{Spec}\big(R^{(\mathbb G_a)^r}\big). 
                        \end{tikzcd}
                    \end{equation}
                \end{proof}
            }
            
            {
                \begin{lemma}\label{lemma of universal geometric quotient by U for affine CDRS length one filtration}
                    Assume Condition \ref{condition of affine CDRS by semi-gradings}. The inclusion $A^U\to A$ induces a universal geometric quotient by $U$ 
                    \begin{equation}
                        \mathrm{Spec}(A)\to \mathrm{Spec}\big(A^U\big). 
                    \end{equation}
                \end{lemma}
                \begin{proof}
                    By Lemma \ref{lemma of functions derivative equal delta_ij for affine CDRS length one filtration}, there exist $f_1,\cdots,f_{r-k}\in A_{\lambda=-w_1}$ such that $\xi_i.f_j=\delta_{ij}$ for $1\leq i,j\leq r-k$. Since $\mathfrak u$ is Abelian, the subspace $\mathfrak u'\subseteq\mathfrak u$ is a Lie sub-algebra. Let $U'\subseteq U$ be the subgroup corresponding to $\mathfrak u'\subseteq\mathfrak u$. The derivations $\xi_1,\cdots,\xi_{r-k}$ on $A$ generate the dual action of $U'$ on $A$. By Lemma \ref{lemma of coincidence of invariants for affine CDRS length one filtration} and Lemma \ref{lemma of additive action with delta-preimages of 1 under derivations imply trivial quotient}, we have: 
                    \begin{itemize}
                        \item $f_1,\cdots,f_{r-k}$ are algebraically independent over $A^U$; 
                        \item $A=A^U[f_1,\cdots,f_{r-k}]$; 
                        \item $\mathrm{Spec}(A)\to\mathrm{Spec}\big(A^U\big)$ is a trivial $U'$-quotient. 
                    \end{itemize}
                    
                    Let $A^U\to B$ be any $A^U$-algebra. Consider the Cartesian diagram 
                    \begin{equation}
                        \begin{tikzcd}
                            \mathrm{Spec}(B')\ar[r]\ar[d,"\pi"]&\mathrm{Spec}(A)\ar[d]\\
                            \mathrm{Spec}(B)\ar[r]&\mathrm{Spec}\big(A^U\big)
                        \end{tikzcd}
                    \end{equation}
                    where 
                    \begin{equation}
                        B'=B\otimes_{A^U}A. 
                    \end{equation}
                    
                    We first show that $B\to B'$ is injective and $B=(B')^U$. As a pull back of a trivial $U'$-quotient, the morphism $\pi:\mathrm{Spec}(B')\to\mathrm{Spec}(B)$ is a trivial $U'$-quotient. In particular, $\pi$ is a geometric quotient by $U'$. Since $\pi$ is surjective, $B\to B'$ is injective. Since $\pi$ is a good quotient, $(B')^{U'}=B$. Consider $U$-invariant elements on both sides, and we have 
                    \begin{equation}
                        (B')^U=B
                    \end{equation}
                    since $B$ is acted trivially by $U$. 
                    
                    We then show $\pi:\mathrm{Spec}(B')\to \mathrm{Spec}(B)$ is a geometric quotient by $U$. The conditions to check are: 
                    \begin{itemize}
                        \item $\pi$ is affine surjective and $U$-invariant. By the above argument, $\pi$ is surjective, and $B=(B')^U$, which implies that $\pi$ is $U$-invariant. 
        
                        \item For any principal affine open $V:=\mathrm{Spec}(B_f)\to\mathrm{Spec}(B)$ with $f\in B$, we can replace $A^U\to B$ by $A^U\to B_f$ and apply the above argument to claim $B_f\subseteq B_f\otimes_{A^U}A\cong B'_f$ and 
                        \begin{equation}
                            B_f=\big(B_f\otimes_{A^U}A\big)^{U}\cong (B'_f)^{U}
                        \end{equation}
                        i.e. 
                        \begin{equation}
                            \mathcal O_{\mathrm{Spec}(B)}(V)\cong \mathcal O_{\mathrm{Spec}(B')}\big(\pi^{-1}(V)\big)^U. 
                        \end{equation}
        
                        \item For any $U$-invariant closed subset $W\subseteq\mathrm{Spec}(B')$, it is also $U'$-invariant, so $\pi(W)$ is closed since $\pi$ is a good quotient by $U'$. 
                        
                        \item For any pair of disjoint $U$-invariant closed subsets $W_1,W_2\subseteq\mathrm{Spec}(B')$, they are also $U'$-invariant, so $\pi(W_1)\cap \pi(W_2)=\emptyset$. 
                        
                        \item For any closed point $z\in\mathrm{Spec}(B)$, the fibre $\pi^{-1}(z)\hookrightarrow\mathrm{Spec}(B')$ is a union of $U$-orbits, since $\pi$ is $U$-invariant. Besides, $\pi^{-1}(z)$ is a single $U'$-orbit, since $\pi$ is a geometric quotient by $U'$. So there is only one $U$-orbit in $\pi^{-1}(z)$. 
                    \end{itemize}
                    
                    We see that any affine base change of $\mathrm{Spec}(A)\to\mathrm{Spec}\big(A^U\big)$ is a geometric quotient by $U$. Since arbitrary base change is covered by affine base changes as open subschemes, and geometric quotients are local on target, we claim any base change of $\mathrm{Spec}(A)\to\mathrm{Spec}\big(A^U\big)$ is a geometric quotient by $U$. This proves that $\mathrm{Spec}(A)\to\mathrm{Spec}\big(A^U\big)$ is a universal geometric quotient by $U$. 
                \end{proof}
            }
    
            {
                \begin{theorem}\label{theorem of affine quotients for affine CDRS length one filtration}
                    Let $Y=\mathrm{Spec}(A)$ be an affine scheme of finite type over $\Bbbk$. Let $\lambda:\mathbb G_m\to\mathrm{Aut}(U)$ be a grading of one weight. Let $\hat U \curvearrowright Y$. Assume Condition \ref{condition of affine CDRS by semi-gradings} for the action $\hat U \curvearrowright Y$. Then: 
                    \begin{itemize}
                        \item[(1)] $Y/U:=\mathrm{Spec}\big(A^U\big)$ is affine of finite type over $\Bbbk$; 
                        \item[(2)] $Y\to Y/U$ is a universal geometric quotient by $U$; 
                        \item[(3)] There exists an affine space $\mathbb A^m$ such that $Y\cong \mathbb A^m\times Y/U$ over $Y/U$, and $m=\dim \mathfrak u-\dim\mathfrak s$; 
                        \item[(4)] The map $\Omega_{A/\Bbbk}\to \mathfrak u^*\otimes A$ induces an injective map $\Omega_{A/A^U}\to\mathfrak u^*\otimes A$ such that the following diagram commutes 
                        \begin{equation}
                            \begin{tikzcd}
                                \Omega_{A/\Bbbk}\ar[r]\ar[d]&\mathfrak u^*\otimes A\\
                                \Omega_{A/A^U}\ar[ru,"\subseteq"{sloped}].
                            \end{tikzcd}
                        \end{equation}
                    \end{itemize}
                \end{theorem}
                \begin{proof}
                    Let $w_1\in\mathbb N_+$ be the unique $\lambda$-weight on $\mathfrak u$. Let $\mathfrak u\to\mathfrak s$ be the surjective linear map in Condition \ref{condition of affine CDRS by semi-gradings}. Let $\mathfrak u'\subseteq\mathfrak u$ be the kernel of $\mathfrak u\to\mathfrak s$. Choose a decomposition $\mathfrak u=\mathfrak u'\oplus\mathfrak s$. Choose a basis $\big(\xi_1,\cdots,\xi_r\big)$ of $\mathfrak u$ such that $\xi_j\in\mathfrak u'$ for $1\leq j\leq r-k$ and $\xi_j\in\mathfrak s$ for $r-k<j\leq r$, where $r=\dim\mathfrak u$, $k=\dim\mathfrak s$. 
    
                    Let $U'\subseteq U$ be the subgroup with Lie algebra $\mathrm{Lie}(U')=\mathfrak u'$. Consider the action of $U'$ on $Y=\mathrm{Spec}(A)$. We have that $\xi_1,\cdots,\xi_{r-k}$ are commuting derivations on $A$ generating the dual action of $U'$ on $A$. By Lemma \ref{lemma of functions derivative equal delta_ij for affine CDRS length one filtration}, there exist $f_1,\cdots,f_{r-k}\in A_{\lambda=-w_1}$ such that $\xi_i.f_j=\delta_{ij}$ for $1\leq i,j\leq r-k$. By Lemma \ref{lemma of coincidence of invariants for affine CDRS length one filtration}, Lemma \ref{lemma of additive action with delta-preimages of 1 under derivations imply trivial quotient} and Lemma \ref{lemma of universal geometric quotient by U for affine CDRS length one filtration}, we have: 
                    \begin{itemize}
                        \item $f_1,\cdots,f_{r-k}\in A$ are algebraically independent over $A^U$; 
                        \item $A=A^U[f_1,\cdots,f_{r-k}]$; 
                        \item $\mathrm{Spec}(A)\to\mathrm{Spec}\big(A^U\big)$ is a universal geometric quotient by $U$. 
                    \end{itemize}
                    Then $(1)$ and $(2)$ are immediate. 
                    
                    For $(3)$, let $x_1,\cdots,x_{r-k}$ be free variables, and let 
                    \begin{equation}
                        S=\Bbbk[x_1,\cdots,x_{r-k}]
                    \end{equation}
                    be a polynomial ring over $\Bbbk$. Then the following map is an isomorphism 
                    \begin{equation}
                        A^U\otimes_{\Bbbk}S\to A,\quad g\otimes p(x_1,\cdots,x_{r-k})\mapsto gp(f_1,\cdots,f_{r-k}). 
                    \end{equation}
                    
                    Let $\mathbb A^m:=\mathrm{Spec}(S)$. Then $m=r-k=\dim\mathfrak u-\dim\mathfrak s$.  Geometrically there is an isomorphism over $Y/U$ 
                    \begin{equation}
                        \begin{tikzcd}
                            Y\ar[r,"\cong"]\ar[rd]&\mathbb A^m\times Y/U\ar[d,"\mathrm{pr}"]\\
                            &Y/U
                        \end{tikzcd}
                    \end{equation}
                    which proves $(3)$. 
                    
                    For $(4)$, consider the diagram 
                    \begin{equation}
                        \begin{tikzcd}
                            \Omega_{A^U/\Bbbk}\otimes_{A^U}A\ar[r,"\phi"]\ar[rd,"\psi\circ\phi"]&\Omega_{A/\Bbbk}\ar[r]\ar[d,"\psi"]&\Omega_{A/A^U}\ar[r]\ar[ld,dashed,"\varphi"]&0\\
                            &\mathfrak u^*\otimes A
                        \end{tikzcd}
                    \end{equation}
                    where the top row is exact by  \cite[\href{https://stacks.math.columbia.edu/tag/00RS}{Lemma 00RS}]{stacks-project} and $\psi\circ\phi=0$ by direct check. Therefore $\varphi$ exists. 
                    
                    We have 
                    \begin{equation}
                        \Omega_{A/A^U}=\bigoplus_{j=1}^{r-k}A\mathrm{d}f_j
                    \end{equation}
                    since $A=A^U[f_1,\cdots,f_{r-k}]$ and $f_1,\cdots,f_{r-k}$ are algebraically independent over $A^U$. Let
                    \begin{equation}
                        \omega:=\sum_{j=1}^{r-k}a_j\mathrm{d}f_j\in\Omega_{A/A^U}
                    \end{equation}
                    be an element in $\ker\varphi$. It suffices to prove $a_j=0$ for $1\leq j\leq r-k$. For any $1\leq j'\leq r-k$, the coefficient of $u_{j'}\otimes 1$ of $\varphi(\omega)$ equals 
                    \begin{equation}
                        0=\sum_{j=1}^{r-k}a_j(\xi_{j'}.f_j)=a_{j'}
                    \end{equation}
                    where we have used $\xi_{j'}.f_j=\delta_{j'j}$. This proves $(4)$. 
                \end{proof}
            }
    
        }
    
        \subsubsection{Induction for general filtrations}
        {
            Consider the general case when $\lambda:\mathbb G_m\to\mathrm{Aut}(U)$ has weights $w_1>\cdots>w_n$. We will prove in this subsubsection that if the action $\hat U \curvearrowright Y=\mathrm{Spec}(A)$ satisfies Condition \ref{condition of affine CDRS by semi-gradings}, then the induced action $(U/U_1)\rtimes \lambda\curvearrowright Y/U_1$ satisfies Condition \ref{condition of affine CDRS by semi-gradings}. This will allow us to iteratively apply Theorem \ref{theorem of affine quotients for affine CDRS length one filtration} (See Theorem \ref{theorem of affine U-quotient for CDRS} below ). 
            
            {
                \begin{proposition}\label{proposition of induction for affine CDRS}
                    If $\hat U\curvearrowright Y:=\mathrm{Spec}(A)$ satisfies Condition \ref{condition of affine CDRS by semi-gradings}, then so does the induced action $(U/U_1)\rtimes \lambda\curvearrowright Y/U_1:=\mathrm{Spec}\big(A^{U_1}\big)$. 
                \end{proposition}
                \begin{proof}
                    If Condition \ref{condition of affine CDRS by semi-gradings} is satisfied for $\hat U \curvearrowright Y$, then this condition is also satisfied for $U_1\rtimes \lambda\curvearrowright Y$. Apply Theorem \ref{theorem of affine quotients for affine CDRS length one filtration} to $U_1\rtimes \lambda\curvearrowright Y$ to conclude: 
                    \begin{itemize}
                        \item $Y/U_1:=\mathrm{Spec}\big(A^{U_1}\big)$ is affine of finite type over $\Bbbk$; 
                        \item $Y\to Y/U_1$ is a universal geometric quotient by $U_1$; 
                        \item $Y\cong \mathbb A^m\times Y/U_1$ over $Y/U_1$; 
                        \item $\Omega_{A/A^{U_1}}\to(\mathfrak u_1)^*\otimes A$ is injective. 
                    \end{itemize}
                    
                    The actions $\lambda\curvearrowright U/U_1$ and $(U/U_1)\rtimes \lambda\curvearrowright Y/U_1$ are induced naturally. We then check Condition \ref{condition of affine CDRS by semi-gradings} for the action $(U/U_1)\rtimes \lambda\curvearrowright Y/U_1$. The following are immediate: 
                    \begin{itemize}
                        \item $\lambda$ is a grading on $U/U_1$ with weights $w_2>\cdots>w_n$; 
                        \item The filtration $\big\{\mathfrak u_i/\mathfrak u_1\big\}_{i=1}^n$ is such that $\mathfrak u_i/\mathfrak u_1=(\mathfrak u/\mathfrak u_1)_{\lambda\geq w_i}$; 
                        \item $A^{U_1}=\bigoplus_{w\in\mathbb Z_{\leq 0}}\big(A^{U_1}\big)_{\lambda=w}$. 
                    \end{itemize}
                    
                    We check the last statement in Condition \ref{condition of affine CDRS by semi-gradings}. Let $\mathfrak u_i/\mathfrak u_{i-1}\to\mathfrak s_i$ be the quotient in Condition \ref{condition of affine CDRS by semi-gradings} for $1\leq i\leq n$. We will prove that $\varphi'$ is an isomorphism in the following diagram for $2\leq i\leq n$ 
                    \begin{equation}\label{equation of diagram of inductive CDRS}
                        \begin{tikzcd}
                            &(\mathfrak s_i)^*\otimes A^{U_1}\ar[d]\ar[rd,"\varphi'"]\\
                            K'_{i-1}\ar[r,"\phi'_i"]&(\mathfrak u_i/\mathfrak u_{i-1})^*\otimes A^{U_1}\ar[r]&\mathrm{coker}(\phi'_i)\ar[r]&0
                        \end{tikzcd}
                    \end{equation}
                    where $K'_i:=\ker\big(\Omega_{A^{U_1}/\Bbbk}\to (\mathfrak u_i/\mathfrak u_{i-1})^*\otimes A^{U_1}\big)$ for $1\leq i\leq n$. 

                    By Theorem \ref{theorem of affine quotients for affine CDRS length one filtration}, the inclusion $A^{U_1}\to A$ is a polynomial ring, so faithfully flat. Then $\varphi$ is an isomorphism if and only if its base change along $A^{U_1}\to A$ is an isomorphism. The base change of diagram \eqref{equation of diagram of inductive CDRS} is 
                    \begin{equation}
                        \begin{tikzcd}
                            &(\mathfrak s_i)^*\otimes A\ar[d]\ar[rd,"\varphi'\otimes 1"]\\
                            K'_{i-1}\otimes_{A^{U_1}}A\ar[r,"\phi'_i\otimes 1"]&(\mathfrak u_i/\mathfrak u_{i-1})^*\otimes A\ar[r]&\mathrm{coker}(\phi'_i)\otimes_{A^{U_1}}A\ar[r]&0. 
                        \end{tikzcd}
                    \end{equation}
                    
                    Condition \ref{condition of affine CDRS by semi-gradings} assumes that $\varphi$ is an isomorphism in the diagram 
                    \begin{equation}
                        \begin{tikzcd}
                            &(\mathfrak s_i)^*\otimes A\ar[d]\ar[rd,"\varphi"]\\
                            K_{i-1}\ar[r,"\phi_i"]&(\mathfrak u_i/\mathfrak u_{i-1})^*\otimes A\ar[r]&\mathrm{coker}(\phi_i)\ar[r]&0.
                        \end{tikzcd}
                    \end{equation}
                    
                    Consider the following diagram with exact rows: 
                    \begin{equation}
                        \begin{tikzcd}
                            &\Omega_{A^{U_1}/\Bbbk}\otimes_{A^{U_1}}A\ar[r]\ar[d]&\Omega_{A/\Bbbk}\ar[r]\ar[d]&\Omega_{A/A^{U_1}}\ar[r]\ar[d,"\subseteq"{sloped}]&0\\
                            0\ar[r]&(\mathfrak u_{i-1}/\mathfrak u_1)^*\otimes A\ar[r]&(\mathfrak u_{i-1})^*\otimes A\ar[r]&(\mathfrak u_1)^*\otimes A\ar[r]&0
                        \end{tikzcd}
                    \end{equation}
                    where the exactness of the top row is by \cite[\href{https://stacks.math.columbia.edu/tag/00RS}{Lemma 00RS}]{stacks-project}. By Theorem \ref{theorem of affine quotients for affine CDRS length one filtration}, the right vertical arrow is injective. Since $A^{U_1}\to A$ is flat, the kernel of the left vertical arrow is $K'_{i-1}\otimes_{A^{U_1}}A$. The following is exact by the Snake lemma 
                    \begin{equation}
                        \begin{tikzcd}
                            K'_{i-1}\otimes_{A^{U_1}}A\ar[r]&K_{i-1}\ar[r]&0. 
                        \end{tikzcd}
                    \end{equation}
                    This surjective map fills the following diagram 
                    \begin{equation}
                        \begin{tikzcd}
                            &(\mathfrak s_i)^*\otimes A\ar[d]\ar[rd,"\varphi'\otimes 1"]\ar[bend right=30,rdd,"\varphi"]\\
                            K'_{i-1}\otimes_{A^{U_1}}A\ar[r,"\phi'_i\otimes1"]\ar[d,->>]&(\mathfrak u_i/\mathfrak u_{i-1})^*\otimes A\ar[r]\ar[d,equal]&\mathrm{coker}(\phi'_i)\otimes_{A^{U_1}}A\ar[r]\ar[d,dashed,"\cong"]&0\\
                            K_{i-1}\ar[r,"\phi_i"]&(\mathfrak u_i/\mathfrak u_{i-1})^*\otimes A\ar[r]&\mathrm{coker}(\phi'_i)\ar[r]&0
                        \end{tikzcd}
                    \end{equation}
                    so an isomorphism between cokernels is induced. Therefore $\varphi'\otimes 1$ is an isomorphism since $\varphi$ is. This proves the last statement in Condition \ref{condition of affine CDRS by semi-gradings}, and hence completes the proof. 
                \end{proof}
            }

            {
                \begin{theorem}\label{theorem of affine U-quotient for CDRS}
                    Let $Y=\mathrm{Spec}(A)$ be an affine scheme of finite type over $\Bbbk$. Assume Condition \ref{condition of affine CDRS by semi-gradings} for the action $\hat U \curvearrowright Y$. Then: 
                    \begin{itemize}
                        \item[(1)] $Y/U:=\mathrm{Spec}\big(A^U\big)$ is affine of finite type over $\Bbbk$; 
                        \item[(2)] $Y\to Y/U$ is a universal geometric quotient by $U$; 
                        \item[(3)] There exists an affine space $\mathbb A^m$ such that $Y\cong\mathbb A^m\times Y/U$ over $Y/U$. 
                    \end{itemize}
                \end{theorem}
                \begin{proof}
                    For $1\leq i\leq n$, we prove by induction on $i$ the following: 
                    \begin{itemize}
                        \item[(1)] $Y/U_i:=\mathrm{Spec}\big(A^{U_i}\big)$ is of finite type over $\Bbbk$; 
                        \item[(2)] $Y\to Y/U_i$ is a universal geometric quotient by $U_i$; 
                        \item[(3)] There exists an affine space $\mathbb A^{m_i}$ such that $Y\cong \mathbb A^{m_i}\times Y/U_i$ over $Y/U_i$; 
                        \item[(4)] The action $(U/U_i)\rtimes \lambda\curvearrowright Y/U_i$ satisfies Condition \ref{condition of affine CDRS by semi-gradings}. 
                    \end{itemize}
                    
                    For $i=1$, Proposition \ref{proposition of induction for affine CDRS} proves $(1)-(4)$. Assume the above are proved for some $1\leq i<n$. Apply Thereom \ref{theorem of affine quotients for affine CDRS length one filtration} to the action $(U_{i+1}/U_i)\rtimes \lambda\curvearrowright Y/U_i$ to conclude: 
                    \begin{itemize}
                        \item[(1')] $Y/U_{i+1}:=\mathrm{Spec}\Big(\big(A^{U_i}\big)^{U_{i+1}/U_i}\Big)=\mathrm{Spec}\big(A^{U_{i+1}}\big)$ is of finite type over $\Bbbk$; 
                        \item[(2')] $Y/U_i\to Y/U_{i+1}$ is a universal geometric quotient by $U_{i+1}/U_i$, so $Y\to Y/U_{i+1}$ is a universal geometric quotient by $U_{i+1}$ as the composition of $Y\to Y/U_i$ and $Y/U_i\to Y/U_{i+1}$; 
                        \item[(3')] There exists an affine $\mathbb A^{m_{i+1}'}$ such that $Y/U_i\cong \mathbb A^{m_{i+1}'}\times Y/U_{i+1}$ over $Y/U_{i+1}$, so $Y\cong \mathbb A^{m_i}\times Y/U_i\cong\mathbb A^{m_i+m_{i+1}'}\times Y/U_{i+1}$. 
                    \end{itemize}
                    
                    Apply Proposition \ref{proposition of induction for affine CDRS} to the action $U/U_i\rtimes \lambda\curvearrowright Y/U_i$ to conclude: 
                    \begin{itemize}
                        \item[(4')] The action $(U/U_{i+1})\rtimes \lambda\curvearrowright Y/U_{i+1}$ satisfies Condition \ref{condition of affine CDRS by semi-gradings}. 
                    \end{itemize}
                    
                    This completes the induction, and the proof. 
                \end{proof}
            }
        }
    
    }
    
    \subsection{U-hat theory for projective schemes}
    {
        Let $\lambda:\mathbb G_m\to \mathrm{Aut}(U)$ be a grading 1PS on $U$ with weights $w_1>\cdots>w_n$ on $\mathfrak u$. Let 
        \begin{equation}
            0=\mathfrak u_0\subseteq\mathfrak u_1\subseteq\cdots\subseteq\mathfrak u_n=\mathfrak u,\quad \mathfrak u_i:=\mathfrak u_{\lambda\geq w_i}
        \end{equation}
        be the corresponding filtration of Lie ideals of $\mathfrak u$ by $\lambda$-weights. Let $U_i\subseteq \hat U$ be the normal subgroup whose Lie algebra is $\mathfrak u_i$. 
    
        Let $X$ be a projective scheme over $\Bbbk$, and let $L$ be an ample line bundle over $X$. Let $\hat U $ act linearly on $X$ with respect to $L$. Then $\bigoplus_{d\geq 0}H^0(X,L^d)$ is finitely generated (\cite[\href{https://stacks.math.columbia.edu/tag/0B5T}{Lemma 0B5T}]{stacks-project}), and $\hat U $ acts on each $H^0(X,L^d)$. 
        
        For $d\in\mathbb N_+$ and $w\in \mathbb Z$, let $H^0(X,L^d)_{\lambda=w}\subseteq H^0(X,L^d)$ be the weight space for $\lambda$ of weight $w$. When $w\in\mathbb Z$ is not a weight, we can define the subspace $H^0(X,L^d)_{\lambda=w}$ as zero. Let $H^0(X,L^d)_{\lambda=\max}$ denote the weight space for $\lambda$ of the maximal weight on $H^0(X,L^d)$. Let\\ $H^0(X,L^d)_{\lambda<\max}$ denote the sum of weight spaces for $\lambda$ of weights that are not the maximal weight on $H^0(X,L^d)$. Then $H^0(X,L^d)=H^0(X,L^d)_{\lambda=\max}\oplus H^0(X,L^d)_{\lambda<\max}$, and $H^0(X,L^d)_{\lambda=\max}\ne0$. 
    
        Let $X^0_{\lambda-\min}\subseteq X$ be the open subscheme where $\bigoplus_{d>0}H^0(X,L^d)_{\lambda=\max}$ does not vanish. For the action of $\lambda$ on $X^0_{\lambda-\min}$, let $Z_{\lambda-\min}\hookrightarrow X^0_{\lambda-\min}$ be the \emph{fixed point subscheme} in the sense of Fogarty (\cite{fogarty1973fixed}). It is easy to see that $Z_{\lambda-\min}\hookrightarrow X$ is the closed subscheme associated to the ideal $\bigoplus_{d>0}H^0(X,L^d)_{\lambda<\max}$. Let $m\in\mathbb N_+$ be such that $H^0(X,L^m)$ generates $\bigoplus_{d\geq0}H^0(X,L^{md})$. Then $X^0_{\lambda-\min}\subseteq X$ has an affine open covering 
        \begin{equation}
            X^0_{\lambda-\min}=\bigcup X_f,\quad f\in H^0(X,L^m)_{\lambda=\max}. 
        \end{equation}
        
        {
            \begin{remark}\label{remark for empty X^0_min}
                If all sections $\bigoplus_{d>0}H^0(X,L^d)_{\lambda=\max}$ are nilpotent, then $X^0_{\lambda-\min}=\emptyset$ by definition. We can assume there exists a section $f\in H^0(X,L^d)_{\lambda=\max}$ for some $d>0$ and $f$ is not nilpotent for $X^0_{\lambda-\min}\ne\emptyset$, or we regard statements as vacuous truth when $X^0_{\lambda-\min}=\emptyset$. 
            \end{remark}
        }
    
        Since $\lambda$ is a grading 1PS on $U$, sections of maximal $\lambda$-weights are $U$-invariant, i.e. for all $d\geq0$ we have $H^0(X,L^d)_{\lambda=\max}\subseteq H^0(X,L^d)^U$ . Then $X_f$ is $\hat U $-invariant for $f\in H^0(X,L^m)_{\lambda=\max}$. 
    
        The ring of functions on $X_f$ is $\mathcal O(X_f):=\Big(\bigoplus_{d\geq0}H^0(X,L^{md})\Big)_{(f)}$. There is an induced action of $\lambda$ on $\mathcal O(X_f)$. Since $f$ has maximal $\lambda$-weight in $H^0(X,L^m)$, the $\lambda$-weights on $\mathcal O(X_f)$ are all non-positive, i.e. 
        \begin{equation}
            \mathcal O(X_f)=\bigoplus_{w\in\mathbb Z_{\leq 0}}\mathcal O(X_f)_{\lambda=w}. 
        \end{equation}
        The closed subscheme $Z_{\lambda-\min,f}\hookrightarrow X_f$ corresponds to $\mathcal O(X_f)\to\mathcal O(X_f)/\mathcal O(X_f)_{\lambda<0}$. Then\\ $\mathcal O(Z_{\lambda-\min,f})\cong\mathcal O(X_f)/\mathcal O(X_f)_{\lambda<0}\cong\mathcal O(X_f)_{\lambda=0}$. The closed immersion $Z_{\lambda-\min,f}\\\hookrightarrow X_f$ has a left inverse $X_f\to Z_{\lambda-\min,f}$, corresponding to $\mathcal O(X_f)_{\lambda=0}\to \mathcal O(X_f)$. They glue over $X_f$ to a morphism $p:X^0_{\lambda-\min}\to Z_{\lambda-\min}$. On $\Bbbk$-points, the morphism $p$ sends $x\in X^0_{\lambda-\min}$ to $p(x):=\lim_{t\to 0}\lambda(t).x$.

        \subsubsection{The CDRS condition}
        {
            The condition Constant-Dimension-of-Relative-Stabilisers (CDRS) is a condition for the action $\hat U \curvearrowright(X,L)$. 
            
            Let $\mathcal K_i:=\ker\big(\Omega_{X^0_{\lambda-\min}}\to\mathfrak u_i^*\otimes_\Bbbk\mathcal O_{X^0_{\lambda-\min}}\big)$ be the kernel of the infinitesimal action of $\mathfrak u_i$ on $X^0_{\lambda-\min}$ for all $0\leq i\leq n$. There exists a morphism for each $1\leq i\leq n$ 
            \begin{equation}\label{equation of phi_i from K_(i-1) to relative free vector bundle}
                \phi_i:\mathcal K_{i-1}\to (\mathfrak u_i/\mathfrak u_{i-1})^*\otimes\mathcal O_{X^0_{\lambda-\min}}. 
            \end{equation}

            {
                We will call the following condition Constant-Dimension-of-Relative-Stabilisers. 
                \begin{condition}\label{condition of CDRS for projective}
                    Let $\hat U \curvearrowright(X,L)$. Let $\lambda:\mathbb G_m\to \mathrm{Aut}(U)$ be a grading 1PS on $\mathfrak u$ with weights $w_1>\cdots>w_n$ on $\mathfrak u$. Let $\{\mathfrak u_i\}_{i=0}^n$ be the filtration $\mathfrak u_i:=\mathfrak u_{\lambda\geq w_i}$. This condition for the action $\hat U \curvearrowright(X,L)$ assumes that $\mathrm{coker}(\phi_i)$ for $\phi_i$ in equation \eqref{equation of phi_i from K_(i-1) to relative free vector bundle} is locally free of constant rank on $X^0_{\lambda-\min}$ for all $1\leq i\leq n$. 
                \end{condition}
                \begin{remark}
                    By Corollary \ref{corollary of flatness of coker(phi_i), Q(g_i) and dim of stab when reduced}, Condition \ref{condition of CDRS for projective} implies that $\dim\mathrm{Stab}_{\mathfrak u_i}(x)$ are constant on $X^0_{\lambda-\min}$. Moreover, the converse is true if $X^0_{\lambda-\min}$ is reduced. 
                \end{remark}
            }
            
            {
                \begin{lemma}\label{lemma of equivalence of coker being locally free and dim stab is constant when reduced}
                    Let $\hat U \curvearrowright (X,L)$. Let $\{\mathfrak u_i\}_{i=0}^n$ be a filtration of $\mathfrak u$. Let $\lambda:\mathbb G_m\to \mathrm{Aut}(U)$ be a grading on $U$. Let $1\leq i\leq n$ and $(r_1,\cdots,r_n)\in\mathbb N^n$. When $X$ is reduced, the following are equivalent: 
                    \begin{itemize}
                        \item[(1)] $\mathrm{coker}(\phi_i)$ for $\phi_i$ in equation $\eqref{equation of phi_i from K_(i-1) to relative free vector bundle}$ is locally free of rank $r_i$ on $X^0_{\lambda-\min}$ for each $1\leq i\leq n$; 
                        \item[(2)] $\dim\mathrm{Stab}_{(\mathfrak u_i:\mathfrak u_{i-1})}(x)=r_i$ for all $\Bbbk$-points $x\in X^0_{\lambda-\min}$ and for each $1\leq i\leq n$. 
                    \end{itemize}
                \end{lemma}
                \begin{proof}
                    This is from Corollary \ref{corollary of flatness of coker(phi_i), Q(g_i) and dim of stab when reduced}. 
                \end{proof}
            }
    
        }
        
        \subsubsection{The Grassmannian}
        {
            Let $W$ be a finite dimensional vector space over $\Bbbk$. Let $r:=\dim W$ and let $0\leq k\leq r$ be an integer. Let $\mathrm{Sch}/\Bbbk$ be the category of schemes of finite type over $\Bbbk$. 
            
            Consider the functor 
            \begin{equation}
                \mathbf{Grass}(W,k):(\mathrm{Sch}/\Bbbk)^{\mathrm{op}}\to\mathrm{Set}
            \end{equation}
            which associates to a scheme $S\in\mathrm{Sch}/\Bbbk$ the set $\mathbf{Grass}(W,k)(S)$ of isomorphism classes of surjections 
            \begin{equation}
                q:W\otimes_{\Bbbk}\mathcal O_S\to \mathcal Q
            \end{equation}
            where $\mathcal Q$ is a locally free sheaf of rank $k$ on $S$. Given $f:T\to S$ a morphism, let $\mathbf{Grass}(W,k)(f):\mathbf{Grass}(W,k)(S)\to \mathbf{Grass}(W,k)(T)$ be the map sending the class of $q:W\otimes\mathcal O_S\to\mathcal Q$ to the class of $f^*q:W\otimes\mathcal O_T\to f^*\mathcal Q$. 
            
            It is well known (cf. \cite[\href{https://stacks.math.columbia.edu/tag/089R}{Section 089R}]{stacks-project}) that the functor above is represented by a projective scheme of finite type over $\Bbbk$. Denote the representative scheme by $\mathrm{Grass}(W,k)$. 
            
            The closed points of $\mathrm{Grass}(W,k)$ are surjective linear maps $W\to Q$ with $\dim Q=k$. Let $V\subseteq W$ be a subspace of dimension $k$. Consider the subset $\mathrm{Grass}(W,k)_V\subseteq \mathrm{Grass}(W,k)$ of quotient spaces $W\to Q$ such that $V\subseteq W\to Q$ induces an isomorphism $V\cong Q$. Then $\mathrm{Grass}(W,k)_V$ is an open subset. The open subscheme $\mathrm{Grass}(W,k)_V\subseteq\mathrm{Grass}(W,k)$ defines an open sub-functor 
            \begin{equation}
                \mathbf{Grass}(W,k)_V:(\mathrm{Sch}/\Bbbk)^{\mathrm{op}}\to\mathrm{Set}
            \end{equation}
            which associates $S\in\mathrm{Sch}/\Bbbk$ to the subset $\mathbf{Grass}(W,k)_V(S)\subseteq\mathbf{Grass}(W,k)(S)$ of classes of surjections $q:W\otimes \mathcal O_S\to \mathcal Q$ such that the composition $V\otimes\mathcal O_S\to W\otimes\mathcal O_S\to \mathcal Q$ is an isomorphism. 
            
            {
                \begin{lemma}\label{lemma of equivariant covering of fixed point scheme of Grassmannian}
                    Let $G$ be a reductive group over $\Bbbk$. Let $W$ be a representation of $G$, inducing an action of $G$ on $\mathrm{Grass}(W,k)$. Let $\mathrm{Grass}(W,k)^G\hookrightarrow\mathrm{Grass}(W,k)$ denote the fixed point subscheme. Then there exist finitely many sub-representations $V_i\subseteq W$ with $\dim V_i=k$ such that $\mathrm{Grass}(W,k)^G\subseteq\bigcup_{i}\mathrm{Grass}(W,k)_{V_i}$. 
                \end{lemma}
                \begin{proof}
                    Consider closed points, i.e. $\Bbbk$-points. A $\Bbbk$-point in $\mathrm{Grass}(W,k)$ is a quotient $q:W\to Q$ with $\dim Q=k$. It is fixed by $G$ if and only if $\ker(q)\subseteq W$ is invariant under $G$. Since $G$ is reductive, there exists a sub-representation $V\subseteq W$ and $V\oplus \ker(q)=W$. Then $\dim V=k$ and $[q]\in\mathrm{Grass}(W,k)_V$. We only need a finite number of such $V_i$ since $\mathrm{Grass}(W,k)^G$ is quasi-compact, as a subscheme of a Noetherian scheme. 
                \end{proof}
            }
        }
    
        \subsubsection{Statement and proof of the first $\hat U$-theorem}
        {
            In this subsubsection, we will state and prove our first $\hat U$-theorem, Theorem \ref{theorem of universal geometric quotient by Un with proj CDRS}. 
            {
                \begin{lemma}\label{lemma of a lambda invariant open contains X_f if it contains Z_f}
                    For any $m\in\mathbb N_+$, let $f\in H^0(X,L^m)_{\lambda=\max}$, and let $Y\subseteq X^0_{\lambda-\min}$ be an open subset. If $Y$ is $\lambda$-invariant and $Z_{\lambda-\min,f}\subseteq Y$, then $X_f\subseteq Y$. 
                \end{lemma}
                \begin{proof}
                    The open affine $X_f$ is $\lambda$-invariant. Let $I\subseteq\mathcal O(X_f)$ be the radical ideal cutting out $X_f\setminus Y$. Then $I$ is $\lambda$-invariant. Since $X_f\setminus Y$ and $Z_{\lambda-\min,f}$ are disjoint, we have $\mathcal O(X_f)=I+\mathcal O(X_f)_{\lambda<0}$. There exists $g\in I$ and $h\in\mathcal O(X_f)_{\lambda<0}$ such that $1=g+h$. The weight zero component of $g$ is $g_{\lambda=0}=1$. Since $I$ is $\lambda$-invariant and $g\in I$, we have $1=g_{\lambda=0}\in I$. So $X_f\setminus Y=\emptyset$. 
                \end{proof}
            }
            
            {
                \begin{lemma}\label{lemma of proj CDRS implies covering of affine CDRS}
                    Assume Condition \ref{condition of CDRS for projective} for the action $\hat U \curvearrowright(X,L)$ with the grading 1PS $\lambda:\mathbb G_m\to \mathrm{Aut}(U)$. Then there exists $m\in\mathbb N_+$ and $f_1,\cdots,f_r\in H^0(X,L^m)_{\lambda=\max}$ such that: 
                    \begin{itemize}
                        \item $X^0_{\lambda-\min}=\bigcup_{k=1}^rX_{f_k}$ is a covering by $\hat U $-invariant affine open subsets; 
                        \item Condition \ref{condition of affine CDRS by semi-gradings} is satisfied for the action $\hat U \curvearrowright X_{f_k}$ for $k=1,\cdots,r$. 
                    \end{itemize}
                \end{lemma}
                \begin{proof}
                    For $1\leq i\leq n$, let $k_i$ be the rank of $\mathrm{coker}(\phi_i)$, and let $W_i:=(\mathfrak u_i/\mathfrak u_{i-1})^*$. Then the surjective morphism 
                    \begin{equation}
                        q_i:W_i\otimes\mathcal O_{X^0_{\lambda-\min}}\to \mathrm{coker}(\phi_i)
                    \end{equation}
                    defines an isomorphism class in $\mathrm{Grass}(W_i,k_i)(X^0_{\lambda-\min})$, corresponding to a morphism 
                    \begin{equation}
                        \Theta_i:X^0_{\lambda-\min}\to \mathrm{Grass}(W_i,k_i). 
                    \end{equation}
                    The product over $\Bbbk$ of the morphisms over $1\leq i\leq n$ is a morphism 
                    \begin{equation}
                        \Theta:X^0_{\lambda-\min}\to \prod_{i=1}^n\mathrm{Grass}(W_i,k_i). 
                    \end{equation}
                    
                    Let $\iota:Z_{\lambda-\min}\hookrightarrow X^0_{\lambda-\min}$. The morphism $q_i$ is $\lambda$-equivariant, so $\Theta_i$ and $\Theta_i\circ\iota$ are $\lambda$-equivariant. In particular $\Theta_i\circ\iota$ is $\lambda$-equivariant. Since $\lambda$ acts trivially on $Z_{\lambda-\min}$, we have that $\Theta_i\circ\iota$ factors through the fixed point subscheme $\mathrm{Grass}(W_i,k_i)^\lambda\hookrightarrow\mathrm{Grass}(W_i,k_i)$. 
                    
                    By Lemma \ref{lemma of equivariant covering of fixed point scheme of Grassmannian}, there exist a finite number of $\lambda$-invariant subspaces $V_{ij}\subseteq W_i$ for $j=1,\cdots,m_i$ such that $\mathrm{Grass}(W_i,k_i)^\lambda\subseteq\bigcup_{j=1}^{m_i}\mathrm{Grass}(W_i,k_i)_{V_{ij}}$. Then 
                    \begin{equation}
                        \prod_{i=1}^n\mathrm{Grass}(W_i,k_i)^\lambda\subseteq\bigcup_{j_1,\cdots,j_n}\prod_{i=1}^n\mathrm{Grass}(W_i,k_i)_{V_{i,j_i}}. 
                    \end{equation}
                    
                    We then have a covering of $Z_{\lambda-\min}$ by open subsets invariant under $\lambda$: 
                    \begin{equation}
                        \begin{split}
                            Z_{\lambda-\min}=&\bigcup_{j_1,\cdots,j_n}\big(\Theta\circ\iota\big)^{-1}\Big(\prod_{i=1}^n\mathrm{Grass}(W_i,k_i)_{V_{i,j_i}}\Big)\\
                            =&\bigcup_{j_1,\cdots,j_n}Z_{\lambda-\min}\cap \Theta^{-1}\Big(\prod_{i=1}^n\mathrm{Grass}(W_i,k_i)_{V_{i,j_i}}\Big). 
                        \end{split}
                    \end{equation}
                    Since $V_{ij}\subseteq W_i$ is $\lambda$-invariant and $\Theta_i$ is $\lambda$-equivariant for each $1\leq i\leq n$ and $1\leq j\leq m_i$, each preimage $\Theta^{-1}\Big(\prod_{i=1}^n\mathrm{Grass}(W_i,k_i)_{V_{i,j_i}}\Big)$ is $\lambda$-invariant open. 
                    
                    This covering can be reduced and refined to a finite covering 
                    \begin{equation}
                        Z_{\lambda-\min}=\bigcup_{k=1}^r Z_{\lambda-\min,f_k}
                    \end{equation}
                    such that $f_1,\cdots,f_r\in H^0(X,L^m)_{\lambda=\max}$ for some $m\in\mathbb N_+$ and each $Z_{\lambda-\min,f_k}$ is contained in some $\Theta^{-1}\Big(\prod_{i=1}^n\mathrm{Grass}(W_i,k_i)_{V_{i,j_i}}\Big)$. The following is an open affine covering of $X^0_{\lambda-\min}$
                    \begin{equation}
                        X^0_{\lambda-\min}=\bigcup_{k=1}^rX_{f_k}
                    \end{equation}
                    and by Lemma \ref{lemma of a lambda invariant open contains X_f if it contains Z_f}, each $X_{f_k}$ is contained in some $\Theta^{-1}\Big(\prod_{i=1}^n\mathrm{Grass}(W_i,k_i)_{V_{i,j_i}}\Big)$. 
                    
                    Let $f$ be one of $f_1,\cdots,f_r$ and assume $X_f\subseteq \Theta^{-1}\Big(\prod_{i=1}^n\mathrm{Grass}(W_i,k_i)_{V_{i,j_i}}\Big)$ for some $j_1,\cdots,j_n$. Then $X_f$ is $\hat U $-invariant and we show that Condition \ref{condition of affine CDRS by semi-gradings} is satisfied for the action $\hat U \curvearrowright X_f$. 
                    
                    The inclusion $\Theta_i(X_f)\subseteq\mathrm{Grass}(W_i,k_i)_{V_{i,j_i}}$ is equivalent to the following diagram 
                    \begin{equation}
                        \begin{tikzcd}
                            &V_{i,j_i}\otimes\mathcal O(X_f)\ar[d]\ar[rd,"\cong"]\\
                            \mathcal K_{i-1}(X_f)\ar[r,"\phi_i(X_f)"]&W_i\otimes\mathcal O(X_f)\ar[r]&\mathrm{coker}(\phi_i)(X_f)\ar[r]&0. 
                        \end{tikzcd}
                    \end{equation}
                    The quotient map $\mathfrak u_i/\mathfrak u_{i-1}\to\mathfrak s_i$ in Condition \ref{condition of affine CDRS by semi-gradings} can be the quotient map $\mathfrak u_i/\mathfrak u_{i-1}\to (V_{i,j_i})^*$. 
    
                    These together show Condition \ref{condition of affine CDRS by semi-gradings} for the action $\hat U \curvearrowright X_f$. 
                \end{proof}
            }
    
            {
                \begin{theorem}\label{theorem of universal geometric quotient by Un with proj CDRS}
    
                    Let $\lambda:\mathbb G_m\to \mathrm{Aut}(U)$ be a grading 1PS on $U$ with weights $w_1>\cdots>w_n$. Let $\{\mathfrak u_i\}_{i=0}^n$ be the filtration $\mathfrak u_i:=\mathfrak u_{\lambda\geq w_i}$. Let $X$ be a projective scheme over $\Bbbk$ and let $L$ be an ample line bundle over $X$. Let $\hat U $ act on $X$ and linearise $L$. Assume that the CDRS condition (Condition \ref{condition of CDRS for projective}) holds for the action $\hat U \curvearrowright(X,L)$. Then the universal geometric quotient of $X^0_{\lambda-\min}$ by $U$ exists 
                    \begin{equation}
                        \phi:X^0_{\lambda-\min}\to X^0_{\lambda-\min}/U
                    \end{equation}
                    and $X^0_{\lambda-\min}/U$ is quasi-projective of finite type over $\Bbbk$. 
                    
                    Moreover, there exists $N\in\mathbb N_+$ such that 
                    \begin{itemize}
                        \item the linear system $W:=H^0(X,L^N)^U$ is base point free on $X^0_{\lambda-\min}$; 
                        \item the morphism $X^0_{\lambda-\min}\to \mathbb P(W^*)$ factors through a closed immersion 
                        \begin{equation}
                            X^0_{\lambda-\min}/U\hookrightarrow \mathbb P(W^*)^0_{\lambda-\min}
                        \end{equation}
                        where $\mathbb P(W^*)^0_{\lambda-\min}\subseteq\mathbb P(W^*)$ is the non-vanishing locus of $H^0(X,L^N)_{\lambda=\max}$ 
                        \begin{equation}
                            \begin{tikzcd}
                                X^0_{\lambda-\min}\ar[r]\ar[d,"\phi"]&\mathbb P(W^*)\\
                                X^0_{\lambda-\min}/U\ar[r,hook]&\mathbb P(W^*)^0_{\lambda-\min}\ar[u,"\subseteq"{sloped}]. 
                            \end{tikzcd}
                        \end{equation}
                    \end{itemize}
                \end{theorem}
                \begin{proof}
                    Condition \ref{condition of CDRS for projective} is unchanged if $L$ is replaced by a tensor power. Assume without loss of generality that $H^0(X,L)$ generates $\bigoplus_{d\geq0}H^0(X,L^d)$. 
    
                    By Lemma \ref{lemma of proj CDRS implies covering of affine CDRS}, there exists $m\in\mathbb N_+$ and non-zero sections 
                    \begin{equation}
                        f_1,\cdots,f_r\in H^0(X,L^m)_{\lambda=\max}
                    \end{equation}
                    satisfying: 
                    \begin{itemize}
                        \item $X^0_{\lambda-\min}=\bigcup_{k=1}^rX_{f_k}$; 
                        \item the action $\hat U \curvearrowright X_{f_k}$ satisfies Condition \ref{condition of affine CDRS by semi-gradings} for $k=1,\cdots,r$. 
                    \end{itemize}
                    
                    By Theorem \ref{theorem of affine quotients for affine CDRS length one filtration}, each $X_{f_k}$ has a universal geometric quotient by $U$ 
                    \begin{equation}
                        \phi_k:X_{f_k}\to Y_k:=X_{f_k}/U
                    \end{equation}
                    and each $Y_k$ is an affine scheme of finite type over $\Bbbk$. 
                    
                    The function $f_{k'}/f_k\in\mathcal O(X_{f_k})$ is $U$-invariant, so it corresponds to a function $g_{k'k}\in\mathcal O(Y_k)\cong\mathcal O(X_{f_k})^U$. The pull back of $\phi_k$ along the open immersion $Y_k\setminus\mathbb V(g_{k'k})\to Y_k$ is also a geometric quotient, i.e. 
                    \begin{equation}
                        X_{f_k}\setminus\mathbb V(f_{k'}/f_k)\to Y_k\setminus\mathbb V(g_{k'k})
                    \end{equation}
                    is a geometric quotient. Since $X_{f_k}\setminus \mathbb V(f_{k'}/f_k)=X_{f_{k'}}\setminus\mathbb V(f_k/f_{k'})$, their quotients by $U$ are isomorphic by a unique isomorphism 
                    \begin{equation}
                        Y_k\setminus\mathbb V(g_{k'k})\cong Y_{k'}\setminus\mathbb V(g_{kk'}). 
                    \end{equation}
                    
                    The affine schemes $\{Y_k\}$ then glue to a scheme $Y$, and the morphisms $\phi_k:X_{f_k}\to Y_k$ glue to a morphism $\phi:X^0_{\lambda-\min}\to Y$, which is a universal geometric quotient by $U$. Moreover $Y$ is a scheme of finite type over $\Bbbk$. \bigskip
    
                    Since $Z_{\lambda-\min}=\bigcup_{k=1}^rZ_{\lambda-\min,f_k}$, for all $m'\gg1$ we have 
                    \begin{equation}\label{equation of f1 ... fk eventually generate homogeneous coordinate ring of Z_min}
                        H^0(X,L^{mm'})_{\lambda=\max}=\sum_{k=1}^rf_kH^0(X,L^{m(m'-1)})_{\lambda=\max}. 
                    \end{equation}
    
                    Let $1\leq k\leq r$. Since $\mathcal O(X_{f_k})^U$ is finitely generated over $\Bbbk$, there exists $m'_k\in\mathbb N_+$ such that $\frac{1}{f_k^{m'_k}}H^0(X,L^{mm'_k})^U$ generates $\mathcal O(X_{f_k})^U$. We can enlarge $m'_k$ to $m'\in\mathbb N_+$ uniform for all $1\leq k\leq r$ and further assume equation \eqref{equation of f1 ... fk eventually generate homogeneous coordinate ring of Z_min}, i.e. $H^0(X,L^{mm'})_{\lambda=\max}\subseteq\langle f_1,\cdots,f_r\rangle$. The generating subspace $\frac{1}{f_k^{m'}}H^0(X,L^{mm'})^U\subseteq\mathcal O(X_{f_k})^U$ induces a surjective ring map 
                    \begin{equation}
                        \mathrm{Sym}^\bullet(W)_{(f_k^{m'})}\to \mathcal O(X_{f_k})^U
                    \end{equation}
                    where $W:=H^0(X,L^{mm'})^U$ and the symmetric product is over $\Bbbk$. 
                    
                    The linear system $W=H^0(X,L^{mm'})^U$ is base point free on $X^0_{\lambda-\min}$ since 
                    \begin{equation}
                        \mathbb V(W)\cap X^0_{\lambda-\min}\subseteq \mathbb V\big(H^0(X,L^{mm'})_{\lambda=\max}\big)\cap X^0_{\lambda-\min}=\emptyset. 
                    \end{equation}
                    Then there is a morphism $X^0_{\lambda-\min}\to \mathbb P(W^*)$, which factors through $X^0_{\lambda-\min}/U\to \mathbb P(W^*)$, since it is $U$-invariant
                    \begin{equation}
                        \begin{tikzcd}
                            X^0_{\lambda-\min}\ar[r]\ar[d,"\phi"]&\mathbb P(W^*)\\
                            X^0_{\lambda-\min}/U\ar[ru]. 
                        \end{tikzcd}
                    \end{equation}

                    We then prove that $X^0_{\lambda-\min}/U\hookrightarrow\mathbb P(W*)^0_{\lambda-\min}$ is a closed immersion. 
                    
                    Observe that for any $g\in H^0(X,L^{m(m'-1)})_{\lambda=\max}$, the open affine subscheme $X_{f_kg}\subseteq X_{f_k}$ is $\hat U $-invariant and Condition \ref{condition of affine CDRS by semi-gradings} is satisfied for the action $\hat U \curvearrowright X_{f_kg}$. Therefore $X_{f_kg}\to X_{f_kg}/U$ is a universal geometric quotient by $U$ and $\mathcal O(X_{f_kg})^U$ is finitely generated by Theorem \ref{theorem of affine U-quotient for CDRS}. We show that $\mathcal O(X_{f_kg})^U$ is generated by $\frac{1}{f_kg}W$ for any $g\in H^0(X,L^{m(m'-1)})_{\lambda=\max}$. Any element in $\mathcal O(X_{f_kg})$ is of the form $\frac{h}{(f_kg)^s}$ for $s\in\mathbb N_+$ and $h\in H^0(X,L^{mm's})$. We can enlarge $s\in\mathbb N_+$ such that $(f_kg)^s$ is $U$-invariant. Assume $\frac{h}{(f_kg)^s}$ is $U$-invariant. Then $\xi.h\in H^0(X,L^{mm's})$ is annihilated by a power of $f_kg$ for all $\xi\in\mathfrak u$. We can test on a basis of $\mathfrak u$, and choose the power of $f_kg$ uniformly. There exists $s'\in\mathbb N_+$ such that $(f_kg)^{s'}h\in H^0(X,L^{mm'(s+s')})^U$. Then $\frac{(f_kg)^{s'}h}{f_k^{m'(s+s')}}\in\mathcal O(X_{f_k})^U$, which is in the image of $\mathrm{Sym}^\bullet(W)_{(f_k^{m'})}$, i.e. 
                    \begin{equation}
                        \frac{(f_kg)^{s'}h}{f_k^{m'(s+s')}}=\frac{\sum_{i_1,\cdots,i_t}w_{i_1}\cdots w_{i_t}}{(f_k^{m'})^t}
                    \end{equation}
                    for some $t\in\mathbb N_+$ and $w_{i_1},\cdots,w_{i_t}\in W$. We also have 
                    \begin{equation}
                        g^{m'}\in H^0(X.L^{mm'(m'-1)})_{\lambda=\max}=\prod^{m'-1}H^0(X,L^{mm'})_{\lambda=\max}\subseteq\prod^{m'-1}W
                    \end{equation}
                    so $g^{m'}=\sum_{j_1,\cdots,j_{m'-1}}w'_{j_1}\cdots w'_{j_{m'-1}}$ for some $w'_{j_1},\cdots, w'_{j_{m'-1}}\in W$. 
                    
                    Then in $\mathcal O(X_{f_kg})^U$ 
                    \begin{equation}
                        \begin{split}
                            \frac{h}{(f_kg)^s}=&\Big(\frac{f_k^{m'}}{f_kg}\Big)^{s+s'}\frac{(f_kg)^{s'}h}{f_k^{m'(s+s')}}\\
                            =&\Big(\frac{f_k^{m'}}{f_kg}\Big)^{s+s'}\frac{\sum_{i_1,\cdots,i_t}w_{i_1}\cdots w_{i_t}}{(f_k^{m'})^t}\\
                            =&\Big(\frac{f_k^{m'}}{f_kg}\Big)^{s+s'}\frac{\big(\sum_{i_1,\cdots,i_t}w_{i_1}\cdots w_{i_t}\big)g^{m't}}{(f_kg)^{m't}}\\
                            =&\Big(\frac{f_k^{m'}}{f_kg}\Big)^{s+s'}\frac{\big(\sum_{i_1,\cdots,i_t}w_{i_1}\cdots w_{i_t}\big)\big(\sum_{j_1,\cdots,j_{m'-1}}w'_{j_1}\cdots w'_{j_{m'-1}}\big)^t}{(f_kg)^{m't}}
                        \end{split}
                    \end{equation}
                    which is generated by $\frac{1}{f_kg}W$. This completes the proof that $\frac{1}{f_kg}W\subseteq\mathcal O(X_{f_kg})^U$ is a generating subspace, i.e. the following ring map is surjective 
                    \begin{equation}
                        \mathrm{Sym}^\bullet(W)_{(f_kg)}\to \mathcal O(X_{f_kg})^U. 
                    \end{equation}
                    
                    The surjective ring map above is equivalent to the closed immersion 
                    \begin{equation}
                        X_{f_kg}/U\hookrightarrow \mathbb P(W^*)_{f_kg}. 
                    \end{equation}
                    
                    The union of $\mathbb P(W^*)_{f_kg}$ over $k=1,\cdots,r$ and $g\in H^0(X,L^{m(m
                    -1)})_{\lambda=\max}$ is the non-vanishing locus of $H^0(X,L^{mm'})_{\lambda=\max}\subseteq W$, since we have assumed equation \eqref{equation of f1 ... fk eventually generate homogeneous coordinate ring of Z_min}. The closed immersions above then glue over $1\leq k\leq r$ to a closed immersion 
                    \begin{equation}
                        X^0_{\lambda-\min}/U\hookrightarrow\mathbb P(W^*)\setminus\mathbb V\big(H^0(X,L^{mm'})_{\lambda=\max}\big)=\mathbb P(W^*)^0_{\lambda-\min}. 
                    \end{equation}
                \end{proof}
            }
    
        }
    
    }

    \subsection{Quotienting by \texorpdfstring{$\hat U$}{U-hat}}
    {
        The residual linearisation of $\lambda$ can be twisted to get a projective geometric quotient of $X$ by $\hat U$ when the CDRS condition (Condition \ref{condition of CDRS for projective}) holds. 
        {
            \begin{theorem}[\cite{berczi2018geometric}]\label{theorem of U-hat quotient with CDRS}
                Let $X$ be a projective scheme over $\Bbbk$ and let $L$ be an ample line bundle over $X$. Let $\hat U$ act on $X$ and linearise $L$. Assume Condition \ref{condition of CDRS for projective} for $\hat U\curvearrowright(X,L)$. Assume $X^0_{\lambda-\min}\ne\emptyset$. Then the following is $\hat U$-invariant open in $X^0_{\lambda-\min}$
                \begin{equation}
                    X^0_{\lambda-\min}\setminus UZ_{\lambda-\min}
                \end{equation}
                which admits a projective universal geometric quotient by $\hat U$. 
            \end{theorem}
            \begin{proof}
                With Condition \ref{condition of CDRS for projective}, by Theorem \ref{theorem of universal geometric quotient by Un with proj CDRS}, the universal geometric quotient of $X^0_{\lambda-\min}$ by $U$ exists 
                \begin{equation}
                    \phi:X^0_{\lambda-\min}\to X^0_{\lambda-\min}/U
                \end{equation}
                and there exists $N\in\mathbb N_+$ such that for $W:=H^0(X,L^N)^{U}$, there is a locally closed immersion
                \begin{equation}
                    \varphi:X^0_{\lambda-\min}/U\to\mathbb P(W^*).
                \end{equation}
                Let $\overline{\varphi}:\overline{X^0_{\lambda-\min}/U}\hookrightarrow \mathbb P(W^*)$ be the scheme theoretic image of $\varphi$. By \cite[\href{https://stacks.math.columbia.edu/tag/01R8}{Lemma 01R8}]{stacks-project}, the underlying closed subset of $\overline{X^0_{\lambda-\min}/U}\subseteq\mathbb P(W^*)$ is the closure of $\varphi:X^0_{\lambda-\min}/U\to \mathbb P(W^*)$. Since $X^0_{\lambda-\min}/U\hookrightarrow\mathbb P(W^*)^0_{\lambda-\min}$ is closed, we have $X^0_{\lambda-\min}/U=\overline{X^0_{\lambda-\min}/U}\cap \mathbb P(W^*)^0_{\lambda-\min}$. 
                
                The action of $\hat U $ on $H^0(X,L^N)$ induces an action of $\lambda$ on $W$. Then $\lambda\curvearrowright\big(\mathbb P(W^*),\mathcal O(1)\big)$ and $\lambda\curvearrowright\big(\overline{X^0_{\lambda-\min}/U},\overline{\varphi}^*\mathcal O(1)\big)$. 
                
                By Theorem 1.19, \cite{mumford1994geometric}, for any ample linearisation $M$ of $\lambda$ on $\mathbb P(W^*)$ 
                \begin{equation}
                    \overline{X^0_{\lambda-\min}/U}^{\lambda,\overline{\varphi}^*M,\mathrm{(s)s}}=\overline{X^0_{\lambda-\min}/U}\cap \mathbb P(W^*)^{\lambda,M,\mathrm{(s)s}}. 
                \end{equation}
                The linearisation $M$ will be chosen as a twist of $\mathcal O(1)$, such that $\mathbb P(W^*)^{\lambda,M,\mathrm{(s)s}}\subseteq \mathbb P(W^*)^0_{\lambda-\min}$. 
                
                Characters of $\lambda$ are integers. Let $\mathrm{Hom}(\lambda,\mathbb G_m)$ be the Abelian group of characters. Call elements in $\mathrm{Hom}(\lambda,\mathbb G_m)\otimes_{\mathbb Z}\mathbb Q$ \emph{rational characters}. If $\theta\in\mathrm{Hom}(\lambda,\mathbb G_m)\otimes_{\mathbb Z}\mathbb Q$ is a rational character, then we denote the twisted linearisation by $\mathcal O(1)_\theta$. 
                
                Denote the maximal weight of $\lambda$ on $H^0(X,L)$ by $\max(\lambda)$. Denote the maximal and next maximal weights of $\lambda$ on $W$ by $\max(\lambda,W)$ and $\mathrm{2nd}\max(\lambda,W)$. Then $\max(\lambda,W)=N\max(\lambda)$. When $\epsilon>0$ is small, we have
                \begin{equation}
                    \mathrm{2nd}\max(\lambda,W)<-N\theta<\max(\lambda,W)
                \end{equation}
                then 
                \begin{equation}
                    \mathbb P(W^*)^{\lambda,\mathcal O(1)_{N\theta},\mathrm{ss}}=\mathbb P(W^*)^{\lambda,\mathcal O(1)_{N\theta},\mathrm{s}}=\mathbb P(W^*)^0_{\lambda-\min}. 
                \end{equation}

                For such $\theta$, we have 
                \begin{equation}
                    \overline{X^0_{\lambda-\min}/U}^{\lambda,L_\theta,\mathrm{(s)s}}= X^0_{\lambda-\min}/U\;\cap\;\Big(\mathbb P(W^*)^0_{\lambda-\min}\setminus \mathbb P\big((W_{\lambda=\max})^*\big)\Big). 
                \end{equation}
                The preimage of $\overline{X^0_{\lambda-\min}/U}^{\lambda,L_\theta,\mathrm{(s)s}}\subseteq X^0_{\lambda-\min}/U$ along $\phi:X^0_{\lambda-\min}\to X^0_{\lambda-\min}/U$ is 
                \begin{equation}
                    X^0_{\lambda-\min}\setminus UZ_{\lambda-\min}
                \end{equation}
                which is $\hat U$-invariant open. As a composition of universal geometric quotients the following is a universal geometric quotient 
                \begin{equation}
                    X^0_{\lambda-\min}\setminus UZ_{\lambda-\min}\to \big(X^0_{\lambda-\min}\setminus UZ_{\lambda-\min}\big)\big/\hat U. 
                \end{equation}
                
                For projectivity, observe that 
                \begin{equation}
                    \big(X^0_{\lambda-\min}\setminus UZ_{\lambda-\min}\big)\big/\hat U\cong \overline{X^0_{\lambda-\min}/U}^{\lambda,L_\theta,\mathrm{(s)s}}/\lambda=\overline{X^0_{\lambda-\min}/U}/\!/\lambda. 
                \end{equation}
            \end{proof}
        }
    }
}
\section{Blow-up}\label{section of blowup}
{
    Let $\hat U\curvearrowright(X,L)$ and $w_1>\cdots>w_n$ and $\{\mathfrak u_i\}_{i=0}^n$ be as above. In this section, we do not assume Condition \ref{condition of CDRS for projective}, and describe a blow-up procedure to achieve Condition \ref{condition of CDRS for projective} so that Theorem \ref{theorem of U-hat quotient with CDRS} applies to the induced action on the blow-up of $X$. A weaker condition is needed for this process (Condition \ref{condition before blow-up for CDRS}), and the main result is Theorem \ref{theorem of blow-up for CDRS}. 
    
    \subsection{Examples with few grading weights}
    {
        The general construction involves a tedious calculation (Lemma \ref{lemma of existence of b^(i)_nu}). It may help to motivate the construction there by first considering some more tractable examples. In this subsection, we consider the simple situations when there are one or two grading weights for $\lambda$ on $U$. Consider the situation: 
        \begin{itemize}
            \item $\lambda:\mathbb G_m\to \mathrm{Aut}(U)$ is a grading on $U$ with grading weights $w_1>\cdots>w_n$ for $n=1$ or $2$; 
            \item $Y=\mathrm{Spec}(A)$ is a \emph{reduced} affine scheme of finite type over $\Bbbk$; 
            \item $\hat U$ acts on $Y$ such that the action of $\lambda$ on $A$ has all weights non-positive (cf. the first requirement in Condition \ref{condition of affine CDRS by semi-gradings}). 
        \end{itemize}
        
        The dimensions of unipotent stabilisers are not assumed constant, and the blow-up is to achieve this. Consider the exact sequence 
        \begin{equation}
            \begin{tikzcd}
                K_{i-1}\ar[r,"\phi_i"]&(\mathfrak u_i/\mathfrak u_{i-1})^*\otimes A\ar[r]&\mathrm{coker}(\phi_i)\ar[r]&0
            \end{tikzcd}
        \end{equation}
        where $K_{i-1}=\ker\big(\Omega_{A/\Bbbk}\to \mathfrak u_{i-1}^*\otimes A\big)$. For simplicity, write $\mathrm{Fit}_k(\phi_i)$ for $\mathrm{Fit}_k(\mathrm{coker}(\phi_i))$. A suitable $\hat U$-equivariant blow-up $\tilde Y\to Y$ with an ample linearisation will satisfy the condition that on $\tilde Y^0_{\lambda-\min}$, the dimensions of $\mathfrak u_i$-stabilisers are constant for each $1\leq i\leq n$. By Corollary \ref{corollary of flatness of coker(phi_i), Q(g_i) and dim of stab when reduced}, the condition that $\mathfrak u_i$-stabilisers have constant dimensions is equivalent to that $\mathrm{coker}(\tilde\phi_i)$ are locally free of constant ranks. We will verify that $\mathrm{Fit}_{k_i-1}(\tilde\phi_i)=0$ and $\mathrm{Fit}_{k_i}(\tilde\phi_i)=\mathcal O_{\tilde Y^0_{\lambda-\min}}$ for each $1\leq i\leq n$ for some $k_1,\cdots,k_n\in\mathbb N$. This condition is equivalent to $\dim\mathrm{Stab}_{\mathfrak u_i}(y)=k_1+\cdots+k_i$ for $\Bbbk$-points $y\in \tilde Y^0_{\lambda-\min}$ and $1\leq i\leq n$ by Proposition \ref{proposition of properties of Fit} and Corollary \ref{corollary of flatness of coker(phi_i), Q(g_i) and dim of stab when reduced}. 
        
        Let $k_1,\cdots,k_n\in\mathbb N$ be such that 
        \begin{equation}
            k_i:=\min\{k\in\mathbb N:\mathrm{Fit}_k(\phi_i)\ne0\},\quad i=1,\cdots,n. 
        \end{equation}

        Recall that $Y=\mathrm{Spec}(A)$ is a reduced scheme of finite type over $\Bbbk$. If we assume $\prod_{i=1}^n\mathrm{Fit}_{k_i}(\phi_i)\ne0$, then the non-vanishing locus of $\prod_{i=1}^n\mathrm{Fit}_{k_i}(\phi_i)$ is non-empty. This non-vanishing locus equals the intersection of the non-vanishing loci of $\mathrm{Fit}_{k_i}(\phi_i)$ for $i=1,\cdots,n$. For $1\leq i\leq n$, the non-vanishing locus of $\mathrm{Fit}_{k_i}(\phi_i)$ is the largest open subscheme such that the restriction of $\mathrm{coker}(\phi_i)$ is locally free of rank $k_i$ by \cite[\href{https://stacks.math.columbia.edu/tag/05P8}{Lemma 05P8}]{stacks-project}. Therefore the non-vanishing locus of $\prod_{i=1}^n\mathrm{Fit}_{k_i}(\phi_i)$ is the largest open subscheme such that $\mathrm{coker}(\phi_i)$ is locally free of rank $k_i$ for each $i=1,\cdots,n$.  By Corollary \ref{corollary of flatness of coker(phi_i), Q(g_i) and dim of stab when reduced}, the closed points of the non-vanishing locus of $\prod_{i=1}^n\mathrm{Fit}_{k_i}(\phi_i)$ are exactly $y\in Y$ such that 
        \begin{equation}
            \dim\mathrm{Stab}_{\mathfrak u_i}(y)=k_1+\cdots+k_i,\quad i=1,\cdots,n. 
        \end{equation}
        
        \subsubsection{One grading weight}
        {
            Assume there is only one grading weight, i.e. $n=1$. Let $w>0$ be the grading weight. Let $r=\dim U$ and choose a basis $(\xi_1,\cdots,\xi_r)$ of $\mathfrak u$ with its dual basis $(u_1,\cdots,u_r)$ on $\mathfrak u^*$. Let $\phi:\Omega_{A/\Bbbk}\to\mathfrak u^*\otimes A$ be the morphism of coherent sheaves and denote by $\mathrm{Fit}_k(\phi)\subseteq A$ the $k$th Fitting ideal of $\mathrm{coker}(\phi)$ for $k\in\mathbb Z$. Let $k\in\mathbb N$ be such that $\mathrm{Fit}_k(\phi)\ne0$ and $\mathrm{Fit}_{k-1}(\phi)=0$. 
            
            We need an assumption on $\mathrm{Fit}_k(\phi)$, that is $\mathrm{Fit}_k(\phi)\not\subseteq A_{\lambda<0}$, which is equivalent to that $\dim\mathrm{Stab}_U(z)=k$ for some $\Bbbk$-points $z\in\mathrm{Spec}(A/A_{\lambda<0})$ (See Condition \ref{condition before blow-up for CDRS} and Lemma \ref{lemma of point-wise description of the nice locus of small unipotent stabilisers} below). Let $I\subseteq A$ be the ideal 
            \begin{equation}
                I:=\mathrm{Fit}_k(\phi)+A_{\lambda<0}=\mathrm{Fit}_k(\phi)_{\lambda=0}\oplus A_{\lambda<0}. 
            \end{equation}
            Let $J\subseteq A$ be the kernel of the composition $A\to\mathcal O(U)\otimes A\to \mathcal O(U)\otimes(A/I)$. The closed subscheme associated to $J\subseteq A$ is the scheme theoretic image of $U\times \mathbb V(I)\to Y$, where $\mathbb V(I)\hookrightarrow Y$ is the closed subscheme associated to $I\subseteq A$. 
            
            {
                \begin{remark}
                    The $\Bbbk$-points in $\mathbb V(I)$ are exactly $z\in \mathbb V(A_{<0})$ such that $\dim\mathrm{Stab}_U(z)>k$. By \cite[\href{https://stacks.math.columbia.edu/tag/01R8}{Lemma 01R8}]{stacks-project}, the closed subset $\mathbb V(J)$ is the closure of the $U$-sweep of $\mathbb V(I)$. 
                \end{remark}
            }
            
            We choose the centre $C\hookrightarrow Y$ of the blow-up as the closed subscheme associated to $J\subseteq A$. The extreme situation when $J=0$ does not happen, since the assumption $\mathrm{Fit}_k(\phi)\not\subseteq A_{\lambda<0}$ implies $0\ne\mathrm{Fit}_k(\phi)_{\lambda=0}\subseteq J$. 
            
            The ideal $J\subseteq A$ is $\hat U$-invariant and $J_{\lambda=\max}=J_{\lambda=0}=\mathrm{Fit}_k(\phi)_{\lambda=0}$. Let $\tilde Y\to Y$ be the blow-up. Then $\tilde Y$ is constructed by 
            \begin{equation}
                \begin{split}
                    \mathrm{Bl}_J(A):=&\bigoplus_{d\geq0}J^d=A\oplus J\oplus J^2\oplus\cdots\\
                    \tilde Y:=&\mathrm{Proj}\Big(\mathrm{Bl}_J(A)\Big)
                \end{split}
            \end{equation}
            where $\mathrm{Bl}_J(A)$ is called the \emph{blowup algebra associated to $J\subseteq A$}. The degree on $\mathrm{Bl}_J(A)$ is the natural one from the direct sum. Let $E\hookrightarrow\tilde Y$ be the exceptional divisor with its associated line bundle $\mathcal O(E)$. The line bundle $\mathcal O(1)$ from the $\mathrm{Proj}$-construction equals $\mathcal O(-E)$ and it is relatively very ample. Since $Y$ ia affine, we have that $\mathcal O(1)$ is ample. The blow-up is $\hat U$-equivariant, and then $\hat U$ acts on $\tilde Y$ linearly with respect to $\mathcal O(1)$. With this ample linearisation, we can define $(\tilde Y)^0_{\lambda-\min}$, which is the union of the following open affine subschemes 
            \begin{equation}
                \begin{split}
                    (\tilde Y)^0_{\lambda-\min}=&\bigcup_{a\in J_{\lambda=\max}}\tilde Y_a\\
                    \tilde Y_a:=&\mathrm{Spec}\Big(A\Big[\frac{J}{a}\Big]\Big)\\
                    A\Big[\frac{J}{a}\Big]:=&\Big(\mathrm{Bl}_J(A)\Big)_{(a)}
                \end{split}
            \end{equation}
            where $A\big[\frac{J}{a}\big]$ is the homogeneous localisation of $\mathrm{Bl}_J(A)$ at $a\in J=\mathrm{Bl}_J(A)_1$, called the \emph{affine blowup algebra}. 
            
            Let $a\in J_{\lambda=\max}=\mathrm{Fit}_k(\phi)_{\lambda=0}$ be a non-zero element of the following form 
            \begin{equation}
                a=\det\begin{pmatrix}
                \xi_1.f_1&\cdots&\xi_1.f_{r-k}\\
                \vdots&\ddots&\vdots\\
                \xi_{r-k}.f_1&\cdots&\xi_{r-k}.f_{r-k}
                \end{pmatrix}
            \end{equation}
            for $f_1,\cdots,f_{r-k}\in A_{\lambda=-w}$. It is easy to see that $(\tilde Y)^0_{\lambda-\min}$ is covered by $\tilde Y_a$ for such $a$. We prove for such $a\in J_{\lambda=\max}$ 
            \begin{equation}
                \mathrm{Fit}_{k-1}(\tilde\phi_a)=0,\quad \mathrm{Fit}_k(\tilde\phi_a)=\mathcal O(\tilde Y_a)
            \end{equation}
            where $\tilde\phi_a$ is the morphism 
            \begin{equation}
                \begin{tikzcd}
                    \Omega_{\mathcal O(\tilde Y_a)/\Bbbk}\ar[r,"\tilde\phi_a"]&\mathfrak u^*\otimes\mathcal O(\tilde Y_a)\ar[r]&\mathrm{coker}(\tilde \phi_a)\ar[r]&0
                \end{tikzcd}
            \end{equation}
            and $\mathrm{Fit}_k(\tilde\phi_a)$ denotes the $k$th Fitting ideal of $\mathrm{coker}(\tilde\phi_a)$. 
            
            {
                \begin{lemma*}
                    $\mathrm{Fit}_{k-1}(\tilde\phi_a)=0$. 
                \end{lemma*}
                \begin{proof}
                    There is a surjective morphism 
                    \begin{equation}
                        A\Big[\frac{J}{a}\Big]\otimes_{\Bbbk}\mathrm{d}\Big(\frac{J}{a}\Big)\to \Omega_{\mathcal O(\tilde Y_a)/\Bbbk},\quad \frac{g}{a^m}\otimes \mathrm{d}\Big(\frac{h}{a}\Big)\mapsto \frac{g}{a^m} \mathrm{d}\Big(\frac{h}{a}\Big)
                    \end{equation}
                    where $m\in\mathbb N,g\in J^m,h\in J$. 
                    
                    Let $h_1,\cdots,h_{r-k+1}\in J$. Then $\tilde\phi_a$ on $\mathrm{d}\big(\frac{h_1}{a}\big),\cdots,\mathrm{d}\big(\frac{h_{r-k+1}}{a}\big)$ is represented by the following matrix 
                    \begin{equation}
                        \begin{split}
                            \tilde\phi_a\Big(\mathrm{d}\Big(\frac{h_1}{a}\Big),\cdots,\mathrm{d}\Big(\frac{h_{r-k+1}}{a}\Big)\Big)=\big(u_1\otimes1,\cdots,u_r\otimes 1\big)\begin{pmatrix}\frac{\xi_1.h_1}{a}&\cdots&\frac{\xi_1.h_{r-k+1}}{a}\\
                            \vdots&\vdots&\vdots\\
                            \frac{\xi_r.h_1}{a}&\cdots&\frac{\xi_r.h_{r-k+1}}{a}\end{pmatrix}. 
                        \end{split}
                    \end{equation}
                    
                    To prove $\mathrm{Fit}_k(\tilde\phi_a)=0$, it suffices to show all $(r-k+1)\times(r-k+1)$-minors of the above matrix are zero. Without loss of generality we consider the top $(r-k+1)\times(r-k+1)$-minor, which is of the form $\frac{\det(M)}{a^{r-k+1}}$, where $\det(M)$ is the determinant of the following square matrix
                    \begin{equation}
                        M=\begin{pmatrix}\xi_1.h_1&\cdots&\xi_1.h_{r-k+1}\\
                        \vdots&\ddots&\vdots\\
                        \xi_{r-k+1}.h_1&\cdots&\xi_{r-k+1}.h_{r-k+1}\end{pmatrix}. 
                    \end{equation}
                    Since $\mathrm{Fit}_{k-1}(\phi)=0$, we have $\det(M)=0$. This proves $\mathrm{Fit}_{k-1}(\tilde\phi_a)=0$. 
                \end{proof}
                
                \begin{lemma*}
                    $\mathrm{Fit}_k(\tilde\phi_a)=\mathcal O(\tilde Y_a)$. 
                \end{lemma*}
                \begin{proof}
                    This is to prove $1\in \mathrm{Fit}_k(\tilde\phi_a)$. For $1\leq i\leq r-k$, let 
                    \begin{equation}
                        b_i:=\det\begin{pmatrix}\xi_1.f_1&\cdots&\xi_1.f_{r-k}\\
                        \vdots&\cdots&\vdots\\
                        f_1&\cdots&f_{r-k}\\
                        \vdots&\cdots&\vdots\\
                        \xi_{r-k}.f_1&\cdots&\xi_{r-k}.f_{r-k}\end{pmatrix}\in A
                    \end{equation}
                    where $(f_1,\cdots,f_{r-k})$ occupies the $i$th row. We see that $b_i\in A_{\lambda=-w}$. 
                    
                    It is easy to see $\xi_i.b_j=\delta_{ij}a$ for $1\leq i,j\leq r-k$. Moreover, we have 
                    $\mathrm{U}(\mathfrak u).b_j\subseteq \mathrm{Fit}_k(\phi)\subseteq I$. By Corollary \ref{corollary of the ideal cutting out U-sweep of Z consists functions whose UEA derivations in the ideal of Z}, the map $A\to\mathcal O(U)\otimes (A/I)$ sends $b_j$ to zero, i.e. $b_j\in J$. Then $\mathrm{d}\big(\frac{b_j}{a}\big)\in\Omega_{\mathcal O(\tilde Y_a)/\Bbbk}$. 
                    
                    The map $\tilde\phi_a:\Omega_{\mathcal O(\tilde Y_a)/\Bbbk}\to\mathfrak u^*\otimes\mathcal O(\tilde Y_a)$ on $\mathrm{d}\big(\frac{b_1}{a}\big),\cdots,\mathrm{d}\big(\frac{b_{r-k}}{a}\big)$ is represented by the following matrix 
                    \begin{equation}
                        \tilde\phi_a\Big(\mathrm{d}\Big(\frac{b_1}{a}\Big),\cdots,\mathrm{d}\Big(\frac{b_{r-k}}{a}\Big)\Big)=\big(u_1\otimes 1,\cdots,u_r\otimes 1\big)\begin{pmatrix}I_{r-k}\\M'\end{pmatrix}
                    \end{equation}
                    where $I_{r-k}$ is the identity matrix of size $(r-k)\times(r-k)$, and $M'\in \mathcal O(\tilde Y_a)^{k\times (r-k)}$. Then the top $(r-k)\times(r-k)$-minor is $1\in\mathrm{Fit}_k(\tilde\phi_a)$. 
                \end{proof}
            }
            
            In the case when there is one grading weight, when proving $\mathrm{Fit}_k(\tilde \phi_a)=\mathcal O(\tilde Y_a)$, we construct $b_1,\cdots,b_{r-k}\in A$ with the following properties: 
            \begin{itemize}
                \item[(1)] $b_1,\cdots,b_{r-k}\in J$, so $\mathrm{d}\big(\frac{b_j}{a}\big)\in\Omega_{\mathcal O(\tilde Y_a)/\Bbbk}$, the domain of definition of $\tilde\phi_a$;  
                \item[(2)] $\xi_i.b_j=\delta_{ij}a$ for $1\leq i,j\leq r-k$, which realises $I_{r-k}$ as a sub-matrix of size $(r-k)\times(r-k)$ of the matrix representing $\tilde\phi_a$. 
            \end{itemize}
            
            In the general case with multiple grading weights, property $(1)$ above, $b_j\in J$, becomes much complicated, while property $(2)$ and $\mathrm{Fit}_{k-1}(\tilde\phi_i)=0$ can be verified essentially the same way. Part of the difficulty of $(1)$ is from the difficulty of describing scheme theoretic images. The blow-up centre is such a scheme theoretic image, and the condition $b_j\in J$ is that $b_j$ pulls back to a zero function on the centre. Another non-trivial part for $(1)$ is to guess what is the form of $b_j$ in general. 
        }
        
        \subsubsection{Two grading weights}
        {
            We focus on the construction of $b_j$ and the condition $b_j\in J$ for cases with 2 grading weights. Assume the condition that $\prod_{i=1}^2\mathrm{Fit}_{k_i}(\phi_i)\not\subseteq A_{\lambda<0}$. 
            
            Let $r_i:=\dim\mathfrak u_{\lambda=w_i}$. Choose a basis $\big(\xi^{(i)}_1,\cdots,\xi^{(i)}_{r_i}\big)$ of $\mathfrak u_{\lambda=w_i}$ and let $\big(u^{(i)}_1,\cdots,u^{(i)}_{r_i}\big)$ be its dual basis on $\mathfrak u_{\lambda=w_i}^*$. Let $I:=\prod_{i=1}^2\mathrm{Fit}_{k_i}(\phi_i)+A_{\lambda<0}$ and let $J\subseteq A$ be the kernel of $A\to \mathcal O(U)\otimes(A/I)$. 
            
            It is easy to see $\prod_{i=1}^2\mathrm{Fit}_{k_i}(\phi_i)$ is generated by $a\in\prod_{i=1}^2\mathrm{Fit}_{k_i}(\phi_i)$ of the form $a=a^{(1)}a^{(2)}$ for $a^{(i)}\in\mathrm{Fit}_{k_i}(\phi_i)$
            \begin{equation}
                \begin{split}
                   a^{(i)}:=&\det\begin{pmatrix}\xi^{(i)}_1.f^{(i)}_1&\cdots&\xi^{(i)}_1.f^{(i)}_{r_i-k_i}\\
                    \vdots&\ddots&\vdots\\
                    \xi^{(i)}_{r_i-k_i}.f^{(i)}_1&\cdots&\xi^{(i)}_{r_i-k_i}.f^{(i)}_{r_i-k_i}\end{pmatrix}\\
                    &f^{(i)}_1,\cdots,f^{(i)}_{r_i-k_i}\in A_{\lambda=-w_i}. 
                \end{split}
            \end{equation}
            
            Recall 
            \begin{equation}
                \begin{split}
                    \tilde \phi_{1,a}:&\Omega_{\mathcal O(\tilde Y_a)/\Bbbk}\to \mathfrak u_1^*\otimes\mathcal O(\tilde Y_a)\\
                    \tilde \phi_{2,a}:&\tilde K_{1,a}\to (\mathfrak u_2/\mathfrak u_1)^*\otimes \mathcal O(\tilde Y_a)\cong(\mathfrak u_{\lambda=w_2})^*\otimes \mathcal O(\tilde Y_a)
                \end{split}
            \end{equation}
            where $\tilde K_{1,a}:=\ker(\tilde\phi_{1,a})\subseteq\Omega_{\mathcal O(\tilde Y_a)/\Bbbk}$. Let $b^{(2)}_\nu$ for $\nu=1,\cdots,r_2-k_2$ be 
            \begin{equation}
                b^{(2)}_\nu:=a^{(1)}\cdot\det\begin{pmatrix}\xi^{(2)}_1.f^{(2)}_1&\cdots&\xi^{(2)}_1.f^{(2)}_{r_2-k_2}\\\vdots&\cdots&\vdots\\f^{(2)}_1&\cdots&f^{(2)}_{r_2-k_2}\\\vdots&\cdots&\vdots\\\xi^{(2)}_{r_2-k_2}.f^{(2)}_1&\cdots&\xi^{(2)}_{r_2-k_2}.f^{(2)}_{r_2-k_2}\end{pmatrix}
            \end{equation}
            where $\big(f^{(2)}_1,\cdots,f^{(2)}_{r_2-k_2}\big)$ occupies the $\nu$th row. It is easy to check: 
            \begin{itemize}
                \item[(1)] $b^{(2)}_\nu\in A_{\lambda=-w_2}$; 
                \item[(2)] $\xi^{(2)}_\mu.b^{(2)}_\nu=\delta_{\mu\nu}a$ for $1\leq \mu,\nu\leq r_2-k_2$, and $\xi^{(2)}_\mu.b^{(2)}_\nu\in a^{(1)}\cdot\mathrm{Fit}_{k_2}(\phi_2)$ for $1\leq \mu\leq r_2$; 
                \item[(3)] $\xi^{(1)}_\mu.b^{(2)}_\nu=0$ for $1\leq \mu\leq r_1$ and $1\leq \nu\leq r_2-k_2$. 
            \end{itemize}
            The last two conditions above show that $b^{(2)}_\nu\in J$ (Corollary \ref{corollary of the ideal cutting out U-sweep of Z consists functions whose UEA derivations in the ideal of Z}). Then $\mathrm{d}\big(\frac{b^{(2)}_\nu}{a}\big)\in\Omega_{\mathcal O(\tilde Y_a)/\Bbbk}$, and since it has $\lambda$-weight $-w_2$, we have $\mathrm{d}\big(\frac{b^{(2)}_\nu}{a}\big)\in\tilde K_{1,a}$. We can write down the matrix representing $\tilde\phi_{2,a}$ on $\mathrm{d}\big(\frac{b^{(2)}_1}{a}\big),\cdots,\mathrm{d}\big(\frac{b^{(2)}_{r_2-k_2}}{a}\big)$ and the top sub-matrix of size $(r_2-k_2)\times(r_2-k_2)$ is $I_{r_2-k_2}$ by condition $(2)$ above. This proves $\mathrm{Fit}_{k_2}(\tilde\phi_{2,a})=\mathcal O(\tilde Y_a)$.

            {
                We need to show $\mathrm{Fit}_{k_1}(\tilde\phi_{1,a})=\mathcal O(\tilde Y_a)$. Before we prove it, we prove an identity. 
                \begin{lemma*}
                    \begin{equation}\label{equation of determinantal sum when 2 weights}
                        \big(\xi^{(2)}_\mu.f\big)a^{(2)}=\sum_{\mu'=1}^{r_2-k_2}\big(\xi^{(2)}_\mu.b^{(2)}_{\mu'}\big)\big(\xi^{(2)}_{\mu'}.f\big),\quad f\in A_{\lambda=-w_2},\;1\leq \mu\leq r_2. 
                    \end{equation}
                \end{lemma*}
                \begin{proof}
                    For $f\in A_{\lambda=-w_2}$ and $1\leq\mu\leq r_2$, consider the following matrix of size $(r_2-k_2+1)\times (r_2-k_2+1)$
                    \begin{equation}
                        \begin{pmatrix}
                            \xi^{(2)}_1.f^{(2)}_1&\cdots&\xi^{(2)}_1.f^{(2)}_{r_2-k_2}&\xi^{(2)}_1.f\\
                            \vdots&\ddots&\vdots&\vdots\\
                            \xi^{(2)}_{r_2-k_2}.f^{(2)}_1&\cdots&\xi^{(2)}_{r_2-k_2}.f^{(2)}_{r_2-k_2}&\xi^{(2)}_{r_2-k_2}.f\\
                            \xi^{(2)}_\mu.f^{(2)}_1&\cdots&\xi^{(2)}_\mu.f^{(2)}_{r_2-k_2}&\xi^{(2)}_\mu.f
                        \end{pmatrix}. 
                    \end{equation}
                    Since $f\in A_{\lambda=-w_2}$, we have $\mathrm{d}f\in K_1=\ker(\phi_1)$. Then the above matrix has determinant zero, since $\mathrm{Fit}_{k_2-1}(\phi_2)=0$. Its determinant can be expanded through the last column 
                    \begin{equation}
                        0=\big(\xi^{(2)}_\mu.f\big)a^{(2)}-\sum_{\mu'=1}^{r_2-k_2}\big(\xi^{(2)}_{\mu'}.f\big)\det\begin{matrix}\begin{pmatrix}\xi^{(2)}_1.f^{(2)}_1&\cdots&\xi^{(2)}_1.f^{(2)}_{r_2-k_2}\\\vdots&\cdots&\vdots\\\xi^{(2)}_\mu.f^{(2)}_1&\cdots&\xi^{(2)}_\mu.f^{(2)}_{r_2-k_2}\\\vdots&\cdots&\vdots\\\xi^{(2)}_{r_2-k_2}.f^{(2)}_1&\cdots&\xi^{(2)}_{r_2-k_2}.f^{(2)}_{r_2-k_2}\end{pmatrix}&\leftarrow\mu'\textrm{th row}\end{matrix}
                    \end{equation}
                    where the summand labelled by $\mu'$ is the determinant $a^{(2)}$ with the $\mu'$th row replaced. Equation \eqref{equation of determinantal sum when 2 weights} is the compact form of the above equation. 
                \end{proof}
            }

            For $1\leq\nu\leq r_1-k_1$, let $\beta^{(1)}_\nu$ be 
            \begin{equation}
                \beta^{(1)}_\nu:=\det\begin{pmatrix}\xi^{(1)}_1.f^{(1)}_1&\cdots&\xi^{(1)}_1.f^{(1)}_{r_1-k_1}\\\vdots&\cdots&\vdots\\f^{(1)}_1&\cdots&f^{(1)}_{r_1-k_1}\\\vdots&\cdots&\vdots\\\xi^{(1)}_{r_1-k_1}.f^{(1)}_1&\cdots&\xi^{(1)}_{r_1-k_1}.f^{(1)}_{r_1-k_1}\end{pmatrix}
            \end{equation}
            where $\big(f^{(1)}_1,\cdots,f^{(1)}_{r_1-k_1}\big)$ occupies the $\nu$th row. For $1\leq \nu\leq r_1-k_1$, let $b^{(1)}_\nu$ be 
            \begin{equation}
                b^{(1)}_\nu:=\beta^{(1)}_\nu a^{(2)}-\frac{w_2}{w_1}\sum_{\rho=1}^{r_2-k_2}\big(\xi^{(2)}_\rho.\beta^{(1)}_\nu\big)b^{(2)}_\rho. 
            \end{equation}
            We have: 
            \begin{itemize}
                \item $b^{(1)}_\nu\in A_{\lambda=-w_1}$; 
                \item $\xi^{(1)}_\mu.b^{(1)}_\nu=\xi^{(1)}_\mu.\big(\beta^{(1)}_\nu a^{(2)}\big)\in\prod_{i=1}^2\mathrm{Fit}_{k_i}(\phi_i)$ for $1\leq \mu\leq r_1$ and $1\leq \nu\leq r_1-k_1$; 
                \item In particular, $\xi^{(1)}_\mu.b^{(1)}_\nu=\delta_{\mu\nu} a$ for $1\leq \mu,\nu\leq r_1-k_1$. 
            \end{itemize}
            
            If we can prove $b^{(1)}_\nu\in J$, then $\mathrm{d}\big(\frac{b^{(1)}_\nu}{a}\big)\in \Omega_{\mathcal O(\tilde Y_a)/\Bbbk}$, and the matrix representing $\tilde\phi_{1,a}$ on $\mathrm{d}\big(\frac{b^{(1)}_1}{a}\big),\cdots,\mathrm{d}\big(\frac{b^{(1)}_{r_1-k_1}}{a}\big)$ is 
            \begin{equation}
                \tilde\phi_{1,a}\Big(\mathrm{d}\Big(\frac{b^{(1)}_1}{a}\Big),\cdots,\mathrm{d}\Big(\frac{b^{(1)}_{r_1-k_1}}{a}\Big)\Big)=\big(u^{(1)}_1\otimes1,\cdots,u^{(1)}_{r_1}\otimes1\big)\begin{pmatrix}I_{r_1-k_1}\\M'\end{pmatrix}
            \end{equation}
            where $M'\in\mathcal O(\tilde Y_a)^{k_1\times(r_1-k_1)}$. Then $1\in\mathrm{Fit}_{k_1}(\tilde\phi_{1,a})$, i.e. $\mathrm{Fit}_{k_1}(\tilde\phi_{1,a})=\mathcal O(\tilde Y_a)$. 
            
            The final step is to verify $b^{(1)}_\nu\in J$. By Corollary \ref{corollary of the ideal cutting out U-sweep of Z consists functions whose UEA derivations in the ideal of Z}, it suffices to check: 
            \begin{itemize}
                \item $\xi^{(1)}_\mu.b^{(1)}_\nu\in \prod_{i=1}^2\mathrm{Fit}_{k_i}(\phi_i)$ for $1\leq \mu\leq r_1$, which is seen above; 
                \item if $w_1/w_2\in\mathbb N_+$, then $\xi^{(2)}_{\mu_1}\cdots\xi^{(2)}_{\mu_m}.b^{(1)}_\nu\in \prod_{i=1}^2\mathrm{Fit}_{k_i}(\phi_i)$ for $m=w_1/w_2\in\mathbb N_+$ and all $1\leq\mu_1,\cdots,\mu_m\leq r_1$. 
            \end{itemize}
            
            We have $m\geq 2$ since $w_1>w_2$. When $m=2$, we have 
            \begin{equation}
                \begin{split}
                    \xi^{(2)}_{\mu_1}\xi^{(2)}_{\mu_2}.b^{(1)}_\nu=&\big(\xi^{(2)}_{\mu_1}\xi^{(2)}_{\mu_2}.\beta^{(1)}_\nu\big)a^{(2)}\\
                    &-\frac{1}{2}\sum_{\rho=1}^{r_2-k_2}\Big(\big(\xi^{(2)}_{\mu_1}\xi^{(2)}_\rho.\beta^{(1)}_\nu\big)\big(\xi^{(2)}_{\mu_2}.b^{(2)}_\rho\big)+\big(\xi^{(2)}_{\mu_2}\xi^{(2)}_\rho.\beta^{(1)}_\nu\big)\big(\xi^{(2)}_{\mu_1}.b^{(2)}_\rho\big)\Big)
                \end{split}
            \end{equation}
            where we have used that $\xi^{(2)}_{\mu_1}\xi^{(2)}_{\mu_2}.\big(\xi^{(2)}_\rho.\beta^{(1)}_\nu\big)=0$ and $\xi^{(2)}_{\mu_1}\xi^{(2)}_{\mu_2}.b^{(2)}_\rho=0$. 
            
            Apply equation \eqref{equation of determinantal sum when 2 weights} to $\big(\xi^{(2)}_{\mu_1}\xi^{(2)}_{\mu_2}.\beta^{(1)}_\nu\big)a^{(2)}$ 
            \begin{equation}
                \big(\xi^{(2)}_{\mu_1}\xi^{(2)}_{\mu_2}.\beta^{(1)}_\nu\big)a^{(2)}=\sum_{\rho=1}^{r_2-k_2}\big(\xi^{(2)}_{\mu_1}.b^{(2)}_{\rho}\big)\big(\xi^{(2)}_{\rho}\xi^{(2)}_{\mu_2}.\beta^{(1)}_\nu\big). 
            \end{equation}
            Then 
            \begin{equation}
                \begin{split}
                    \xi^{(2)}_{\mu_1}\xi^{(2)}_{\mu_2}.b^{(1)}_\nu=&\frac{1}{2}\sum_{\rho=1}^{r_2-k_2}\big(\big[\xi^{(2)}_\rho,\xi^{(2)}_{\mu_2}\big].\beta^{(1)}_\nu\big)\big(\xi^{(2)}_{\mu_1}.b^{(2)}_\rho\big)\\
                    &+\frac{1}{2}\sum_{\rho=1}^{r_2-k_2}\big(\xi^{(2)}_\rho\xi^{(2)}_{\mu_2}.\beta^{(1)}_\nu\big)\big(\xi^{(2)}_{\mu_1}.b^{(2)}_\rho\big)-\frac{1}{2}\sum_{\rho=1}^{r_2-k_2}\big(\xi^{(2)}_{\mu_1}\xi^{(2)}_\rho.\beta^{(1)}_\nu\big)\big(\xi^{(2)}_{\mu_2}.b^{(2)}_\rho\big). 
                \end{split}
            \end{equation}
            We can exchange $\mu_1$ and $\mu_2$ to get a similar expression for $\xi^{(2)}_{\mu_2}\xi^{(2)}_{\mu_1}.b^{(1)}_\nu$. Then 
            \begin{equation}
                \begin{split}
                    2\xi^{(2)}_{\mu_1}\xi^{(2)}_{\mu_2}.b^{(1)}_\nu=&\big[\xi^{(2)}_{\mu_1},\xi^{(2)}_{\mu_2}\big].b^{(1)}_\nu+\big(\xi^{(1)}_{\mu_1}\xi^{(1)}_{\mu_2}+\xi^{(1)}_{\mu_2}\xi^{(1)}_{\mu_1}\big).b^{(1)}_\nu\\
                    =&\big[\xi^{(2)}_{\mu_1},\xi^{(2)}_{\mu_2}\big].b^{(1)}_\nu\\
                    &+\sum_{\rho=1}^{r_2-k_2}\Big(\big(\big[\xi^{(2)}_\rho,\xi^{(2)}_{\mu_2}\big].\beta^{(1)}_\nu\big)\big(\xi^{(2)}_{\mu_1}.b^{(2)}_\rho\big)+\big(\big[\xi^{(2)}_\rho,\xi^{(2)}_{\mu_1}\big].\beta^{(1)}_\nu\big)\big(\xi^{(2)}_{\mu_2}.b^{(2)}_\rho\big)\Big). 
                \end{split}
            \end{equation}
            Since $w_1=2w_2$, we have $\big[\xi^{(2)}_\rho,\xi^{(2)}_{\mu_1}\big]\in\mathfrak u_{\lambda=w_1}$. We see that $\xi^{(2)}_{\mu_1}\xi^{(2)}_{\mu_2}.b^{(1)}_\nu\in\prod_{i=1}^2\mathrm{Fit}_{k_i}(\phi_i)$. This proves $b^{(1)}_\nu\in J$ when $w_1=2w_2$. 
            
            When $w_1/w_2=m\geq 3$, we have that $\mathfrak u$ is Abelian, and apply equation \eqref{equation of determinantal sum when 2 weights} to get 
            \begin{equation}
                \begin{split}
                    \xi^{(2)}_{\mu_1}\cdots\xi^{(2)}_{\mu_m}.b^{(1)}_\nu=&\big(\xi^{(2)}_{\mu_1}\cdots\xi^{(2)}_{\mu_m}.\beta^{(1)}_\nu\big)a^{(2)}\\
                    &-\frac{1}{m}\sum_{\rho=1}^{r_2-k_2}\sum_{j=1}^m\big(\xi^{(2)}_{\mu_1}\cdots\widehat{\xi^{(2)}_{\mu_j}}\cdots\xi^{(2)}_{\mu_m}\xi^{(2)}_\rho.\beta^{(1)}_\nu\big)\big(\xi^{(2)}_{\mu_j}.b^{(2)}_\rho\big)\\
                    =&\frac{1}{m}\sum_{j=1}^m\big(\xi^{(2)}_{\mu_j}\xi^{(2)}_{\mu_1}\cdots\widehat{\xi^{(2)}_{\mu_j}}\cdots\xi^{(2)}_{\mu_m}.\beta^{(1)}_\nu\big)a^{(2)}\\
                    &-\frac{1}{m}\sum_{\rho=1}^{r_2-k_2}\sum_{j=1}^m\big(\xi^{(2)}_{\mu_1}\cdots\widehat{\xi^{(2)}_{\mu_j}}\cdots\xi^{(2)}_{\mu_m}\xi^{(2)}_\rho.\beta^{(1)}_\nu\big)\big(\xi^{(2)}_{\mu_j}.b^{(2)}_\rho\big)\\
                    =&\frac{1}{m}\sum_{j=1}^m\sum_{\rho=1}^{r_2-k_2}\big(\xi^{(2)}_\rho\xi^{(2)}_{\mu_1}\cdots\widehat{\xi^{(2)}_{\mu_j}}\cdots\xi^{(2)}_{\mu_m}.\beta^{(1)}_\nu\big)\big(\xi^{(2)}_{\mu_j}.b^{(2)}_\rho\big)\\
                    &-\frac{1}{m}\sum_{\rho=1}^{r_2-k_2}\sum_{j=1}^m\big(\xi^{(2)}_{\mu_1}\cdots\widehat{\xi^{(2)}_{\mu_j}}\cdots\xi^{(2)}_{\mu_m}\xi^{(2)}_\rho.\beta^{(1)}_\nu\big)\big(\xi^{(2)}_{\mu_j}.b^{(2)}_\rho\big)\\
                    =&0\in\prod_{i=1}^2\mathrm{Fit}_{k_i}(\phi_i). 
                \end{split}
            \end{equation}
            This proves $b^{(1)}_\nu\in J$. 
            
            In the case with 2 grading weights, the construction is no longer simple. In Lemma \ref{lemma of existence of b^(i)_nu}, the general construction of $b^{(i)}_\nu$ has to be formulated inductively, and equation \eqref{equation of determinantal sum when 2 weights} in the general situation is Lemma \ref{lemma of determinantal sum}. The assumption $\prod_{i=0}^2\mathrm{Fit}_{k_i}(\phi_i)\not\subseteq A_{\lambda<0}$ generalises to Condition \ref{condition before blow-up for CDRS}. 
        }
    }

    \subsection{The Weak Unipotent Upstairs condition}
    {
        In \cite{hoskins2021quotients}, the so-called \emph{Weak Unipotent Upstairs Assumption} (\cite{hoskins2021quotients} Assumption 4.42) was assumed for the sequence of blow-ups described in that paper. In this subsection, we describe a similar assumption for our blow-up procedure and call it the \emph{Weak Unipotent Upstairs Condition} as well. 
        
        Recall $\hat U=U\rtimes_\lambda\mathbb G_m$ for a grading 1PS $\lambda:\mathbb G_m\to \mathrm{Aut}(U)$ with grading weights $w_1>\cdots>w_n$. For a linear action $\hat U\curvearrowright(X,L)$, the open subscheme $X^0_{\lambda-\min}\subseteq X$ is defined as the non-vanishing locus of $\bigoplus_{d>0}H^0(X,L^d)_{\lambda=\max}$. Recall the notation 
        \begin{equation}
            \begin{split}
                &\mathcal K_i:=\ker\big(\Omega_{X^0_{\lambda-\min}/\Bbbk}\to \mathfrak u_i^*\otimes_\Bbbk\mathcal O_{X^0_{\lambda-\min}}\big)\\
                &\phi_i:\mathcal K_{i-1}\to (\mathfrak u_i/\mathfrak u_{i-1})^*\otimes_\Bbbk\mathcal O_{X^0_{\lambda-\min}}\\
                &\mathrm{Fit}_k(\phi_i):=\mathrm{Fit}_k(\mathrm{coker}(\phi_i)). 
            \end{split}
        \end{equation}

        {
            We call the following condition \emph{Weak Unipotent Upstairs}. The name comes from \cite{hoskins2021quotients} Assumption 4.42. 
            \begin{condition}\label{condition before blow-up for CDRS}
                Let $\hat U \curvearrowright(X,L)$, $\{\mathfrak u_i\}_{i=0}^n$ and $\phi_i$ be as above. Let $k_i:=\min\big\{k\in\mathbb N:\mathrm{Fit}_k(\phi_i)\ne0\big\}$. This condition for the action $\hat U \curvearrowright (X,L)$ assumes that there exists a closed point $z\in Z_{\lambda-\min}$ such that $\dim\mathrm{Stab}_{\mathfrak u_i}(z)=k_1+\cdots+k_i$ for all $1\leq i\leq n$.  
            \end{condition}
            
            \begin{remark}
                Condition \ref{condition of CDRS for projective} implies  Condition \ref{condition before blow-up for CDRS} by Corollary \ref{corollary of flatness of coker(phi_i), Q(g_i) and dim of stab when reduced}.
            \end{remark}
        }
        
        {
            \begin{lemma}\label{lemma of point-wise description of the nice locus of small unipotent stabilisers}
                Let $\hat U\curvearrowright (X,L)$ with the grading 1PS $\lambda:\mathbb G_m\to \mathrm{Aut}(U)$ with weights $w_1>\cdots>w_n$. Let $\{\mathfrak u_i\}_{i=0}^n$ be the filtration of $\mathfrak u$ such that $\mathfrak u_i:=\mathfrak u_{\lambda\geq w_i}$. Let $k_i:=\min\big\{k\in\mathbb N:\mathrm{Fit}_k(\phi_i)\ne0\big\}$ for $1\leq i\leq n$. Then for each $1\leq i\leq n$, the following subsets of closed points of $X^0_{\lambda-\min}$ are the same: 
                \begin{equation}
                    \begin{split}
                        S^{(1)}_i:=&\big\{x\in X^0_{\lambda-\min}:x\textrm{ is a closed point in the non-vanishing locus of }\prod_{i'=1}^i\mathrm{Fit}_{k_{i'}}(\phi_{i'})\big\}\\
                        S^{(2)}_i:=&\Big\{x\in X^0_{\lambda-\min}:\begin{matrix}x\textrm{ is a closed point and for }i'=1,\cdots,i\\ \dim\mathrm{Stab}_{\mathfrak u_{i'}}(x)\leq k_1+\cdots+k_{i'}\end{matrix}\Big\}\\
                        S^{(3)}_i:=&\Big\{x\in X^0_{\lambda-\min}:\begin{matrix}x\textrm{ is a closed point and for }i'=1,\cdots,i\\ \dim\mathrm{Stab}_{\mathfrak u_{i'}}(x)= k_1+\cdots+k_{i'}\end{matrix}\Big\}. 
                    \end{split}
                \end{equation}
                Moreover the subsets are all non-empty if $X^0_{\lambda-\min}$ is integral, i.e. reduced and irreducible. 
            \end{lemma}
            \begin{proof}
                If $X^0_{\lambda-\min}$ is reduced, then the non-vanishing locus of $\mathrm{Fit}_{k_i}(\phi_i)$ is an non-empty open subset for each $1\leq i\leq n$. Finite intersections of non-empty open subsets are non-empty if $X^0_{\lambda-\min}$ is irreducible. This proves $S^{(1)}_i\ne\emptyset$ for all $1\leq i\leq n$ if $X^0_{\lambda-\min}$ is integral. 

                We will then prove $S^{(1)}_i=S^{(2)}_i=S^{(3)}_i$ by induction on $i=1,\cdots,n$. 
                
                Since $X$ is of finite type over $\Bbbk$, open subsets are determined uniquely by their closed points (\cite[\href{https://stacks.math.columbia.edu/tag/02J6}{Lemma 02J6}]{stacks-project}, \cite[\href{https://stacks.math.columbia.edu/tag/005Z}{Lemma 005Z}]{stacks-project}). We will also denote the non-vanishing locus of $\prod_{i'=1}^i\mathrm{Fit}_{k_{i'}}(\phi_{i'})$ by $S^{(1)}_i$, which is open in $X^0_{\lambda-\min}$. 
                
                Initially consider when $i=1$. Let $x:\mathrm{Spec}(\Bbbk)\to X^0_{\lambda-\min}$ be a closed point in $X^0_{\lambda-\min}$ and let $x^\#:\mathcal O_{X^0_{\lambda-\min}}\to x_*(\Bbbk)$ be the associated morphism of sheaves of rings. Then $x\not\in S^{(1)}_1$ if and only if $\mathrm{Fit}_{k_1}(\phi_1)\subseteq\ker x^\#$, which is equivalent to the condition that $\dim\mathrm{Stab}_{\mathfrak u_1}(x)>k_1$ by Lemma \ref{lemma of dim of relative stab and Fit when flat}, i.e. $x\notin S^{(2)}_1$. This proves $S^{(1)}_1=S^{(2)}_1$. Obviously $S^{(3)}_1\subseteq S^{(2)}_1$. Assume the closed point $x$ is in $S^{(2)}_1\setminus S^{(3)}_1$. Let $l_1:=\dim\mathrm{Stab}_{\mathfrak u_1}(x)<k_1$. Then $\mathrm{Fit}_{l_1}(\phi_1)=0$ since $k_1$ is the smallest integer such that $\mathrm{Fit}_k(\phi_1)\ne0$. We have $\mathrm{Fit}_{l_1}(\phi_1)\subseteq\ker x^\#$ and then $l_1=\dim\mathrm{Stab}_{\mathfrak u_1}(x)>l_1$ by Lemma \ref{lemma of dim of relative stab and Fit when flat}, which is a contradiction. This proves $S^{(2)}_1=S^{(3)}_1$, so $S^{(1)}_1=S^{(2)}_1=S^{(3)}_1$. 
                
                Assume $S^{(1)}_i=S^{(2)}_i=S^{(3)}_i$ for some $1\leq i< n$. If $S^{(1)}_i=S^{(2)}_i=S^{(3)}_i=\emptyset$, then $S^{(1)}_{i+1}=S^{(2)}_{i+1}=S^{(3)}_{i+1}=\emptyset$. Without loss of generality, we assume $S^{(1)}_i\ne\emptyset$. We have that $\mathrm{coker}(\phi_{i'})|_{S^{(1)}_i}$ is locally free of rank $k_{i'}$ for each $1\leq i'\leq i$ by Proposition \ref{proposition of properties of Fit}. By Corollary \ref{corollary of flatness of coker(phi_i), Q(g_i) and dim of stab when reduced}, we have that $Q(\mathfrak u_{i'})|_{S^{(1)}_i}$ is locally free of rank $k_1+\cdots+k_{i'}$ for each $1\leq i'\leq i$, where $Q(\mathfrak u_{i'})$ is the cokernel of $\Omega_{X^0_{\lambda-\min}/\Bbbk}\to\mathfrak u_i^*\otimes_{\Bbbk}\mathcal O_{X^0_{\lambda-\min}}$. In particular $Q(\mathfrak u_i)|_{S^{(1)}_i}$ is flat. Let $x:\mathrm{Spec}(\Bbbk)\to X^0_{\lambda-\min}$ be a closed point. If $x\in S^{(1)}_i$, then the condition that $\dim\mathrm{Stab}_{(\mathfrak u_{i+1}:\mathfrak u_i)}(x)\leq k_{i+1}$ is equivalent to the condition that $\mathrm{Fit}_{k_{i+1}}(\phi_{i+1})\not\subseteq \ker x^\#$ by Lemma \ref{lemma of dim of relative stab and Fit when flat}. Then we have 
                \begin{equation}
                    \begin{split}
                        x\in S^{(2)}_{i+1}\quad\textrm{if and only if}\quad &\begin{cases}x\in S^{(1)}_i=S^{(2)}_i=S^{(3)}_i\\\dim\mathrm{Stab}_{\mathfrak u_{i+1}}(x)\leq k_1+\cdots+k_i+ k_{i+1}\end{cases}\\
                        \quad\textrm{if and only if}\quad &\begin{cases}x\in S^{(1)}_i=S^{(2)}_i=S^{(3)}_i\\\dim\mathrm{Stab}_{(\mathfrak u_{i+1}:\mathfrak u_i)}(x)\leq k_{i+1}\end{cases}\\
                        \quad\textrm{if and only if}\quad &\begin{cases}x\in S^{(1)}_i=S^{(2)}_i=S^{(3)}_i\\\mathrm{Fit}_{k_{i+1}}(\phi_{i+1})\not\subseteq\ker x^\#\end{cases}\\
                        \quad\textrm{if and only if}\quad & x\in S^{(1)}_{i+1}. 
                    \end{split}
                \end{equation}
                This proves $S^{(1)}_{i+1}=S^{(2)}_{i+1}$. 
                
                Obviously we have $S^{(3)}_{i+1}\subseteq S^{(2)}_{i+1}$. Let $x:\mathrm{Spec}(\Bbbk)\to X^0_{\lambda-\min}$ be a closed point in $S^{(2)}_{i+1}\setminus S^{(3)}_{i+1}$. Then $x\in S^{(3)}_i$ and $\dim\mathrm{Stab}_{(\mathfrak u_{i+1}:\mathfrak u_i)}(x)<k_{i+1}$. As above, we have that $Q(\mathfrak u_i)|_{S^{(1)}_i}$ is flat and $x\in S^{(1)}_i$. Let $l_{i+1}:=\dim\mathrm{Stab}_{(\mathfrak u_{i+1}:\mathfrak u_i)}(x)<k_{i+1}$. We have that $\mathrm{Fit}_{l_{i+1}}(\phi_{i+1})=0$ since $k_{i+1}$ is the smallest integer such that $\mathrm{Fit}_k(\phi_{i+1})\ne0$. In particular $\mathrm{Fit}_{l_{i+1}}(\phi_{i+1})\subseteq\ker x^\#$. Then $\dim\mathrm{Stab}_{(\mathfrak u_{i+1}:\mathfrak u_i)}(x)>l_{i+1}$ by Lemma \ref{lemma of dim of relative stab and Fit when flat}, which is a contradiction. This proves $S^{(2)}_{i+1}=S^{(3)}_{i+1}$, so $S^{(1)}_{i+1}=S^{(2)}_{i+1}=S^{(3)}_{i+1}$. 
            \end{proof}
            
            \begin{remark}
                We can define integers $d_1,\cdots,d_n\in\mathbb N$ by 
                \begin{equation}
                    \begin{split}
                        d_1:=&\min\big\{\dim\mathrm{Stab}_{\mathfrak u_1}(x):x\in X^0_{\lambda-\min}\textrm{ is a closed point}\big\}\\
                        d_i:=&\min\Big\{\dim\mathrm{Stab}_{\mathfrak u_i}(x):\begin{matrix}x\in X^0_{\lambda-\min}\textrm{ is a closed point and } \\\dim\mathrm{Stab}_{\mathfrak u_{i'}}(x)=d_1+\cdots+d_{i'}\textrm{ for }1\leq i'<i\end{matrix}\Big\}. 
                    \end{split}
                \end{equation}
                If the sets of closed points described in Lemma \ref{lemma of point-wise description of the nice locus of small unipotent stabilisers} are non-empty, then $k_i=d_i$ for $1\leq i\leq n$. In particular, we can recover $k_1,\cdots,k_n$ from the underlying topological spaces when $X^0_{\lambda-\min}$ is integral. 
            \end{remark}
            
            \begin{remark}
                Let $k_i,d_i,S^{(3)}_i$ be as above. Condition \ref{condition before blow-up for CDRS} implies the existence of a closed point $z\in Z_{\lambda-\min}$ such that $\dim\mathrm{Stab}_{\mathfrak u_i}(z)=d_1+\cdots+d_i$ for all $1\leq i\leq n$, which is the condition in \cite{berczi2016projective} Remark 7.12. However this implication is not an equivalence, which becomes one if $S^{(1)}_n=S^{(2)}_n=S^{(3)}_n\ne\emptyset$. 
            \end{remark}
        }
    }
    
    \subsection{Centre of the blow-up}
    {
        Assume Condition \ref{condition before blow-up for CDRS}. Let $\mathcal I_{k_1,\cdots,k_n}:=\prod_{i=1}^n\mathrm{Fit}_{k_i}(\phi_i)\subseteq\mathcal O_{X^0_{\lambda-\min}}$. By Lemma \ref{lemma of point-wise description of the nice locus of small unipotent stabilisers}, Condition \ref{condition before blow-up for CDRS} is equivalent to the condition that $Z_{\lambda-\min}\not\subseteq\mathbb V(\mathcal I_{k_1,\cdots,k_n})$. The blow-up centre $\overline{C}\hookrightarrow X$ is a closed subscheme supported on the closure of the $U$-sweep of $Z_{\lambda-\min}\cap \mathbb V(\mathcal I_{k_1,\cdots,k_n})$. In this subsection $\overline C\hookrightarrow X$ will be defined and its intersection $C$ with the open subscheme $X^0_{\lambda-\min}\subseteq X$ will be studied locally. 
        
        Let $m\in\mathbb N_+$ be such that $H^0(X,L^m)$ generates $\bigoplus_{d\geq0}H^0(X,L^{md})$. Let $f\in H^0(X,L^m)_{\lambda=\max}$ be non-nilpotent. Then $X_f\subseteq X^0_{\lambda-\min}$ is a non-empty affine open subscheme. On $X_f$, let $I_f:=\mathcal I_{k_1,\cdots,k_n}(X_f)+\mathcal O(X_f)_{\lambda<0}$. Then $I_f\subseteq\mathcal O(X_f)$ cuts out $Z_{\lambda-\min,f}\cap \mathbb V(\mathcal I_{k_1,\cdots,k_n})$. If $Z_{\lambda-\min,f}\not\subseteq\mathbb V(\mathcal I_{k_1,\cdots,k_n})$, then $\mathcal O(X_f)_{\lambda<0}\subsetneq I_f$, i.e. $(I_f)_{\lambda=\max}=(I_f)_{\lambda=0}$. Else $Z_{\lambda-\min,f}\subseteq\mathbb V(\mathcal I_{k_1,\cdots,k_n})$, and then $I_f=\mathcal O(X_f)_{\lambda<0}$. Let $J_f\subseteq \mathcal O(X_f)$ be the kernel of the composition 
        \begin{equation}
            \mathcal O(X_f)\to\mathcal O(U)\otimes\mathcal O(X_f)\to\mathcal O(U)\otimes\big(\mathcal O(X_f)/I_f\big) 
        \end{equation}
        where $\mathcal O(X_f)\to\mathcal O(U)\otimes \mathcal O(X_f)$ is the co-action map of $U\curvearrowright X_f$. Let $C_f\hookrightarrow X_f$ be the closed subscheme associated to the ideal $J_f\subseteq\mathcal O(X_f)$. The ideals $J_f\subseteq\mathcal O(X_f)$ for $f\in H^0(X,L^m)_{\lambda=\max}$ glue to a sheaf of ideals $\mathcal J\subseteq\mathcal O_{X^0_{\lambda-\min}}$. The closed subschemes $C_f\hookrightarrow X_f$ glue to a closed subscheme $C\hookrightarrow X^0_{\lambda-\min}$, which is the closed subscheme associated to $\mathcal J$. Let $\overline C\hookrightarrow X$ be the scheme theoretic image of $C\hookrightarrow X^0_{\lambda-\min}\subseteq X$, which is the blow-up centre. 
        {
            \begin{remark}
                Let $Z\hookrightarrow X^0_{\lambda-\min}$ be the scheme theoretic intersection of $Z_{\lambda-\min}\hookrightarrow X^0_{\lambda-\min}$ and $\mathbb V(\mathcal I_{k_1,\cdots,k_n})\hookrightarrow X^0_{\lambda-\min}$, where the subscheme structure on $\mathbb V(\mathcal I_{k_1,\cdots,k_n})$ is associated to $\mathcal I_{k_1,\cdots,k_n}$. Consider the morphism $U\times Z\to X^0_{\lambda-\min}$, sending a $\Bbbk$-point $(u,z)$ to $u.z\in X^0_{\lambda-\min}$. Then $C\hookrightarrow X^0_{\lambda-\min}$ is the scheme theoretic image of this morphism. By \cite[\href{https://stacks.math.columbia.edu/tag/01R8}{Lemma 01R8}]{stacks-project}, the closed subset of $C$ is the closure of the image, i.e. the $U$-sweep of $Z$. Also by \cite[\href{https://stacks.math.columbia.edu/tag/01R8}{Lemma 01R8}]{stacks-project}, scheme theoretic images can be constructed locally. Then the scheme theoretic image of $U\times Z_f\to X_f$ is $C_f\hookrightarrow X_f$. 
            \end{remark}
        }
        
        Let $\pi:\tilde X\to X$ be the blow-up of $X$ along $\overline C\hookrightarrow X$. Let $E\hookrightarrow X$ be the exceptional divisor. For any $x,y\in\mathbb N_+$, the line bundle $\pi^*L^x\otimes \mathcal O(-yE)$ is ample on $\tilde X$. Naturally $\hat U $ acts on $\tilde X$ linearly with respect to $\pi^*L^x\otimes\mathcal O(-yE)$. Then $(\tilde X)^0_{\lambda-\min}$ is defined. Blow-ups are constructed locally. Define the \emph{blowup algebra associated to $J_f\subseteq\mathcal O(X_f)$} 
        \begin{equation}
            \mathrm{Bl}_{J_f}\big(\mathcal O(X_f)\big):=\bigoplus_{k\in\mathbb N}(J_f)^k
        \end{equation}
        where $(J_f)^k$ has degree $k$. For a degree one element $a\in J_f$, define the \emph{affine blowup algebra} as the homogeneous localisation at $a$ 
        \begin{equation}
            \mathcal O(X_f)\Big[\frac{J_f}{a}\Big]:=\Big(\mathrm{Bl}_{J_f}\big(\mathcal O(X_f)\big)\Big)_{(a)}. 
        \end{equation}
        
        Let $\tilde X_f\to X_f$ be the blow-up of $X_f$ along $C_f\hookrightarrow X_f$. It has an open affine covering 
        \begin{equation}
            \tilde X_f=\bigcup_{a\in J_f}\tilde X_{f,a},\quad \tilde X_{f,a}:=\mathrm{Spec}\Big(\mathcal O(X_f)\Big[\frac{J_f}{a}\Big]\Big). 
        \end{equation}
        
        The open subscheme $(\tilde X)^0_{\min}\subseteq \tilde X$ is the union of the following affine open subschemes 
        \begin{equation}
            (\tilde X)^0_{\lambda-\min}=\bigcup_{\substack{f\in H^0(X,L^m)_{\lambda=\max}\\a\in (J_f)_{\lambda=0}}}\tilde X_{f,a}. 
        \end{equation}
        The open subscheme $(\tilde X)^0_{\lambda-\min}\subseteq\tilde X$ is covered by those $\tilde X_{f,a}$ such that $f\in H^0(X,L^m)_{\lambda=\max}$ with $Z_{\lambda-\min,f}\not\subseteq\mathbb V(\mathcal I_{k_1,\cdots,k_n})$ and $a\in (J_f)_{\lambda=0}$ of the form $a=a^{(1)}\cdots a^{(n)}$ for $a^{(i)}\in \mathrm{Fit}_{k_i}(\phi_i)(X_f)_{\lambda=0}$. We focus on one such $\tilde X_{f,a}$. For each $1\leq i\leq n$ choose a suitable basis $\big(\xi^{(i)}_1,\cdots,\xi^{(i)}_{r_i}\big)$ of $\mathfrak u_{\lambda=w_i}$ and assume without loss of generality that 
        \begin{equation}
            a^{(i)}=\det\begin{pmatrix}\xi^{(i)}_1.f^{(i)}_1&\cdots&\xi^{(i)}_1.f^{(i)}_{r_i-k_i}\\
            \vdots&\ddots&\vdots\\
            \xi^{(i)}_{r_i-k_i}.f^{(i)}_1&\cdots&\xi^{(i)}_{r_i-k_i}.f^{(i)}_{r_i-k_i}\end{pmatrix}\in \mathrm{Fit}_{k_i}(\phi_i)(X_f)_{\lambda=0}
        \end{equation}
        where $f^{(i)}_j\in\mathcal O(X_f)_{\lambda=-w_i}$. 

        Let $\mathrm{U}(\mathfrak u)$ denote the \emph{universal enveloping algebra of $\mathfrak u$}. For $1\leq i\leq n$ and $1\leq \mu_i\leq r_i-k_i$, define a linear map 
        \begin{equation}
            \begin{split}
                E^{(i)}_{\mu_i}:\mathrm{U}(\mathfrak u)&\to \mathcal O(X_f)\\
                A&\mapsto \det\begin{pmatrix}\xi^{(i)}_1.f^{(i)}_1&\cdots&\xi^{(i)}_1.f^{(i)}_{r_i-k_i}\\
                \vdots&\cdots&\vdots\\
                A.f^{(i)}_1&\cdots&A.f^{(i)}_{r_i-k_i}\\
                \vdots&\cdots&\vdots\\
                \xi^{(i)}_{r_i-k_i}.f^{(i)}_1&\cdots&\xi^{(i)}_{r_i-k_i}.f^{(i)}_{r_i-k_i}\end{pmatrix}
            \end{split}
        \end{equation}
        where $\big(A.f^{(i)}_1,\cdots,A.f^{(i)}_{r_i-k_i}\big)$ occupies the $\mu_i$th row. For $A,B\in \mathrm{U}(\mathfrak u)$, we have $B.E^{(i)}_{\mu_i}(A)=E^{(i)}_{\mu_i}(BA)$, since all entries in the matrix above are in $\mathcal O(X_f)_{\lambda=0}\subseteq\mathcal O(X_f)^U$ except the $\mu_i$th row. 
        
        {
            \begin{lemma}\label{lemma of determinantal sum}
                Let $\tilde X_{f,a},a^{(i)},\xi^{(i)}_j$ be as above. Let $1\leq i\leq n$ and $h\in\mathcal O(X_f)_{\lambda=-w_i}$ and $A\in\mathfrak u_{\lambda=w_i}$. We have the following equation
                \begin{equation}
                    \sum_{\mu_i=1}^{r_i-k_i}\big(\xi^{(i)}_{\mu_i}.h\big) E^{(i)}_{\mu_i}(A)=\big(A.h\big)a^{(i)}.
                \end{equation}
            \end{lemma}
            \begin{proof}
                Since $f^{(i)}_j,h\in\mathcal O(X_f)_{\lambda=-w_i}$, we have 
                \begin{equation}
                    \mathrm{d}f^{(i)}_1,\cdots,\mathrm{d}f^{(i)}_{r_i-k_i},\mathrm{d}h\in \mathcal K_{i-1}(X_f). 
                \end{equation}
                Then the determinant of the following $(r_i-k_i+1)\times(r_i-k_i+1)$ matrix is in $\mathrm{Fit}_{k_i-1}(\phi_i)$ 
                \begin{equation}
                    M=\begin{pmatrix}\xi^{(i)}_1.f^{(i)}_1&\cdots&\xi^{(i)}_1.f^{(i)}_{r_i-k_i}&\xi^{(i)}_1.h\\
                    \vdots&\ddots&\vdots&\vdots\\
                    \xi^{(i)}_{r_i-k_i}.f^{(i)}_1&\cdots&\xi^{(i)}_{r_i-k_i}.f^{(i)}_{r_i-k_i}&\xi^{(i)}_{r_i-k_i}.h\\
                    A.f^{(i)}_1&\cdots&A.f^{(i)}_{r_i-k_i}&A.h
                    \end{pmatrix}
                \end{equation}
                
                Expand $\det(M)$ along the right column 
                \begin{equation}
                    \det(M)=(A.h\big)a^{(i)}-\sum_{\mu_i=1}^{r_i-k_i}\big(\xi^{(i)}_{\mu_i}.h\big)E^{(i)}_{\mu_i}(A). 
                \end{equation}

                By the definition of $k_i$ we have $\det(M)=0$. This completes the proof. 
            \end{proof}
        }

        {
            \begin{lemma}\label{lemma of existence of b^(i)_nu}
                Let $\tilde X_{f,a},a^{(i)},\xi^{(i)}_j,I_f,J_f$ be as above. Then for $1\leq i\leq n$ and $1\leq \mu\leq r_i-k_i$, there exists $b^{(i)}_\mu\in\mathcal O(X_f)_{\lambda=-w_i}$ satisfying: 
                \begin{itemize}
                    \item[(1)] $\xi^{(i)}_\mu.b^{(i)}_\nu=w_i\delta_{\mu\nu}\prod_{i'=i}^na^{(i')}$ for $1\leq i\leq n$ and $1\leq \mu,\nu\leq r_i-k_i$; 
                    \item[(2)] $A.b^{(i)}_\mu\in\prod_{i'=i}^n\mathrm{Fit}_{k_{i'}}(\phi_{i'})(X_f)_{\lambda=0}$ for all $1\leq i\leq n$ and $A\in\mathrm{U}(\mathfrak u)_{\lambda=-w_i}$; 
                    \item[(3)] $\Big(\prod_{i'<i}a^{(i')}\Big)\cdot b^{(i)}_\mu\in J_f$ for $1\leq i\leq n$ and $1\leq \mu\leq r_i-k_i$. 
                \end{itemize}
            \end{lemma}
            \begin{proof}
                If $(2)$ is proved, then $A.\Big(\prod_{i'<i}a^{(i')}\cdot b^{(i)}_\mu\Big)\in I_f$ for all $A\in\mathrm{U}(\mathfrak u)$. By Lemma \ref{corollary of the ideal cutting out U-sweep of Z consists functions whose UEA derivations in the ideal of Z}, it follows that $\prod_{i'<i}a^{(i')}\cdot b^{(i)}_\mu\in J_f$, which is $(3)$. 
                
                For $1\leq \mu_n\leq r_n-k_n$, we regard the smallest grading weight $w_n\in\mathbb N_+$ as a scalar $w_n\in \mathrm{U}(\mathfrak u)$ and define 
                \begin{equation}
                    b^{(n)}_{\mu_n}:=E^{(n)}_{\mu_n}(w_n). 
                \end{equation}
                
                We construct $b^{(i)}_{\mu_i}$ inductively for $i=n,\cdots,1$. Assume $b^{(i')}_{\mu_{i'}}$ is defined for $i<i'\leq n$ and $1\leq \mu_{i'}\leq r_{i'}-k_{i'}$. Define  $b^{(i)}_{\mu_i}$ as 
                \begin{equation}
                    \begin{split}
                        b^{(i)}_{\mu_i}:=E^{(i)}_{\mu_i}(w_i)\cdot\prod_{i'=i+1}^na^{(i')}-\sum_{i'=i+1}^n\sum_{\mu_{i'}=1}^{r_{i'}-k_{i'}}E^{(i)}_{\mu_i}\big(\xi^{(i')}_{\mu_{i'}}\big)\cdot \prod_{i<i''<i'}a^{(i'')}\cdot b^{(i')}_{\mu_{i'}}. 
                    \end{split}
                \end{equation}
                We see $b^{(i)}_{\mu_i}\in\mathcal O(X_f)_{\lambda=-w_i}$ and $(1)$ is checked immediately. 

                By the Poincar\'e-Birkhoff-Witt theorem, the following elements form a basis of $\mathrm{U}(\mathfrak u)$ 
                \begin{equation}
                    \big(\xi^{(1)}_1\big)^{p^{(1)}_1}\cdots\big( \xi^{(1)}_{r_1}\big)^{p^{(1)}_{r_1}}\cdots \big(\xi^{(n)}_1\big)^{p^{(n)}_1}\cdots\big( \xi^{(n)}_{r_n}\big)^{p^{(n)}_{r_n}}\in\mathrm{U}(\mathfrak u),\quad p^{(i)}_j\in\mathbb N. 
                \end{equation}
                For simplicity denote 
                \begin{equation}
                    \begin{split}
                        \xi^p:=&\big(\xi^{(1)}\big)^{p^{(1)}}\cdots\big(\xi^{(n)}\big)^{p^{(n)}},\quad\big(\xi^{(i)}\big)^{p^{(i)}}:=\big(\xi^{(i)}_1\big)^{p^{(i)}_1}\cdots\big(\xi^{(i)}_{r_i}\big)^{p^{(i)}_{r_i}}\\
                        p:=&\big(p^{(1)},\cdots,p^{(n)}\big),\quad p^{(i)}:=\big(p^{(i)}_1,\cdots,p^{(i)}_{r_i}\big)\in\mathbb N^{r_i},\quad|p^{(i)}|:=p^{(i)}_1+\cdots+p^{(i)}_{r_i}. 
                    \end{split}
                \end{equation}
                For $(2)$ it suffices to consider $A=\xi^p\in\mathrm{U}(\mathfrak u)_{\lambda=-w_i}$. 
                
                We construct a family of functions $\beta^{(i)}_{\mu_i}(\xi^p)\in \mathcal O(X_f)_{\lambda=0}$ for $1\leq i\leq n$ and $1\leq \mu_i\leq r_i-k_i$ and $\xi^p\in\mathrm{U}(\mathfrak u)_{\lambda=w_i}$ inductively by 
                \begin{equation}
                    \begin{split}
                        \beta^{(n)}_{\mu_n}(\xi^p):=&w_{\max(p)}E^{(n)}_{\mu_n}\big([\xi^p]\!]\big)\\
                        \beta^{(i)}_{\mu_i}(\xi^p):=&w_{\max(p)}E^{(i)}_{\mu_i}\big([\xi^p]\!]\big)\cdot\prod_{i'=i+1}^na^{(i')}\\
                        &-\sum_{i'=i+1}^n\sum_{\substack{0\leq q\leq p\\\mathrm{weight}(\xi^q)=w_{i'}}}\sum_{\mu_{i'}=1}^{r_{i'}-k_{i'}}\binom{p}{q}E^{(i)}_{\mu_i}([\xi^{p-q}\xi^{(i')}_{\mu_{i'}}]\!])\cdot\prod_{i<i''<i'}a^{(i'')}\cdot\beta^{(i')}_{\mu_{i'}}(\xi^q)
                    \end{split}
                \end{equation}
                where $w_{\max(p)}:=w_{\max\{j:p^{(j)}\ne0\}}$ and $[-]\!]:\mathrm{U}(\mathfrak u)\to\mathfrak u$ is the complete bracket defined in Appendix \ref{appendix of complete brackets} respecting $\lambda$-weights. We inductively see that $\beta^{(i)}_{\mu_i}(\xi^p)\in \prod_{i'=i}^n\mathrm{Fit}_{k_{i'}}(\phi_{i'})(X_f)_{\lambda=0}$. We then prove for $1\leq i\leq n$ and $1\leq \mu_i\leq r_i-k_i$ and $\xi^p\in\mathrm{U}(\mathfrak u)_{\lambda=w_i}$ 
                \begin{equation}\label{equation of UEA on elements and complete bracket}
                    \xi^p.b^{(i)}_{\mu_i}= \beta^{(i)}_{\mu_i}(\xi^p). 
                \end{equation}

                Let $1\leq i<n$, and assume equation \eqref{equation of UEA on elements and complete bracket} is proved for $i'>i$. Let $\xi^p\in\mathrm{U}(\mathfrak u)_{\lambda=w_i}$, i.e. $\sum_{i'=1}^nw_{i'}|p^{(i')}|=w_i$. We calculate $\xi^p.b^{(i)}_{\mu_i}$ 
                \begin{equation}
                    \begin{split}
                        \xi^p.b^{(i)}_{\mu_i}=&\xi^p.\bigg(E^{(i)}_{\mu_i}(w_i)\cdot\prod_{i'=i+1}^na^{(i')}-\sum_{i'=i+1}^n\sum_{\mu_{i'}=1}^{r_{i'}-k_{i'}}E^{(i)}_{\mu_i}\big(\xi^{(i')}_{\mu_{i'}}\big)\cdot \prod_{i<i''<i'}a^{(i'')}\cdot b^{(i')}_{\mu_{i'}}\bigg)\\
                        =&E^{(i)}_{\mu_i}(w_i\xi^p)\cdot\prod_{i'>i}a^{(i')}-\xi^p.\bigg(\sum_{i'=i+1}^n\sum_{\mu_{i'}=1}^{r_{i'}-k_{i'}}E^{(i)}_{\mu_i}\big(\xi^{(i')}_{\mu_{i'}}\big)\cdot \prod_{i<i''<i'}a^{(i'')}\cdot b^{(i')}_{\mu_{i'}}\bigg)
                    \end{split}
                \end{equation}
                where we have used that $a^{(i)}\in\mathcal O(X_f)_{\lambda=0}$ is $U$-invariant, and $A.E^{(i)}_{\mu_i}(B)=E^{(i)}_{\mu_i}(AB)$ for all $A,B\in\mathrm{U}(\mathfrak u)$. Write $\xi^p.b^{(i)}_{\mu_i}=T_1-T_2$ for
                \begin{equation}
                    \begin{split}
                        T_1:=&E^{(i)}_{\mu_i}(w_i\xi^p)\cdot\prod_{i'>i}a^{(i')}\\
                        T_2:=&\xi^p.\bigg(\sum_{i'=i+1}^n\sum_{\mu_{i'}=1}^{r_{i'}-k_{i'}}E^{(i)}_{\mu_i}\big(\xi^{(i')}_{\mu_{i'}}\big)\cdot \prod_{i<i''<i'}a^{(i'')}\cdot b^{(i')}_{\mu_{i'}}\bigg). 
                    \end{split}
                \end{equation}
                
                For the term $T_1$ we have by Lemma \ref{lemma of identity of weighted complete bracket}
                \begin{equation}
                    \begin{split}
                        T_1\equiv&E^{(i)}_{\mu_i}(w_i\xi^p)\cdot\prod_{i'>i}a^{(i')}\\
                        =&\sum_{0<q\leq p}\binom{p}{q}w_{\max(q)}E^{(i)}_{\mu_i}\big([\xi^q]\!]\xi^{p-q}\big)\cdot\prod_{i'>i}a^{(i')}\\
                        =&w_{\max(p)}E^{(i)}_{\mu_i}\big([\xi^p]\!]\big)\cdot\prod_{i'>i}a^{(i')}+\sum_{0<q<p}\binom{p}{q}w_{\max(q)}E^{(i)}_{\mu_i}
                        \big([\xi^q]\!]\xi^{p-q}
                        \big)\cdot\prod_{i'>i}a^{(i')}\\
                        =&w_{\max(p)}E^{(i)}_{\mu_i}\big([\xi^p]\!]\big)\cdot\prod_{i'>i}a^{(i')}\\&+\sum_{i'=i+1}^n\sum_{\substack{0\leq q\leq p\\\mathrm{weight}(\xi^q)=w_{i'}}}\binom{p}{q}w_{\max(q)}E^{(i)}_{\mu_i}\big([\xi^q]\!]\xi^{p-q}\big)\cdot\prod_{i''>i}a^{(i'')}\\
                        =&T_{1,1}+T_{1,2}
                    \end{split}
                \end{equation}
                where 
                \begin{equation}
                    \begin{split}
                       T_{1,1}:=&w_{\max(p)}E^{(i)}_{\mu_i}\big([\xi^p]\!]\big)\cdot\prod_{i'>i}a^{(i')}\\ T_{1,2}:=&\sum_{\substack{i'>i\\0\leq q\leq p\\\mathrm{weight}(\xi^q)=w_{i'}}}\binom{p}{q}w_{\max(q)}E^{(i)}_{\mu_i}\big([\xi^q]\!]\xi^{p-q}\big)\cdot\prod_{i''>i}a^{(i'')}. 
                    \end{split}
                \end{equation}
                
                For the term $T_2$ we have 
                \begin{equation}
                    \begin{split}
                        T_2\equiv&\xi^p.\bigg(\sum_{i'=i+1}^n\sum_{\mu_{i'}=1}^{r_{i'}-k_{i'}}E^{(i)}_{\mu_i}\big(\xi^{(i')}_{\mu_{i'}}\big)\cdot \prod_{i<i''<i'}a^{(i'')}\cdot b^{(i')}_{\mu_{i'}}\bigg)\\
                        =&\sum_{i'=i+1}^n\sum_{0\leq q\leq p}\sum_{\mu_{i'}=1}^{r_{i'}-k_{i'}}\binom{p}{q}E^{(i)}_{\mu_i}\big(\xi^{p-q}\xi^{(i')}_{\mu_{i'}}\big)\cdot \prod_{i<i''<i'}a^{(i'')}\cdot \beta^{(i')}_{\mu_{i'}}(\xi^q)\\
                        =&\sum_{i'=i+1}^n\sum_{\substack{0\leq q\leq p\\\mathrm{weight}(\xi^q)=w_{i'}}}\sum_{\mu_{i'}=1}^{r_{i'}-k_{i'}}\binom{p}{q}E^{(i)}_{\mu_i}\big(\xi^{p-q}\xi^{(i')}_{\mu_{i'}}\big)\cdot \prod_{i<i''<i'}a^{(i'')}\cdot \beta^{(i')}_{\mu_{i'}}(\xi^q)\\
                        =&\sum_{\substack{i'>i\\0\leq q\leq p\\\mathrm{weight}(\xi^q)=w_{i'}}}\sum_{\mu_{i'}=1}^{r_{i'}-k_{i'}}\binom{p}{q}E^{(i)}_{\mu_i}\bigg(\sum_{0\leq r\leq p-q}\binom{p-q}{r}[\xi^{p-q-r}\xi^{(i')}_{\mu_{i'}}]\!]\xi^r\bigg)\cdot \prod_{i<i''<i'}a^{(i'')}\cdot \beta^{(i')}_{\mu_{i'}}(\xi^q)\\
                        =&\sum_{\substack{i'>i\\0\leq q\leq p\\\mathrm{weight}(\xi^q)=w_{i'}}}\sum_{\mu_{i'}=1}^{r_{i'}-k_{i'}}\binom{p}{q}E^{(i)}_{\mu_i}\big([\xi^{p-q}\xi^{(i')}_{\mu_{i'}}]\!]\big)\cdot \prod_{i<i''<i'}a^{(i'')}\cdot \beta^{(i')}_{\mu_{i'}}(\xi^q)\\
                        &+\sum_{\substack{i'>i\\0\leq q\leq p\\\mathrm{weight}(\xi^q)=w_{i'}}}\sum_{\mu_{i'}=1}^{r_{i'}-k_{i'}}\binom{p}{q}E^{(i)}_{\mu_i}\big(\xi^{(i')}_{\mu_{i'}}\xi^{p-q}\big)\cdot \prod_{i<i''<i'}a^{(i'')}\cdot \beta^{(i')}_{\mu_{i'}}(\xi^q)\\
                        +&\sum_{\substack{i'>i\\0\leq q\leq p\\\mathrm{weight}(\xi^q)=w_{i'}}}\sum_{\mu_{i'}=1}^{r_{i'}-k_{i'}}\binom{p}{q}E^{(i)}_{\mu_i}\bigg(\sum_{0< r< p-q}\binom{p-q}{r}[\xi^{p-q-r}\xi^{(i')}_{\mu_{i'}}]\!]\xi^r\bigg)\cdot \prod_{i<i''<i'}a^{(i'')}\cdot \beta^{(i')}_{\mu_{i'}}(\xi^q)
                    \end{split}
                \end{equation}
                where we have used: 
                \begin{itemize}
                    \item the generalised Leibniz rule 
                    \begin{equation}
                        \xi^p.(fg)=\sum_{0\leq q\leq p}\binom{p}{q}(\xi^{p-q}.f)(\xi^q.g); 
                    \end{equation}
                    \item the induction hypothesis for $i'>i$
                    \begin{equation}
                        \xi^q.b^{(i')}_{\mu_{i'}}=\beta^{(i')}_{\mu_{i'}}(\xi^q); 
                    \end{equation}
                    \item if $\mathrm{weight}(\xi^q)>w_{i'}$, then $\xi^q.b^{(i')}_{\mu_{i'}}=0$, since $b^{(i')}_{\mu_{i'}}\in\mathcal O(X_f)_{\lambda=-w_{i'}}$; 
                    \item if $\mathrm{weight}(\xi^q)<w_{i'}$, then $\mathrm{weight}(\xi^{p-q})>w_i-w_{i'}$, and then $E^{(i)}_{\mu_i}\big(\xi^{p-q}\xi^{(i')}_{\mu_{i'}}\big)=0$, since $E^{(i)}_{\mu_i}\big(\xi^{(i')}_{\mu_{i'}}\big)\in\mathcal O(X_f)_{-w_i+w_{i'}}$; 
                    \item the identity in Lemma \ref{lemma of commutator identity with complete bracket}
                    \begin{equation}
                        \xi^{p-q}\xi^{(i')}_{\mu_{i'}}=\sum_{0\leq r\leq p-q}\binom{p-q}{r}[\xi^{p-q-r}\xi^{(i')}_{\mu_{i'}}]\!]\xi^r. 
                    \end{equation}
                \end{itemize}

                We write $T_2=T_{2,1}+T_{2,2}+T_{2,3}$ for 
                \begin{equation}
                    \begin{split}
                        T_{2,1}:=&\sum_{\substack{i'>i\\0\leq q\leq p\\\mathrm{weight}(\xi^q)=w_{i'}}}\sum_{\mu_{i'}=1}^{r_{i'}-k_{i'}}\binom{p}{q}E^{(i)}_{\mu_i}\big([\xi^{p-q}\xi^{(i')}_{\mu_{i'}}]\!]\big)\cdot \prod_{i<i''<i'}a^{(i'')}\cdot \beta^{(i')}_{\mu_{i'}}(\xi^q)\\
                        T_{2,2}:=&\sum_{\substack{i'>i\\0\leq q\leq p\\\mathrm{weight}(\xi^q)=w_{i'}}}\sum_{\mu_{i'}=1}^{r_{i'}-k_{i'}}\binom{p}{q}E^{(i)}_{\mu_i}\big(\xi^{(i')}_{\mu_{i'}}\xi^{p-q}\big)\cdot \prod_{i<i''<i'}a^{(i'')}\cdot \beta^{(i')}_{\mu_{i'}}(\xi^q)\\
                        T_{2,3}:=&\\
                        \sum_{\substack{i'>i\\0\leq q\leq p\\\mathrm{weight}(\xi^q)=w_{i'}}}&\sum_{\mu_{i'}=1}^{r_{i'}-k_{i'}}\binom{p}{q}E^{(i)}_{\mu_i}\bigg(\sum_{0< r< p-q}\binom{p-q}{r}[\xi^{p-q-r}\xi^{(i')}_{\mu_{i'}}]\!]\xi^r\bigg)\cdot \prod_{i<i''<i'}a^{(i'')}\cdot \beta^{(i')}_{\mu_{i'}}(\xi^q). 
                    \end{split}
                \end{equation}
                
                Observe that $\beta^{(i)}_{\mu_i}(\xi^p)=T_{1,1}-T_{2,1}$. Then equation \eqref{equation of UEA on elements and complete bracket} is equivalent to $T_{1,2}-T_{2,2}-T_{2,3}=0$. For $T_{2,3}$ by a change of indices we have 
                \begin{equation}
                    \begin{split}
                        T_{2,3}\equiv&\sum_{\substack{i'>i\\0\leq q\leq p\\\mathrm{weight}(\xi^q)=w_{i'}}}\sum_{\mu_{i'}=1}^{r_{i'}-k_{i'}}\binom{p}{q}E^{(i)}_{\mu_i}\bigg(\sum_{0< r< p-q}\binom{p-q}{r}[\xi^{p-q-r}\xi^{(i')}_{\mu_{i'}}]\!]\xi^r\bigg)\cdot \prod_{i<i''<i'}a^{(i'')}\cdot \beta^{(i')}_{\mu_{i'}}(\xi^q)\\
                        =&\sum_{\substack{i'>i\\0\leq q\leq p\\\mathrm{weight}(\xi^q)=w_{i'}}}\sum_{q<s<p}\sum_{\mu_{i'}=1}^{r_{i'}-k_{i'}}\binom{p}{q}\binom{p-q}{p-s}E^{(i)}_{\mu_i}\big([\xi^{s-q}\xi^{(i')}_{\mu_{i'}}]\!]\xi^{p-s}\big)\cdot \prod_{i<i''<i'}a^{(i'')}\cdot \beta^{(i')}_{\mu_{i'}}(\xi^q)\\
                        =&\sum_{\substack{i'>i\\0\leq q\leq p\\\mathrm{weight}(\xi^q)=w_{i'}}}\sum_{\substack{i<i''<i'\\q\leq s\leq p\\\mathrm{weight}(\xi^s)=w_{i''}}}\sum_{\mu_{i'}=1}^{r_{i'}-k_{i'}}\binom{p}{q}\binom{p-q}{p-s}E^{(i)}_{\mu_i}\big([\xi^{s-q}\xi^{(i')}_{\mu_{i'}}]\!]\xi^{p-s}\big)\cdot \prod_{i<i''<i'}a^{(i'')}\cdot \beta^{(i')}_{\mu_{i'}}(\xi^q)\\
                        =&\sum_{\substack{i''>i\\0\leq s\leq p\\\mathrm{weight}(\xi^s)=w_{i''}}}\sum_{\substack{i'>i''\\0\leq q\leq s\\\mathrm{weight}(\xi^q)=w_{i'}}}\sum_{\mu_{i'}=1}^{r_{i'}-k_{i'}}\binom{p}{q}\binom{p-q}{p-s}E^{(i)}_{\mu_i}\big([\xi^{s-q}\xi^{(i')}_{\mu_{i'}}]\!]\xi^{p-s}\big)\cdot \prod_{i<i''<i'}a^{(i'')}\cdot \beta^{(i')}_{\mu_{i'}}(\xi^q)\\
                        =&\sum_{\substack{i'>i\\0\leq q\leq p\\\mathrm{weight}(\xi^q)=w_{i'}}}\sum_{\substack{i''>i'\\0\leq r\leq q\\\mathrm{weight}(\xi^r)=w_{i''}}}\sum_{\mu_{i''}=1}^{r_{i''}-k_{i''}}\binom{p}{q}\binom{q}{r}E^{(i)}_{\mu_i}\big([\xi^{q-r}\xi^{(i'')}_{\mu_{i''}}]\!]\xi^{p-q}\big)\cdot\prod_{i<i'''<i''}a^{(i''')}\cdot \beta^{(i'')}_{\mu_{i''}}(\xi^r). 
                    \end{split}
                \end{equation}
                
                For $T_{2,2}$ by the inductive definition of $\beta^{(i')}_{\mu_{i'}}(\xi^q)$, we have 
                \begin{equation}
                    \begin{split}
                        &T_{2,2}\\
                        \equiv&\sum_{\substack{i'>i\\0\leq q\leq p\\\mathrm{weight}(\xi^q)=w_{i'}}}\sum_{\mu_{i'}=1}^{r_{i'}-k_{i'}}\binom{p}{q}E^{(i)}_{\mu_i}\big(\xi^{(i')}_{\mu_{i'}}\xi^{p-q}\big)\cdot \prod_{i<i''<i'}a^{(i'')}\cdot\beta^{(i')}_{\mu_{i'}}(\xi^q)\\
                        =&\sum_{\substack{i'>i\\0\leq q\leq p\\\mathrm{weight}(\xi^q)=w_{i'}}}\sum_{\mu_{i'}=1}^{r_{i'}-k_{i'}}\binom{p}{q}E^{(i)}_{\mu_i}\big(\xi^{(i')}_{\mu_{i'}}\xi^{p-q}\big)\cdot \prod_{i<i''<i'}a^{(i'')}\\
                        &\times\begin{pmatrix}w_{\max(q)}E^{(i')}_{\mu_{i'}}\big([\xi^q]\!]\big)\cdot\prod_{i''>i'}a^{(i'')}\\-\sum_{i''>i'}\sum_{\substack{0\leq r\leq q\\\mathrm{weight}(\xi^r)=w_{i''}}}\sum_{\mu_{i''}=1}^{r_{i''}-k_{i''}}\binom{q}{r}E^{(i')}_{\mu_{i'}}\big([\xi^{q-r}\xi^{(i'')}_{\mu_{i''}}]\!]\big)\cdot\prod_{i'<i'''<i''}a^{(i''')}\cdot\beta^{(i'')}_{\mu_{i''}}(\xi^r)\end{pmatrix}\\
                        =&\sum_{\substack{i'>i\\0\leq q\leq p\\\mathrm{weight}(\xi^q)=w_{i'}}}\sum_{\mu_{i'}=1}^{r_{i'}-k_{i'}}\binom{p}{q}w_{\max(q)}E^{(i)}_{\mu_i}\big(\xi^{(i')}_{\mu_{i'}}\xi^{p-q}\big) E^{(i')}_{\mu_{i'}}\big([\xi^q]\!]\big)\cdot\prod_{\substack{i''>i\\i''\ne i'}}a^{(i'')}\\
                        &-\sum_{\substack{i'>i\\0\leq q\leq p\\\mathrm{weight}(\xi^q)=w_{i'}}}\sum_{\mu_{i'}=1}^{r_{i'}-k_{i'}}\sum_{\substack{i''>i'\\0\leq r\leq q\\\mathrm{weight}(\xi^r)=w_{i''}}}\sum_{\mu_{i''}=1}^{r_{i''}-k_{i''}}\begin{pmatrix}\binom{p}{q}\binom{q}{r}E^{(i)}_{\mu_i}\big(\xi^{(i')}_{\mu_{i'}}\xi^{p-q}\big)E^{(i')}_{\mu_{i'}}\big([\xi^{q-r}\xi^{(i'')}_{\mu_{i''}}]\!]\big)\\
                        \times \prod_{\substack{i<i'''<i''\\i'''\ne i'}}a^{(i''')}\cdot \beta^{(i'')}_{\mu_{i''}}(\xi^r)\end{pmatrix}. 
                    \end{split}
                \end{equation}
                
                By Lemma \ref{lemma of determinantal sum}, we have the following equations
                \begin{equation}
                    \begin{split}
                        \sum_{\mu_{i'}=1}^{r_{i'}-k_{i'}}E^{(i)}_{\mu_i}\big(\xi^{(i')}_{\mu_{i'}}\xi^{p-q}\big)E^{(i')}_{\mu_{i'}}\big([\xi^q]\!]\big)&=E^{(i)}_{\mu_i}\big([\xi^q]\!]\xi^{p-q}\big)a^{(i')}\\
                        \sum_{\mu_{i'}=1}^{r_{i'}-k_{i'}}E^{(i)}_{\mu_i}\big(\xi^{(i')}_{\mu_{i'}}\xi^{p-q}\big)E^{(i')}_{\mu_{i'}}\big([\xi^{q-r}\xi^{(i'')}_{\mu_{i''}}]\!]\big)&=E^{(i)}_{\mu_i}\big([\xi^{q-r}\xi^{(i'')}_{\mu_{i''}}]\!]\xi^{p-q}\big)a^{(i')}. 
                    \end{split}
                \end{equation}
                Then for $T_{2,2}$ we have 
                \begin{equation}
                    \begin{split}
                        T_{2,2}=&\sum_{\substack{i'>i\\0\leq q\leq p\\\mathrm{weight}(\xi^q)=w_{i'}}}\binom{p}{q}w_{\max(q)}E^{(i)}_{\mu_i}\big([\xi^q]\!]\xi^{p-q}\big)\cdot\prod_{i''>i}a^{(i'')}\\
                        &-\sum_{\substack{i'>i\\0\leq q\leq p\\\mathrm{weight}(\xi^q)=w_{i'}}}\sum_{\substack{i''>i'\\0\leq r\leq q\\\mathrm{weight}(\xi^r)=w_{i''}}}\sum_{\mu_{i''}=1}^{r_{i''}-k_{i''}}\binom{p}{q}\binom{q}{r}E^{(i)}_{\mu_i}\big([\xi^{q-r}\xi^{(i'')}_{\mu_{i''}}]\!]\xi^{p-q}\big)\cdot\prod_{i<i'''<i''}a^{(i''')}\cdot \beta^{(i'')}_{\mu_{i''}}(\xi^r)\\
                        =&T_{1,2}-T_{2,3}. 
                    \end{split}
                \end{equation}
                This completes the proof. 
            \end{proof}
        }
    }
    
    \subsection{Statement and proof of the main theorem}
    {
        We will state and prove Theorem \ref{theorem of blow-up for CDRS} and Corollary \ref{corollary of consequence of blow-up}. Let $\hat U \curvearrowright (X,L)$ with grading weights $w_1>\cdots>w_n$. Let $\{\mathfrak u_i\}_{i=0}^n$ be the filtration $\mathfrak u_i:=\mathfrak u_{\lambda\geq w_i}$. On $X^0_{\lambda-\min}$, consider coherent sheaves and morphisms 
        \begin{equation}
            \begin{split}
                &\mathcal K_i=\ker\Big(\Omega_{X^0_{\lambda-\min}/\Bbbk}\to (\mathfrak u_i)^*\otimes \mathcal O_{X^0_{\lambda-\min}}\Big)\\
                &\phi_i:\mathcal K_{i-1}\to (\mathfrak u_i/\mathfrak u_{i-1})^*\otimes\mathcal O_{X^0_{\lambda-\min}}\\
                &\mathrm{Fit}_k(\phi_i):=\mathrm{Fit}_k\big(\mathrm{coker}(\phi_i)\big). 
            \end{split}
        \end{equation}
        
        Let $(k_1,\cdots,k_n)\in\mathbb N^n$ be such that 
        \begin{equation}
            k_i:=\min\big\{k\in\mathbb Z:\mathrm{Fit}_k(\phi_i)\ne0\big\},\quad i=1,\cdots,n.
        \end{equation}
        and let $\mathcal I_{k_1,\cdots,k_n}:=\prod_{i=1}^n\mathrm{Fit}_{k_i}(\phi_i)$, which is a coherent sheaf of ideals on $X^0_{\lambda-\min}$. 
        
        Assume Condition \ref{condition before blow-up for CDRS}, which is equivalent to the condition that $Z_{\lambda-\min}\not\subseteq \mathbb V(\mathcal I_{k_1,\cdots,k_n})$ by Lemma \ref{lemma of point-wise description of the nice locus of small unipotent stabilisers}. Let $Z\hookrightarrow X^0_{\lambda-\min}$ be the scheme theoretic intersection $Z_{\lambda-\min}\cap \mathbb V(\mathcal I_{k_1,\cdots,k_n})$. Let $C\hookrightarrow X^0_{\lambda-\min}$ be the scheme theoretic image of $U\times Z\to X^0_{\lambda-\min}$. Let $\overline C\hookrightarrow X$ be the scheme theoretic image of $C\hookrightarrow X^0_{\lambda-\min}\subseteq X$. Let $\pi:\tilde X\to X$ be the blow-up of $X$ along $\overline C\hookrightarrow X$. For any $x,y\in\mathbb N_+$, consider the ample line bundle $\tilde L:=\pi^*L^x\otimes\mathcal O(-yE)$ on $\tilde X$. There is an induced action and linearisation $\hat U \curvearrowright\big(\tilde X,\tilde L\big)$. Then similar constructions apply on $\tilde X^0_{\lambda-\min}:=\big(\tilde X\big)^0_{\lambda-\min}$
        \begin{equation}
            \begin{split}
                &\tilde{\mathcal K}_i:=\ker \Big(\Omega_{\tilde X^0_{\lambda-\min}/\Bbbk}\to (\mathfrak u_i)^*\otimes\mathcal O_{\tilde X^0_{\lambda-\min}}\Big)\\
                &\tilde\phi_i:\tilde{\mathcal K}_{i-1}\to (\mathfrak u_i/\mathfrak u_{i-1})^*\otimes\mathcal O_{\tilde X^0_{\lambda-\min}}\\
                &\mathrm{Fit}_k(\tilde\phi_i):=\mathrm{Fit}_k\big(\mathrm{coker}(\tilde\phi_i)\big). 
            \end{split}
        \end{equation}
        
        {
            \begin{theorem}\label{theorem of blow-up for CDRS}
                Let $X$ be a projective scheme over $\Bbbk$ and let $L$ be an ample line bundle on $X$. Let $\hat U \curvearrowright (X,L)$ and let $\lambda:\mathbb G_m\to \mathrm{Aut}(U)$ be a grading on $U$. Assume that Condition \ref{condition before blow-up for CDRS} holds. Let $\pi:\tilde X\to X$ be the blow-up described above with the ample line bundle $\tilde L$. Let $(k_1,\cdots,k_n)$, $\tilde\phi_i$ be as above. Then $\mathrm{coker}(\tilde\phi_i)$ is locally free of rank $k_i$ on $\tilde X^0_{\lambda-\min}$ for $1\leq i\leq n$. 
            \end{theorem}
            \begin{proof}
                The statement is local. Let $m\in\mathbb N_+$ be such that $H^0(X,L^m)$ generates $\bigoplus_{d\geq0}H^0(X,L^{md})$. Let $f\in H^0(X,L^m)_{\lambda=\max}$ be such that $Z_{\lambda-\min,f}\not\subseteq \mathbb V(\mathcal I_{k_1,\cdots,k_n})$. Let $I_f:=\mathcal I_{k_1,\cdots,k_n}(X_f)+\mathcal O(X_f)_{\lambda<0}\subseteq \mathcal O(X_f)$ and let $J_f\subseteq\mathcal O(X_f)$ be the kernel of $\mathcal O(X_f)\to \mathcal O(U)\otimes\big(\mathcal O(X_f)/I_f\big)$. Then $J_f$ is the ideal associated to $C_f\hookrightarrow X_f$. 
                
                A set of generators of $\mathrm{Fit}_{k_i}(\phi_i)(X_f)_{\lambda=0}$ as an ideal of $\mathcal O(X_f)_{\lambda=0}$ consists of $(r_i-k_i)\times(r_i-k_i)$-minors of the following matrices
                \begin{equation}
                    \begin{pmatrix}\langle \xi^{(i)}_1,\omega^{(i)}_1\rangle &\cdots&\langle \xi^{(i)}_1,\omega^{(i)}_{r_i-k_i}\rangle\\
                    \vdots&\vdots&\vdots\\
                    \langle \xi^{(i)}_{r_i},\omega^{(i)}_1\rangle &\cdots&\langle \xi^{(i)}_{r_i},\omega^{(i)}_{r_i-k_i}\rangle
                    \end{pmatrix},\quad \omega^{(i)}_1,\cdots\omega^{(i)}_{r_i-k_i}\in\mathcal K_{i-1}(X_f)_{\lambda=-w_i}. 
                \end{equation}
                Choose an $(r_i-k_i)\times(r_i-k_i)$-minor of the above matrix. We can change the basis of $\mathfrak u_{\lambda=w_i}$ and then assume this $(r_i-k_i)\times(r_i-k_i)$-minor has the form 
                \begin{equation}
                    a^{(i)}=\det\begin{pmatrix}\langle \xi^{(i)}_1,\omega^{(i)}_1\rangle &\cdots&\langle \xi^{(i)}_1,\omega^{(i)}_{r_i-k_i}\rangle\\
                    \vdots&\ddots&\vdots\\
                    \langle \xi^{(i)}_{r_i-k_i},\omega^{(i)}_1\rangle &\cdots&\langle \xi^{(i)}_{r_i-k_i},\omega^{(i)}_{r_i-k_i}\rangle
                    \end{pmatrix}. 
                \end{equation}
                If $\omega^{(i)}_{j}=\sum_\mu a^{(i)}_{j,\mu}\mathrm{d}b^{(i)}_{j,\mu}$, then let $f^{(i)}_j=\sum_\mu\big( a^{(i)}_{j,\mu}\big)_{\lambda=0}\big(b^{(i)}_{j,\mu}\big)_{\lambda=-w_i}\in\mathcal O(X_f)_{\lambda=-w_i}$. We then have that every $\mathrm{d}f^{(i)}_j\in\mathcal K_{i-1}(X_f)_{\lambda=-w_i}$ and 
                \begin{equation}
                    a^{(i)}=\det\begin{pmatrix}\xi^{(i)}_1.f^{(i)}_1&\cdots&\xi^{(i)}_1.f^{(i)}_{r_i-k_i}\\
                    \vdots&\ddots&\vdots\\
                    \xi^{(i)}_{r_i-k_i}.f^{(i)}_1&\cdots&\xi^{(i)}_{r_i-k_i}.f^{(i)}_{r_i-k_i}
                    \end{pmatrix}. 
                \end{equation}
                
                Let $a=a^{(1)}\cdots a^{(n)}$ with $a^{(i)}\in \mathrm{Fit}_{k_i}(\phi_i)_{\lambda=0}$ of the above form for $f^{(i)}_j\in\mathcal O(X_f)_{\lambda=-w_i}$. Then $a\in (J_f)_{\lambda=0}$ and such elements generate $(J_f)_{\lambda=0}$, i.e. $\tilde X^0_{\lambda-\min}\cap \tilde X_f$ is covered by such $\tilde X_{f,a}$.  
                
                The theorem is proved if we can show that on $\tilde X_{f,a}$, for $1\leq i\leq n$
                \begin{equation}
                    \mathrm{Fit}_{k_i-1}(\tilde\phi_i)(\tilde X_{f,a})=0,\quad \mathrm{Fit}_{k_i}(\tilde\phi_i)(\tilde X_{f,a})=\mathcal O(\tilde X_{f,a}). 
                \end{equation}
                
                There is a natural surjective morphism 
                \begin{equation}
                    \begin{tikzcd}
                        \mathcal O(\tilde X_{f,a})\otimes_{\Bbbk}\mathrm{d}\big(\frac{J_f}{a}\big)\ar[r]&\Omega_{\mathcal O(\tilde X_{f,a})/\Bbbk}\ar[r]&0. 
                    \end{tikzcd}
                \end{equation}
                
                Let $\tilde\omega^{(i)}_1,\cdots,\tilde\omega^{(i)}_{r_i-k_i+1}\in\tilde{\mathcal K}_{i-1}(X_f)$. By the above surjective morphism, we can write $\tilde\omega^{(i)}_j=\sum_\mu \frac{h^{(i)}_{j,\mu}}{a^N}\mathrm{d}\big(\frac{g^{(i)}_{j,\mu}}{a}\big)$ for $N\in\mathbb N_+$, $h^{(i)}_{j,\mu}\in(J_f)^N$ and $g^{(i)}_{j,\mu}\in J_f$. The fact that $\tilde\omega^{(i)}_j\in\tilde{\mathcal K}_{i-1}(\tilde X_{f,a})$ implies that there exists $M\in\mathbb N_+$ such that $\omega^{(i)}_j:=a^M\sum_\mu h^{(i)}_{j,\mu}\mathrm{d} g^{(i)}_{j,\mu}\in \mathcal K_{i-1}(X_f)$. 
                
                The matrix representing $\tilde\phi_i$ on $\big(\tilde\omega^{(i)}_1,\cdots,\tilde\omega^{(i)}_{r_i-k_i+1}\big)$ is
                \begin{equation}
                    \begin{split}
                        \tilde\phi_i&\big(\tilde\omega^{(i)}_1,\cdots,\tilde\omega^{(i)}_{r_i-k_i+1}\big)\\
                        &=\big(u^{(i)}_1\otimes1,\cdots,u^{(i)}_{r_i}\otimes 1\big)\begin{pmatrix}\sum_\mu\frac{h^{(i)}_{1,\mu}}{a^N}\frac{\xi^{(i)}_1.g^{(i)}_{1,\mu}}{a}&\cdots&\sum_\mu\frac{h^{(i)}_{r_i-k_i+1,\mu}}{a^N}\frac{\xi^{(i)}_1.g^{(i)}_{r_i-k_i+1,\mu}}{a}\\
                        \vdots&\vdots&\vdots\\
                        \sum_\mu\frac{h^{(i)}_{1,\mu}}{a^N}\frac{\xi^{(i)}_{r_i}.g^{(i)}_{1,\mu}}{a}&\cdots&\sum_\mu\frac{h^{(i)}_{r_i-k_i+1,\mu}}{a^N}\frac{\xi^{(i)}_{r_i}.g^{(i)}_{r_i-k_i+1,\mu}}{a}
                        \end{pmatrix}\\
                        &=\big(u^{(i)}_1\otimes1,\cdots,u^{(i)}_{r_i}\otimes 1\big)\begin{pmatrix}\frac{\langle\xi^{(i)}_1,\omega^{(i)}_1\rangle}{a^{M+N+1}}&\cdots&\frac{\langle\xi^{(i)}_1,\omega^{(i)}_{r_i-k_i+1}\rangle}{a^{M+N+1}}\\
                        \vdots&\vdots&\vdots\\
                        \frac{\langle\xi^{(i)}_{r_i},\omega^{(i)}_1\rangle}{a^{M+N+1}}&\cdots&\frac{\langle\xi^{(i)}_{r_i},\omega^{(i)}_{r_i-k_i+1}\rangle}{a^{M+N+1}}
                        \end{pmatrix}. 
                    \end{split}
                \end{equation}
                
                Let 
                \begin{equation}
                    \begin{split}
                        \tilde A:=&\begin{pmatrix}\frac{\langle\xi^{(i)}_1,\omega^{(i)}_1\rangle}{a^{M+N+1}}&\cdots&\frac{\langle\xi^{(i)}_1,\omega^{(i)}_{r_i-k_i+1}\rangle}{a^{M+N+1}}\\
                        \vdots&\vdots&\vdots\\
                        \frac{\langle\xi^{(i)}_{r_i},\omega^{(i)}_1\rangle}{a^{M+N+1}}&\cdots&\frac{\langle\xi^{(i)}_{r_i},\omega^{(i)}_{r_i-k_i+1}\rangle}{a^{M+N+1}}
                        \end{pmatrix}\\
                        A:=&\begin{pmatrix}\langle\xi^{(i)}_1,\omega^{(i)}_1\rangle&\cdots&\langle\xi^{(i)}_1,\omega^{(i)}_{r_i-k_i+1}\rangle\\\vdots&\vdots&\vdots\\\langle\xi^{(i)}_{r_i},\omega^{(i)}_1\rangle&\cdots&\langle\xi^{(i)}_{r_i},\omega^{(i)}_{r_i-k_i+1}\rangle\end{pmatrix}. 
                    \end{split}
                \end{equation}
                
                Let $\tilde m$ be an $(r_i-k_i+1)\times(r_i-k_i+1)$-minor of $\tilde A$. Then it is of the form 
                \begin{equation}
                    \tilde m=\frac{m}{a^{(r_i-k_i+1)(M+N+1)}}
                \end{equation}
                where $m\in (J_f)^{(r_i-k_i+1)(M+N+1)}$ is the $(r_i-k_i+1)\times(r_i-k_i+1)$-minor of $A$ with the same columns and rows as $\tilde m$. By the definition of $k_i$, any $(r_i-k_i+1)\times(r_i-k_i+1)$-minor of the matrix $A$ is zero. In particular $m=0$, so 
                \begin{equation}
                    \tilde m=\frac{m}{a^{(r_i-k_i+1)(M+N+1)}}=0
                \end{equation}
                which proves $\mathrm{Fit}_{k_i-1}(\tilde \phi_i)(\tilde X_{f,a})=0$. 
                
                We then prove $\mathrm{Fit}_{k_i}(\tilde \phi_i)(\tilde X_{f,a})=\mathcal O(\tilde X_{f,a})$. By Lemma \ref{lemma of existence of b^(i)_nu}, for $1\leq i\leq n$ and $1\leq \mu\leq r_i-k_i$, there exists $c^{(i)}_\mu\in(J_f)_{\lambda=-w_i}$ such that 
                \begin{equation}
                    \xi^{(i)}_{\mu}.c^{(i)}_\nu=w_i\delta_{\mu\nu}a,\quad 1\leq i\leq n,\;1\leq \mu\leq r_i-k_i. 
                \end{equation}
                Then $\mathrm{d}\big(\frac{c^{(i)}_\mu}{a}\big)\in\big(\Omega_{\mathcal O(\tilde X_{f,a})/\Bbbk}\big)_{\lambda=-w_i}\subseteq \tilde{\mathcal K}_{i-1}(\tilde X_{f,a})$. For each $1\leq i\leq n$ the matrix representing $\tilde\phi_i$ on $\Big(\mathrm{d}\big(\frac{c^{(i)}_1}{a}\big),\cdots,\mathrm{d}\big(\frac{c^{(i)}_{r_i-k_i}}{a}\big)\Big)$ is 
                \begin{equation}
                    \begin{split}
                        \tilde\phi_i&\Big(\mathrm{d}\big(\frac{c^{(i)}_1}{a}\big),\cdots,\mathrm{d}\big(\frac{c^{(i)}_{r_i-k_i}}{a}\big)\Big)\\
                        &=\big(u^{(i)}_1\otimes1,\cdots,u^{(i)}_{r_i}\otimes 1\big)\begin{pmatrix}\begin{matrix}w_i&\cdots&0\\\vdots&\ddots&\vdots\\0&\cdots&w_i\end{matrix}\\M'\end{pmatrix}. 
                    \end{split}
                \end{equation}
                Then the top $(r_i-k_i)\times(r_i-k_i)$-minor of the matrix is $w_i^{r_i-k_i}\in \mathrm{Fit}_{k_i}(\tilde\phi_i)(\tilde X_{f,a})$. Therefore $\mathrm{Fit}_{k_i}(\tilde\phi_i)(\tilde X_{f,a})=\mathcal O(\tilde X_{f,a})$. This completes the proof. 
            \end{proof}
            
            \begin{corollary}\label{corollary of consequence of blow-up}
                Let $X$ be a projective scheme over $\Bbbk$ and let $L$ be an ample line bundle on $X$. Let $\hat U\curvearrowright(X,L)$ and let $\lambda:\mathbb G_m\to \mathrm{Aut}(U)$ be a grading 1PS. Assume that Condition \ref{condition before blow-up for CDRS} holds. Then there exists a $\hat U$-equivariant blow-up $\pi:\tilde X\to X$ and the action $\hat U\curvearrowright\tilde X$ lifts to an ample line bundle $\tilde L$, such that: 
                \begin{itemize}
                    \item $\tilde X^0_{\lambda-\min}\to \tilde X^0_{\lambda-\min}/U$ is a quasi-projective universal geometric quotient by $U$; 
                    \item $\tilde X^0_{\lambda-\min}\setminus U\tilde Z_{\lambda}\to \big(\tilde X^0_{\lambda-\min}\setminus U\tilde Z_{\lambda}\big)\big/\hat U$ is a projective universal geometric quotient by $\hat U$. 
                \end{itemize}
            \end{corollary}
            \begin{proof}
                This is obvious from Theorem \ref{theorem of blow-up for CDRS}, Theorem \ref{theorem of universal geometric quotient by Un with proj CDRS} and Theorem \ref{theorem of U-hat quotient with CDRS}. 
            \end{proof}
        }
    }
    
    \subsection{Appendix}
    {
        \subsubsection{Associative algebra and complete brackets}\label{appendix of complete brackets}
        {
            Let $A$ be an associative algebra. Let $m\in\mathbb N_+$. Define the \emph{complete bracket of $m$ elements} by 
            \begin{equation}
                \begin{split}
                    [-]\!]_m:&\prod_{i=1}^mA\to A\\
                    [a_1,\cdots,a_m]\!]_m:=&\mathrm{ad}_{a_1}\cdots\mathrm{ad}_{a_{m-1}}(a_m)= \Big[a_1,\big[a_2,\cdots[a_{m-1},a_m]\cdots\big]\Big]. 
                \end{split}
            \end{equation}
            When $m=1$, the map $[-]\!]_1:A\to A$ is the identity map. 
            
            Let $\mathfrak g$ be a Lie algebra over $\Bbbk$ and let $\mathrm{U}(\mathfrak g)$ be its universal enveloping algebra. Let $\big(\xi_1,\cdots,\xi_r\big)$ be an ordered basis of $\mathfrak g$. Then $\big\{\xi_1^{p_1}\cdots\xi_r^{p_r}:(p_1,\cdots,p_r)\in\mathbb N^r\big\}$ is a basis of $\mathrm{U}(\mathfrak g)$ by the Poincar\'e-Birkhoff-Witt theorem. We can define the \emph{complete bracket with respect to the basis $\big(\xi_1,\cdots,\xi_n\big)$} by 
            \begin{equation}
                \begin{split}
                    [-]\!]:\mathrm{U}(\mathfrak g)&\to\mathfrak g\\
                    [\xi^p]\!]:=&\Big[\underbrace{\xi_1,\cdots,\xi_1}_{p_1},\cdots,\underbrace{\xi_r,\cdots,\xi_r}_{p_r}\Big]\!\Big]_{p_1+\cdots+p_r}. 
                \end{split}
            \end{equation}
            
            Moreover if there is a representation of $\mathbb G_m$ on $\mathfrak g$, which induces a representation of $\mathbb G_m$ on $\mathrm{U}(\mathfrak g)$, and each $\xi_i$ in the basis $(\xi_1,\cdots,\xi_r)$ is a weight vector, then $[-]\!]:\mathrm{U}(\mathfrak g)\to \mathfrak g$ preserves the weights, i.e. it is a morphism of $\mathbb G_m$-representations. 
    
            {
                \begin{lemma}\label{lemma of identity of weighted complete bracket}
                    Let $n\in\mathbb N_+$ and let $\Bbbk\langle y_1,\cdots,y_n\rangle$ be the free associative algebra over $\Bbbk$ generated by $n$ variables. Let $w_1,\cdots,w_n\in\Bbbk$ be scalars. Let $k_1,\cdots,k_n\in\mathbb N$. Then 
                    \begin{equation}
                        \Big(\sum_{i=1}^nk_iw_i\Big)y^k=\sum_{0<s\leq k}\binom{k}{s}w_{\max\{i:s_i\ne0\}}[y^s]\!]y^{k-s}
                    \end{equation}
                    where 
                    \begin{equation}
                        \begin{split}
                            s:=&(s_1,\cdots,s_n)\ne(0,\cdots,0),\quad s_i\leq k_i\;\textrm{for }i=1,\cdots,n\\
                            \binom{k}{s}:=&\binom{k_1}{s_1}\cdots\binom{k_n}{s_n}\\
                            y^s:=&y_1^{s_1}\cdots y_n^{s_n}\\
                            [y^s]]:=&\Big[\underbrace{y_1,\cdots,y_1}_{s_1},\cdots,\underbrace{y_n,\cdots,y_n}_{s_n}\Big]\!\Big]_{s_1+\cdots+s_n}
                        \end{split}
                    \end{equation}
                \end{lemma}
            }
            
            {
                \begin{lemma}\label{lemma of commutator identity with complete bracket}
                    Let $n\in\mathbb N$ and let $\Bbbk\langle x_1,\cdots,x_n,y\rangle$ be the free associative algebra over $\Bbbk$ generated by $n+1$ variables. Let $k:=(k_1,\cdots,k_n)\in\mathbb N^n$. Then 
                    \begin{equation}
                        x^ky=\sum_{0\leq s\leq k}\binom{k}{s}[x^{k-s}y]\!]x^s
                    \end{equation}
                    where 
                    \begin{equation}
                        [x^{k-s}y]\!]:=\Big[\underbrace{x_1,\cdots,x_1}_{k_1-s_1},\cdots,\underbrace{x_n,\cdots,x_n}_{k_n-s_n},y\Big]\!\Big]_{1+\sum_{i=1}^n(k_i-s_i)}. 
                    \end{equation}
                \end{lemma}
            }
        }
    
        \subsubsection{Unipotent action on affine schemes}
        {
            Let $U$ be an affine algebraic unipotent group over $\Bbbk$. Let $A$ be a commutative $\Bbbk$-algebra and let $U$ acts on $A$ dually by $\varphi:A\to \mathcal O(U)\otimes A$, where the tensor product is over $\Bbbk$. Identify $\mathcal O(U)\cong \Bbbk[u_1,\cdots,u_n]$ and let $\mu^*:\Bbbk[u]\to \Bbbk[u]\otimes\Bbbk[u]$ be the co-multiplication. Let $\mathrm{U}(\mathfrak u)$ be the universal enveloping algebra of $\mathfrak u:=\mathrm{Lie}(U)$. The aim is to prove the following. 
    
            {
                \begin{proposition}\label{proposition of dual action of unipotent groups}
                    For $f\in A$, if 
                    \begin{equation}
                        A\to \Bbbk[u]\otimes A,\quad f\mapsto \sum_{\alpha\in\mathbb N^n}u^\alpha\otimes f_\alpha
                    \end{equation}
                    then $f_\alpha\in\mathrm{U}(\mathfrak u).f$ for all $\alpha$. 
                \end{proposition}
            }

            {
                Before the proof we state a useful corollary. 
                \begin{corollary}\label{corollary of the ideal cutting out U-sweep of Z consists functions whose UEA derivations in the ideal of Z}
                    Let $U$ be a unipotent group over $\Bbbk$ and let $A$ be a $\Bbbk$-algebra of finite type. Let $\varphi:A\to \mathcal O(U)\otimes A$ be a dual action of $U$ on $A$. Let $I\subseteq A$ be an ideal. Then the following are equivalent for $a\in A$: 
                    \begin{itemize}
                        \item[(1)] $X.a\in I$ for all $X\in\mathrm{U}(\mathfrak u)$; 
                        \item[(2)] the map $A\to\mathcal O(U)\otimes A\to \mathcal O(U)\otimes\big(A/I\big)$ sends $a\in A$ to zero. 
                    \end{itemize}
                \end{corollary}
                \begin{proof}
                    Identify $\mathcal O(U)\cong \Bbbk[u_1,\cdots,u_n]$. For $a\in A$ and $\alpha\in\mathbb N^n$, let $a_\alpha\in A$ be such that 
                    \begin{equation}
                        A\to \mathcal O(U)\otimes A,\quad a\mapsto \sum_{\alpha\in\mathbb N^n}u^\alpha\otimes a_\alpha. 
                    \end{equation}
                    
                    Assume $(1)$. We have $\mathrm{U}(\mathfrak u).a\subseteq I$, and then $a_\alpha\in I$ by Proposition \ref{proposition of dual action of unipotent groups}, which proves $(2)$. 
                    
                    Assume $(2)$. Let $q:A\to A/I$ be the quotient ring map. Let $J\subseteq A$ be the kernel of $A\to \mathcal O(U)\otimes\big(A/I\big)$. To prove $(1)$, it suffices to show: 
                    \begin{itemize}
                        \item $J\subseteq I$; 
                        \item $J\subseteq A$ is $U$-invariant. 
                    \end{itemize}
                    Assume they are proved. For $a\in J$, since $J$ is $U$-invariant, we have $X.a\in J$ for all $X\in\mathrm{U}(\mathfrak u)$. Since $J\subseteq I$, we have $X.a\in I$, which is $(1)$. 
                    
                    We prove that $J\subseteq I$. Consider the diagram 
                    \begin{equation}
                        \begin{tikzcd}
                            A\ar[r,"\varphi"]\ar[rd,"1"]&\mathcal O(U)\otimes A\ar[r,"1\otimes q"]\ar[d,"e^*\otimes 1"]&\mathcal O(U)\otimes\big(A/I\big)\ar[d,"e^*\otimes 1"]\\
                            &A\ar[r,"q"]&A/I
                        \end{tikzcd}
                    \end{equation}
                    where $e^*:\mathcal O(U)\to\Bbbk$ is the co-unit map. Since $J=\ker\big((1\otimes q)\circ\varphi\big)$, we have $q(J)=q\circ 1(J)=(e^*\otimes 1)\circ(1\otimes q)\circ\varphi(J)=0$, i.e. $J\subseteq I$. 
    
                    We prove that $J\subseteq A$ is $U$-invariant. Consider the following diagram 
                    \begin{equation}
                        \begin{tikzcd}
                            A\ar[r,"\varphi"]\ar[d,"\varphi"]&\mathcal O(U)\otimes A\ar[d,"1\otimes \varphi"]\\
                            \mathcal O(U)\otimes A\ar[r,"\mu^*\otimes1"]\ar[d,"1\otimes q"]&\mathcal O(U)\otimes\mathcal O(U)\otimes A\ar[d,"1\otimes 1\otimes q"]\\
                            \mathcal O(U)\otimes\big(A/I\big)\ar[r,"\mu^*\otimes 1"]&\mathcal O(U)\otimes\mathcal O(U)\otimes\big(A/I\big). 
                        \end{tikzcd}
                    \end{equation}
                    Since $J=\ker \big((1\otimes q)\circ\varphi\big)$, by the diagram we have $\varphi(J)\subseteq \ker\big((1\otimes1\otimes q)\circ(1\otimes\varphi)\big)$. The composition of the right vertical arrows, $(1\otimes1\otimes q)\circ(1\otimes\varphi)$, is the base change of $(1\otimes q)\circ \varphi$ along $\Bbbk\to\mathcal O(U)$. By the flatness of $\Bbbk\to \mathcal O(U)$, we have $\mathcal O(U)\otimes J=\ker\big((1\otimes1\otimes q)\circ(1\otimes\varphi)\big)$. So $\varphi(J)\subseteq\mathcal O(U)\otimes J$, which proves that $J$ is $U$-invariant. 
                \end{proof}
            }
            
            The action of $\mathfrak u$ on $A$ is defined as follows 
            \begin{equation}
                \begin{tikzcd}
                    A\ar[r,"\varphi"]\ar[rrr,bend right=30,"\xi"]&\Bbbk[u]\otimes A\ar[r,"\psi"]&\mathfrak m/\mathfrak m^2\otimes A\ar[r,"(-\xi)\otimes 1"]&A
                \end{tikzcd}
            \end{equation}
            where $\mathfrak m=(u_1,\cdots,u_n)\subseteq\Bbbk[u]$ is the maximal ideal associated to the identity, and $\psi$ is the linear map 
            \begin{equation}
                \psi:\Bbbk[u]\to \mathfrak m/\mathfrak m^2,\quad p(u)\mapsto \overline{p(u)-p(0)}. 
            \end{equation}
            The action of $\mathrm{U}(\mathfrak u)$ on $A$ is induced from that of $\mathfrak u\curvearrowright A$. 
            
            There exists $c^\alpha_{\beta,\gamma}\in\Bbbk$ for $\alpha,\beta,\gamma\in\mathbb N^n$ such that 
            \begin{equation}
                \mu^*:\Bbbk[u]\to\Bbbk[u]\otimes\Bbbk[u],\quad u^\alpha\mapsto\sum_{\beta,\gamma}c^\alpha_{\beta,\gamma}u^\beta\otimes u^\gamma. 
            \end{equation}
            They satisfy the following properties: 
            \begin{itemize}
                \item Since $\mu^*$ respects multiplication, we have 
                \begin{equation}
                    c^{\alpha_1+\alpha_2}_{\beta,\gamma}=\sum_{\substack{\beta_1+\beta_2=\beta\\\gamma_1+\gamma_2=\gamma}}c^{\alpha_1}_{\beta_1,\gamma_1}c^{\alpha_2}_{\beta_2,\gamma_2}; 
                \end{equation}
                \item The existence of identity in $U$ implies 
                \begin{equation}
                    c^\alpha_{0,\gamma}=\delta^\alpha_\gamma,\quad c^\alpha_{\beta,0}=\delta^\alpha_\beta. 
                \end{equation}
            \end{itemize}
            
            {
                \begin{lemma}\label{lemma of coefficient nonzero imply |alpha|<=|beta|+|gamma|}
                    If $c^\alpha_{\beta,\gamma}\ne0$, then $|\alpha|\leq |\beta|+|\gamma|$. 
                \end{lemma}
                \begin{proof}
                    We prove $c^\alpha_{\beta,\gamma}=0$ if $|\alpha|>|\beta|+|\gamma|$ by induction on $|\alpha|$. When $|\alpha|=1$, we have $|\beta|+|\gamma|=0$, i.e. $\beta=\gamma=0$. Then $c^\alpha_{\beta,\gamma}=c^\alpha_{0,0}=\delta^\alpha_0=0$. Assume $c^\alpha_{\beta,\gamma}=0$ if $k\geq|\alpha|>|\beta|+|\gamma|$ for some $k\in\mathbb N_+$. If $k+1=|\alpha|>|\beta|+|\gamma|$, write $\alpha=\alpha'+e_i$ for some $|\alpha'|=k$ and $1\leq i\leq n$. Then 
                    \begin{equation}
                        \begin{split}
                            c^\alpha_{\beta,\gamma}=&c^{\alpha'+e_i}_{\beta,\gamma}\\
                            =&\sum_{\substack{\beta'+\beta''=\beta\\\gamma'+\gamma''=\gamma}}c^{\alpha'}_{\beta',\gamma'}c^{e_i}_{\beta'',\gamma''}. 
                        \end{split}
                    \end{equation}
                    If there exists a nonzero summand 
                    \begin{equation}
                        c^{\alpha'}_{\beta',\gamma'}\ne0,\quad c^{e_i}_{\beta'',\gamma''}\ne0
                    \end{equation}
                    then by the induction hypothesis for $|\alpha'|=k$ and $|e_i|=1$ 
                    \begin{equation}
                        |\alpha'|\leq|\beta'|+|\gamma'|,\quad 1\leq|\beta''|+|\gamma''|
                    \end{equation}
                    which implies $|\alpha|=|\alpha'|+|e_j|\leq |\beta|+|\gamma|$, a contradiction. Then $c^\alpha_{\beta,\gamma}=0$ as a sum of zeros. 
                \end{proof}
            }
            
            {
                \begin{lemma}\label{lemma of equivalent conditions of coefficient when gamma = e_j}
                    If $|\alpha|=|\beta|+1$ and $1\leq j\leq n$, then the following are equivalent: 
                    \begin{itemize}
                        \item[(1)] $c^\alpha_{\beta,e_j}=1+\max\{k\in\mathbb N:ke_j\leq \beta\}$; 
                        \item[(2)] $c^\alpha_{\beta,e_j}\ne0$; 
                        \item[(3)] $\alpha=\beta+e_j$. 
                    \end{itemize}
                \end{lemma}
                \begin{proof}
                    We prove the equivalence by induction on $|\alpha|$. When $|\alpha|=1$, we have $\beta=0$ and then $c^\alpha_{0,e_j}=\delta^\alpha_{e_j}$, which proves the equivalence. 
        
                    Assume for some $k\in\mathbb N_+$ the equivalence is proved for $|\alpha|\leq k$. Let $|\alpha|=k+1$. Obviously $(1)\implies(2)$. 
                    
                    For $(2)\implies(3)$, assume $c^\alpha_{\beta,e_j}\ne0$. Write $\alpha=\alpha'+e_i$ for some $|\alpha'|=k$ and $1\leq i\leq n$. Then 
                    \begin{equation}
                        \begin{split}
                            0\ne c^\alpha_{\beta,e_j}=&\sum_{\beta'+\beta''=\beta}c^{\alpha'}_{\beta',e_j}c^{e_i}_{\beta'',0}+\sum_{\beta'+\beta''=\beta}c^{\alpha'}_{\beta',0}c^{e_i}_{\beta'',e_j}\\
                            =&\sum_{\beta'+\beta''=\beta}c^{\alpha'}_{\beta',e_j}\delta^{e_i}_{\beta''}+\sum_{\beta'+\beta''=\beta}\delta^{\alpha'}_{\beta'}c^{e_i}_{\beta'',e_j}. 
                        \end{split}
                    \end{equation}
                    If there exist $\beta'+\beta''=\beta$ such that $c^{\alpha'}_{\beta',e_j}\delta^{e_i}_{\beta''}\ne0$, and then 
                    \begin{equation}
                        \beta''=e_i,\quad |\beta'|=|\beta|-|\beta''|=k-1. 
                    \end{equation}
                    So $k=|\alpha'|=|\beta'|+1$, which implies $\alpha'=\beta'+e_j$ by the induction hypothesis. So $\alpha=\alpha'+e_i=\beta'+e_j+\beta''=\beta+e_j$, which is $(3)$. 
        
                    Else there exist $\beta'+\beta''=\beta$ such that $\delta^{\alpha'}_{\beta'}c^{e_i}_{\beta'',e_j}\ne0$, and then $\alpha'=\beta'$. In particular $|\beta''|=|\beta|-|\beta'|=0$, i.e. $\beta''=0$. Then $c^{e_i}_{\beta'',e_j}=\delta^{e_i}_{e_j}\ne0$, i.e. $e_i=e_j$. We have $\alpha=\alpha'+e_i=\beta'+e_j=\beta+e_j$, which is $(3)$. 
                    
                    For $(3)\implies (1)$, assume $\alpha=\beta+e_j$. Then 
                    \begin{equation}
                        \begin{split}
                            c^\alpha_{\beta,e_j}=&c^{\beta+e_j}_{\beta,e_j}\\
                            =&\sum_{\beta'+\beta''=\beta}c^\beta_{\beta',e_j}c^{e_j}_{\beta'',0}+\sum_{\beta'+\beta''=\beta}c^\beta_{\beta',0}c^{e_j}_{\beta'',e_j}\\
                            =&\sum_{\beta'+\beta''=\beta}c^\beta_{\beta',e_j}\delta^{e_j}_{\beta''}+\sum_{\beta'+\beta''=\beta}\delta^\beta_{\beta'}c^{e_j}_{\beta'',e_j}\\
                            =&\sum_{\beta'+\beta''=\beta}c^\beta_{\beta',e_j}\delta^{e_j}_{\beta''}+1. 
                        \end{split}
                    \end{equation}
                    If $e_j\leq \beta$, then 
                    \begin{equation}
                        c^\alpha_{\beta,e_j}=c^\beta_{\beta-e_j,e_j}+1
                    \end{equation}
                    and $c^\beta_{\beta-e_j,e_j}=1+\max\{k\in\mathbb N:ke_j\leq \beta-e_j\}=\max\{k\in\mathbb N:ke_j\leq \beta\}$ by the induction hypothesis, which proves $(1)$. Else $e_j\not\leq \beta$, then $c^\alpha_{\beta,e_j}=1$, which is $(1)$. 
                \end{proof}
            }
        
            {
                We can now prove Proposition \ref{proposition of dual action of unipotent groups}. 
                \begin{proof}[Proof of Proposition \ref{proposition of dual action of unipotent groups}]
                    We prove this by induction on $|\alpha|$. When $\alpha=0$, we have $f_0=f$. When $|\alpha|=1$, we have $f_{e_i}=-\xi_i.f\in\mathrm{U}(\mathfrak u).f$. 
                    
                    Assume for some $k\in\mathbb N_+$, we have $f_\alpha\in\mathrm{U}(\mathfrak u).f$ when $|\alpha|\leq k$. Let $|\alpha|=k+1$. 
                    
                    Consider the ring map 
                    \begin{equation}
                        \varphi_{k+1}:=A\to \underbrace{\Bbbk[u]\otimes\cdots\otimes\Bbbk[u]}_{k+1}\otimes A
                    \end{equation}
                    which maps $f\in A$ to 
                    \begin{equation}
                        \varphi_{k+1}(f)=\sum_{\beta_1,\cdots,\beta_{k+1}}\sum_{\gamma_1,\cdots,\gamma_{k+1}}c^{\beta_1}_{0,\gamma_1}c^{\beta_2}_{\beta_1,\gamma_2}c^{\beta_3}_{\beta_2,\gamma_3}\cdots c^{\beta_{k+1}}_{\beta_k,\gamma_{k+1}}u^{\gamma_1}\otimes\cdots\otimes u^{\gamma_{k+1}}\otimes f_{\beta_{k+1}}. 
                    \end{equation}
                    
                    Write $\alpha=e_{i_1}+\cdots+e_{i_{k+1}}$ for $1\leq i_1,\cdots,i_{k+1}\leq n$. Consider the map $X$
                    \begin{equation}
                        \begin{tikzcd}
                            \Bbbk[u]\otimes\cdots\otimes\Bbbk[u]\otimes A\ar[r,equal]\ar[dd,"X"]&\Bbbk[u]\otimes\cdots\otimes\Bbbk[u]\otimes A\ar[d,"\psi\otimes\cdots\otimes\psi\otimes 1"]\\
                            &\mathfrak m/\mathfrak m^2\otimes\cdots\otimes \mathfrak m/\mathfrak m^2\otimes A\ar[d,"(-\xi_{i_1})\otimes\cdots\otimes(-\xi_{i_{k+1}})\otimes 1"]\\
                            A\ar[r,equal]&A
                        \end{tikzcd}
                    \end{equation}
                    which maps $\varphi_{k+1}(f)$ to $\xi_{i_{k+1}}\cdots\xi_{i_1}.f\in\mathrm{U}(\mathfrak u).f$ by the definition of $\mathrm{U}(\mathfrak u)\curvearrowright A$. We can evaluate $X\circ\varphi_{k+1}(f)$ using the expression of $\varphi_{k+1}(f)$ 
                    \begin{equation}
                        \begin{split}
                            &X\circ\varphi_{k+1}(f)\\
                            =&X\bigg(\sum_{\beta_1,\cdots,\beta_{k+1}}\sum_{\gamma_1,\cdots,\gamma_{k+1}}c^{\beta_1}_{0,\gamma_1}c^{\beta_2}_{\beta_1,\gamma_2}\cdots c^{\beta_{k+1}}_{\beta_k,\gamma_{k+1}}u^{\gamma_1}\otimes\cdots\otimes u^{\gamma_{k+1}}\otimes f_{\beta_{k+1}}\bigg)\\
                            =&\sum_{\beta,\gamma}\Big(\prod_{j=1}^{k+1}c^{\beta_j}_{\beta_{j-1},\gamma_j}\Big)\Big(\prod_{j=1}^{k+1}\big\langle -\xi_{i_j},\psi(u^{\gamma_j})\big\rangle\Big) f_{\beta_{k+1}}
                        \end{split}
                    \end{equation}
                    where $\beta_0=0$. We have 
                    \begin{equation}
                        \big\langle -\xi_{i_j},\psi(u^{\gamma_j})\big\rangle=-\delta^{\gamma_j}_{e_{i_j}},\quad j=1,\cdots,k+1
                    \end{equation}
                    so 
                    \begin{equation}
                        \sum_\beta \Big(\prod_{j=1}^{k+1}c^{\beta_j}_{\beta_{j-1},e_{i_j}}\Big)f_{\beta_{k+1}}=(-1)^{k+1}X\circ\varphi_{k+1}(f)\in\mathrm{U}(\mathfrak u).f. 
                    \end{equation}
                    
                    If $\beta_1,\cdots,\beta_{k+1}$ are such that $\prod_{j=1}^{k+1}c^{\beta_j}_{\beta_{j-1},e_{i_j}}\ne0$, then 
                    \begin{equation}
                        \beta_1=e_{i_1},\quad c^{\beta_j}_{\beta_{j-1},e_{i_i}}\ne0\;\textrm{for }j\geq 2. 
                    \end{equation}
                    
                    By Lemma \ref{lemma of coefficient nonzero imply |alpha|<=|beta|+|gamma|}, we have $|\beta_j|\leq |\beta_{j-1}|+1$ for $j\geq 2$. In particular $|\beta_{k+1}|\leq k+1$, and $|\beta_{k+1}|=k+1$ implies $|\beta_j|=j$ for all $j=1,\cdots,k+1$. By the induction hypothesis $f_{\beta_{k+1}}\in\mathrm{U}(\mathfrak u).f$ if $|\beta_{k+1}|\leq k$. Then 
                    \begin{equation}
                        \sum_{\substack{\beta_1,\cdots,\beta_{k+1}\\|\beta_j|=j\;\textrm{for all }j}}\Big(\prod_{j=1}^{k+1}c^{\beta_j}_{\beta_{j-1},e_{i_j}}\Big)f_{\beta_{k+1}}\in\mathrm{U}(\mathfrak u).f
                    \end{equation}
                    
                    By Lemma \ref{lemma of equivalent conditions of coefficient when gamma = e_j}, when $|\beta_j|=j$ for all $j$, we have $c^{\beta_j}_{\beta_{j-1},e_{i_j}}\ne0$ if and only if $\beta_j=\beta_{j-1}+e_{i_j}$. Therefore only one summand above is non-zero
                    \begin{equation}
                        \Big(\prod_{j=1}^{k+1}c^{\sigma_j}_{\sigma_{j-1},e_{i_j}}\Big) f_{\sigma_{k+1}}\in\mathrm{U}(\mathfrak u).f
                    \end{equation}
                    where $\sigma_j=e_{i_1}+\cdots+e_{i_j}$. The coefficient $\prod_{j=1}^{k+1}c^{\sigma_j}_{\sigma_{j-1},e_{i_j}}$ is nonzero since $c^{\sigma_j}_{\sigma_{j-1},e_{i_j}}\ne0$ by Lemma \ref{lemma of equivalent conditions of coefficient when gamma = e_j}. Notice that $\sigma_{k+1}=e_{i_1}+\cdots+e_{i_{k+1}}=\alpha$. We have 
                    \begin{equation}
                        f_\alpha\in\mathrm{U}(\mathfrak u).f
                    \end{equation}
                    which completes the induction. 
                \end{proof}
            }
            
        }
    }
}

\addcontentsline{toc}{section}{Bibliography}

\printbibliography
\end{document}